\documentclass[10pt,leqno]{amsart}

\usepackage{latexsym,amssymb,amsthm,amsmath,amscd}

\newtheorem{Theorem}{Theorem}[section]
\newtheorem{Lemma}{Lemma}[section]
\newtheorem{Proposition}{Proposition}[section]
\newtheorem{Corollary}{Corollary}[section]
\newtheorem{Definition}{Definition}[section]
\newtheorem{Example}{Example}[section]

\newtheorem{Remark}{Remark}

\def\({\left (}
\def\){\right )}
\def\<{\left<}
\def\> { \right>}
\def\e{\eqref}

\def\<{\left<}
\def\>{\right>}

\def\d{\delta }
\def\v{\varphi }
\def\g{\gamma }

\begin{document}

\title[]{A comprehensive survey on parallel submanifolds in Riemannian and pseudo-Riemannian manifolds}

\author[B.-Y. Chen]{ Bang-Yen Chen}
\address{B.-Y. Chen: Department of Mathematics, Michigan State University, East Lansing, Michigan 48824--1027, U.S.A.}
\email{chenb@msu.edu}

\begin{abstract}{A submanifold of a Riemannian manifold is called a parallel submanifold if its second fundamental form is parallel with respect to the van der Waerden-Bortolotti connection.  From submanifold point of view, parallel submanifolds are the simplest Riemannian submanifolds next to totally geodesic ones. Parallel submanifolds form an important class of Riemannian submanifolds since extrinsic invariants of a parallel submanifold do not vary from point to point.
In this paper we provide a  comprehensive survey on this important class of submanifolds.}
\end{abstract}


\maketitle

\noindent {\bf Table of contents}
\vskip.1in

\noindent Section 1. Introduction \dotfill 3

\noindent Section  2. Preliminaries \dotfill 4

2.1. Basic definitions, formulas and equations \dotfill 4

 2.2. Indefinite real space forms \dotfill 5

 2.3.  Gauss image  \dotfill 6

\noindent Section 3. Some general properties of parallel submanifolds \dotfill 7

\noindent Section 4. Parallel submanifolds of Euclidean spaces \dotfill 7

4.1. Gauss map and parallel submanifolds \dotfill 7

4.2. Normal sections and parallel submanifolds \dotfill 8

4.3. Symmetric submanifolds and parallel submanifolds\dotfill 9

4.4. Extrinsic $k$-symmetric submanifolds as $\nabla^c$-parallel submanifolds\dotfill 9

\noindent Section  5. Symmetric $R$-spaces and parallel submanifolds of real space forms \dotfill 10

5.1. Symmetric $R$-spaces and Borel subgroups \dotfill 10

5.2. Classification of symmetric $R$-spaces \dotfill 10

5.3. Ferus' theorem \dotfill 11

5.4. Parallel submanifolds in spheres \dotfill 11

5.5. Parallel submanifolds in hyperbolic spaces \dotfill 11

\noindent Section 6. Parallel Kaehler submanifolds\dotfill 11

6.1. Segre and Veronese maps\dotfill 12

6.2. Classification of parallel Kaehler submanifolds of $CP^m$ and $CH^m$ \dotfill 12

6.3. Parallel Kaehler submanifolds of Hermitian symmetric spaces \dotfill 13

6.4. Parallel Kaehler manifolds in complex Grassmannian manifolds  \dotfill 13

\noindent Section 7. Parallel totally real submanifolds\dotfill 14

7.1. Basics on totally real submanifolds \dotfill 14

7.2.  Parallel Lagrangian submanifolds in $CP^m$ \dotfill 14

7.3. Parallel  surfaces of  $CP^2$ and $CH^2$ \dotfill 15

7.4. Parallel totally real submanifolds in nearly Kaehler $S^6$ \dotfill 16

\noindent Section 8.  Parallel slant submanifolds of complex space forms \dotfill 16

8.1. Basics on slant submanifolds  \dotfill 16

8.2. Classification of parallel slant submanifolds  \dotfill 17

\noindent Section 9. Parallel submanifolds of quaternionic space forms and Cayley plane \dotfill 17

9.1. Parallel submanifolds of quaternionic space forms \dotfill 17

9.2. Parallel submanifolds of the Cayley plane\dotfill 18

\noindent Section 10. Parallel spatial submanifolds in pseudo-Euclidean spaces\dotfill 18

10.1. Marginally trapped surfaces \dotfill 18

10.2. Classification of parallel spatial surfaces in $\mathbb E^m_s$\dotfill 18

10.3. Special case: parallel spatial surfaces in $\mathbb E^3_1$ \dotfill 19

\noindent Section  11. Parallel spatial surfaces in $S^m_s$ \dotfill 19

11.1. Classification of parallel spatial surfaces in $S^m_s$ \dotfill 19

11.2. Special case: parallel spatial surfaces in $S^3_1$ \dotfill 21

\noindent Section 12. Parallel spatial surfaces in $H^m_s$ \dotfill 21

12.1. Classification of parallel spatial surfaces in $H^m_s$ \dotfill 21

12.2. A parallel spatial surfaces in $H^4_2$ \dotfill 23

12.3. Special case: parallel surfaces in $H^3_1$ \dotfill 24

\noindent Section 13. Parallel Lorentz surfaces in pseudo-Euclidean spaces \dotfill 24

13.1. Classification of parallel Lorentzian surfaces in $\mathbb E^m_s$ \dotfill 24

13.2. Classification of parallel Lorentzian surfaces in $E^3_1$ \dotfill 26

\noindent Section 14. Parallel surfaces in a light cone ${\mathcal LC}$ \dotfill  26

14.1. Light cones in general relativity \dotfill 26

14.2. Parallel surfaces in ${\mathcal LC}^3_1\subset \mathbb E^4_1$ \dotfill 27

14.3. Parallel surfaces in ${\mathcal LC}^3_2\subset \mathbb E^4_2$\dotfill 27

\noindent Section  15. Parallel surfaces in de Sitter space-time $S^4_1$\dotfill 27

15.1. Classification of parallel spatial surfaces in de Sitter space-time $S^4_1$
\dotfill 28

15.2. Classification of parallel Lorentzian surfaces in de Sitter space-time $S^4_1$ \dotfill 29

\noindent Section 16. Parallel surfaces in anti de Sitter space-time $H^4_1$ \dotfill 29

16.1. Classification of parallel spatial surfaces in $H^4_1$ \dotfill 29

16.2. Classification of parallel Lorentzian surfaces in anti de Sitter space-time $H^4_1$\dotfill 30

16.3. Special case: parallel Lorentzian surfaces in $H^3_1$ \dotfill 31

\noindent Section  17. Parallel spatial surfaces in $S^4_2$\dotfill 31

17.1. Four-dimensional manifolds with neutral metrics \dotfill 31

17.2. Classification of parallel Lorentzian surfaces in $S^4_2$ \dotfill 32

17.3. Classification of parallel Lorentzian surfaces in $H^4_2$ \dotfill 33

\noindent Section  18. Parallel spatial surfaces in  $S^4_3$ and in $H^4_3$ \dotfill 34

18.1. Classification of parallel spatial surfaces in $S^4_3$ \dotfill 34

18.2. Classification of parallel spatial surfaces in $H^4_3$ \dotfill 34

\noindent Section  19. Parallel Lorentzian surfaces in ${\mathbb C}^n$, $CP^2_1$ and  $CH^2_1$ \dotfill 35

19.1. Hopf fibration \dotfill 35

19.2. Classification of parallel spatial surfaces in $\mathbb C^2_1$ \dotfill 36

19.3. Classification of parallel Lorentzian surface in $CP^2_1$  \dotfill 36

19.4. Classification of parallel Lorentzian surface in $CH^2_1$  \dotfill 38

\noindent Section 20. Parallel Lorentz surfaces in $I\times_f R^n(c)$ \dotfill 38

20.1. Basics on Robertson-Walker space-times \dotfill 38

20.2. Parallel submanifolds of Robertson-Walker space-times \dotfill 39

\noindent Section 21. Thurston's eight 3-dimensional model geometries \dotfill 39

\noindent Section 22. Parallel surfaces in three-dimensional Lie groups\dotfill 40

22.1. Milnor's classification of 3-dimensional unimodular Lie groups \dotfill 40

22.2. Parallel  surfaces in  the motion group $E(1,1)$ \dotfill 41

22.3. Parallel surfaces in the Heisenberg group $Sol_3$ \dotfill 42

22.4. Parallel  surfaces in  the motion group  $E(2)$  \dotfill 42

22.5. Parallel surfaces in $SU(2)$ \dotfill 43

22.6. Parallel surfaces in the real special linear group $SL(2,{\mathbb R})$ \dotfill 44

22.7. Parallel surfaces in in non-unimodular three-dimensional Lie groups \dotfill 45

22.8. Parallel surfaces in the Heisenberg group $Nil_3$ \dotfill 45

\noindent Section 23. Parallel  surfaces in three-dimensional Lorentzian Lie groups  \dotfill  46

23.1. Three-dimensional Lorentzian Lie groups \dotfill 46

23.2. Classification of parallel surfaces in three-dimensional  Lorentzian Lie groups  \dotfill 47

\noindent Section 24. Parallel surfaces in reducible three-spaces  \dotfill 49

24.1. Classification of parallel surfaces in reducible three-spaces \dotfill 49

24.2.  Parallel surfaces in Walker three-manifolds  \dotfill 50

\noindent Section 25.  Bianchi-Cartan-Vranceasu spaces \dotfill 51

25.1. Basics on Bianchi-Cartan-Vranceasu spaces \dotfill 51

25.2. $B$-scrolls \dotfill 51

25.3. Parallel surfaces in Bianchi-Cartan-Vranceasu spaces \dotfill 52

\noindent Section 26. Parallel surfaces in homogeneous three-spaces  \dotfill 52

26.1. Homogeneous  three-spaces \dotfill 52

26.2. Classification of parallel surfaces in  homogeneous  Lorentzian three-spaces \dotfill 52

\noindent Section 27. Parallel surfaces in Lorentzian symmetric three-spaces  \dotfill 53

27.1. Lorentzian symmetric three-spaces \dotfill 53

27.2. Classification of parallel surfaces in  homogeneous  Lorentzian three-spaces \dotfill 54

\noindent Section 28.  Three natural extensions of parallel submanifolds \dotfill 56

28.1. Submanifolds with parallel mean curvature vector \dotfill 56

28.2. Higher order parallel submanifolds \dotfill 56

28.3. Semi-parallel submanifolds \dotfill 56

\noindent {\bf References} \dotfill 57

\section{Introduction}  

In Riemannian geometry, parallel transport is a way of transporting geometrical data along smooth curves in a Riemannian manifold. Following an important idea of T. Levi-Civita \cite{Levi17} in 1917,  one can transport vectors of a Riemannian manifold along curves so that they stay parallel with respect to the Levi-Civita connection (or Riemannian connection). Afterwards,  a general theory of parallel transportation of tensor fields in Riemannian geometry was studied in the  1920s by T. Levi-Civita, J. A. Schouten, J. D. Struik, H. Weyl, E. Cartan, B. L. van der Waerden and E. Bortolotti  among others (cf. e.g., \cite{Lu00}).

For an immersed submanifold $M$ of a Riemannian manifold $(N, \tilde g)$, there exist two important symmetric tensor fields; namely,   the first fundamental form which is the induced metric tensor field $g$ of $M$ and the second fundamental form  $h$ which is a normal bundle valued $(1,2)$-tensor field.

It is well known that the first fundamental form $g$ is a parallel tensor field with respect to the Levi-Civita connection. The submanifold $M$ is called a {\it parallel submanifold} if its second fundamental form $h$ is a parallel tensor field with respect to the van der Waerden-Bortolotti  connection. Thus the extrinsic invariants of a parallel submanifold $M$ do not vary from point to point. Obviously, parallel submanifolds are natural extensions of totally geodesic submanifolds for which the second fundamental form vanishes identically. 

Parallel surfaces in a Euclidean 3-space $\mathbb E^3$ are classified in 1948 by V. F. Kagan in \cite{K48}. Kagan's result states that open parts of planes $\mathbb E^2$, of spheres $S^2$ and of round cylinders $S^1\times \mathbb E^1$ are the only parallel surfaces in $\mathbb E^3$.  For $n>2$, parallel hypersurfaces in Euclidean spaces are classified by U. Simon and A. Weinstein in \cite{SW69}. A general classification theorem of parallel submanifolds in any Euclidean space is archived in 1974 by D. Ferus \cite{F74}.  Since then the study of parallel submanifolds became a very interesting and important research subject in differential geometry.

In this paper, we provide a  comprehensive survey on this important subject in differential geometry  from classical  results to the most recent ones.

\section{Preliminaries}

An immersion from a manifold $M$ into a pseudo-Riemannian manifold $(N,\tilde g)$ is called a {\it pseudo-Riemannian submanifold} if the induced metric $g$ on $M$ is a pseudo-Riemannian metric. For a pseudo-Riemannian submanifold $M$ of $N$,  let $\nabla$ and $\tilde\nabla$ be the Levi-Civita connection of $g$ and $\tilde g$, respectively. Let us denote  the Riemann curvature tensors of $M$ and $N$ by $R$ and $\tilde R$,  respectively, and let  $\<\;\, ,\;\>$ denote  the  associated inner product for both $g$ and $\tilde g$. A pseudo-Riemannian manifold is called a {\it Lorentzian manifold} if its index is one at each point.
 
 A tangent vector $v$ of a pseudo-Riemannian manifold is called {\it space-like} (respectively, {\it time-like}) if $v= 0$ or $\left<v,v\right>>0$ (respectively, $\left<v,v\right><0$). A vector $v$ is called {\it light-like} or {\it null} if $\left<v,v\right>=0$ and $v\ne 0$. A pseudo-Riemannian submanifold $M$ is called {\it spatial} (or {\it space-like}) if each tangent vector vector of $M$ is space-like. 

A submanifold $M$ of a pseudo-Riemannian manifold is called {\it non-degenerate} if the induced metric on $M$ is non-degenerate. In particular, a {\it non-degenerate surface}  of a pseudo-Riemannian manifold is either spatial or Lorentzian.
{\it Throughout this article, we assume that every parallel surface $M$ is {\it non-degenerate}, i.e., the induced metric on $M$ is non-degenerate.}

\subsection{Basic definitions, formulas and equations}

 The formulas of Gauss and Weingarten  of a pseudo-Riemannian submanifold $M$ of a pseudo- Riemannian manifold $(N,\tilde g)$ are given respectively by (cf. \cite{c2011,c1973,c1981})
\begin{align} &\label{1}\tilde \nabla_XY=\nabla_XY+h(X,Y),\\& \label{2}\tilde\nabla_X\xi=-A_\xi X+D_X \xi\end{align}
for vector fields $X,Y$ tangent to $M$ and $\xi$ normal to $M$,
where $h,A$ and $D$ are the second fundamental form, the shape operator and the normal connection of $M$. 
The shape operator and the second fundamental form are related by
\begin{align}\label{3} \tilde g(h(X,Y),\xi)=g(A_{\xi}X,Y)\end{align}
for vector fields $X,Y$ tangent to $M$ and $\xi$ normal to $M$. 
The equations of Gauss, Codazzi and Ricci of $M$ in $N$ are given respectively by
\begin{align}\label{4} &g(R(X,Y)Z,W) =\tilde g(\tilde R(X,Y)Z,W)+  \<h(X,W),h(Y,Z)\>  - \<h(X,Z),h(Y,W)\>, \\ &  \label{5} (\tilde R(X,Y)Z)^\perp=(\bar\nabla_X h)(Y,Z) - (\bar\nabla_Y h)(X,Z),
\\&\label{6} \tilde g(R^D(X,Y)\xi,\eta)=\tilde R(X,Y;\xi,\eta)+\tilde g([A_\xi,A_\eta]X,Y),
\end{align}
for  vectors $X,Y,Z,W$ tangent to $M$ and vector $\xi,\eta$ normal to $M$, where $R^D$ is the normal curvature tensor defined by
\begin{align}\label{7} R^D(X,Y)=[D_X,D_Y]-D_{[X,Y]},\end{align}
and $\bar\nabla h$ denotes the covariant derivative of $h$ with respect to the {\it van der Waerden-Bortolotti connection} $\bar \nabla=\nabla\oplus D$, defined by
\begin{align}\label{8}&(\bar \nabla_X h)(Y,Z) = D_X h(Y,Z) - h(\nabla_X Y,Z) - h(Y,\nabla_X Z).\end{align}

The mean curvature vector $H$ of $M$ in $N$ is given  by 
\begin{align}\label{9}&H= \(\frac{1}{n}\) {\rm Trace}\, h=\(\frac{1}{n}\)\sum_{i=1}^n \epsilon_i h(e_i,e_i),\;\; n=\dim M,\end{align}
where $\{e_1,\ldots,e_n\}$ is an orthonormal frame of $M$ such that $\<e_j,e_k\>=\epsilon_j \delta_{jk}$.

The {\it relative null subspace} ${\mathcal N}_p$ of a pseudo-Riemannian submanifold $M$ in $N$ at  $p\in M$ is defined by
$${\mathcal N}_p=\{X\in T_pM: h(X,Y)=0\; \forall  Y\in T_pM\}.$$
The dimension $\nu_p$ of ${\mathcal N}_p$ is called the {\it relative nullity} at $p$.

The {\it first normal space} at a point $p$ of a pseudo-Riemannian submanifold $M$ in $\tilde M$ is, by definition, the image space, ${\rm Im}\,h(p)$, of the second fundamental form of $M$ at $p$, i.e., 
$$ {\rm Im}\,h(p)=\{h(X,Y): X,Y\in T_pM\}.$$

\subsection{Indefinite real space forms}\label{susec2.2}

Let $(N,\tilde g)$ be a pseudo-Riemannian manifold.
At a point $p\in N$,  a 2-dimensional linear subspace $\pi$ of the tangent space $T_pN$ is called a {\it plane section}. For a given basis  $\{v,w\}$ of a plane section $ \pi$, we define a real number by
$$Q(v,w)=\<v,v\>\<w,w\>-\<v,w\>^2.$$
A plane section $\pi$ is called {\it nondegenerate}  if $Q(u,v)\ne 0$. For a nondegenerate plane section $\pi\subset T_pN$ at $p$, the number
\begin{align} \tilde K(u,v)=\frac{\langle\tilde R(u,v)v,u\rangle}{Q(u,v)}\end{align}
is independent of the choice of basis $\{u,v\}$ for $\pi$ and is called the {\it sectional curvature} $\tilde K(\pi)$ of $\pi$.

A pseudo-Riemannian manifold is said to have {\it constant curvature} if its sectional curvature function is constant. It is well known that if a pseudo-Riemannian manifold $N$ is of constant curvature $c$, then its curvature tensor $\tilde R$ satisfies
\begin{align} \label{11} \tilde R(X,Y)Z=c\{\<Y,Z\>X-\<X,Z\>Y\}.\end{align}

\begin{Example} {\rm (see, e.g., \cite{c2011})
Let $\mathbb E_{t}^n$ denote the  pseudo-Euclidean $n$-space equipped with the canonical pseudo-Euclidean metric of  index $t$ 
given by
\begin{align}\label{12}g_{0}= -\sum_{i=1}^{t} du_{i}^{2} +\sum_{j=t+1}^{n}du_{j}^{2}, \end{align}
where $(u_{1},\ldots,u_{n})$ is a rectangular coordinate system of  $\mathbb E_{t}^n$. 
 For a nonzero real number $c$, we put
\begin{align} &\label{13} S^k_s({\bf x}_0,c)=\left\{{\bf x}\in \mathbb E^{k+1}_s: \<{\bf x}-{\bf x}_0,{\bf x}-{\bf x}_0\>=c^{-1}>0\right\},\;\; s>0,\\&\label{14} H^k_s({\bf x}_0,c)=\left\{{\bf x}\in \mathbb E^{k+1}_{s+1}: \<{\bf x}-{\bf x}_0,{\bf x}-{\bf x}_0\>=c^{-1}<0\right\},\;\; s>0,\\&\label{15} H^k(c)=\left\{{\bf x}\in \mathbb E^{k+1}_{1}: \<{\bf x},{\bf x}\>=c^{-1}<0\;{\rm and}\; x_1>0\right\},
\end{align}
where $\<\;\,,\:\>$ is the associated scalar product. $S^k_s({\bf x}_0,c)$ and $H^k_s({\bf x}_0,c)$ are pseudo-Riemannian manifolds of constant curvature $c$ with index $s$, known as a {\it pseudo-sphere} and a {\it pseudo-hyperbolic space}, respectively. The point ${\bf x}_0$ is called the center of $S^m_s({\bf x}_0,c)$ and  $H^m_s({\bf x}_0,c)$. If ${\bf x}_0$ is the origin $o$ of the pseudo-Euclidean spaces, we denote $S^k_s(o,c)$ and $H^k_s(o,c)$ by  $S^k_s(c)$ and $H^k_s(c)$, respectively.
The pseudo-Riemannian manifolds $\,\mathbb E^k_s,\; S^k_s(c),\; H^k_s(c)$ are the standard
models  of the {\it indefinite real space forms}.  In particular, $\mathbb E^k_1,\, S^k_1(c),\,H^k_1(c)$ are the standard models  of  {\it Lorentzian space forms}.
 
The Riemannian manifolds  $\mathbb E^k,\, S^k(c)$ and $H^k(c)$ (with $s=0$) are of constant curvature, called {\it real space forms}. The Euclidean $k$-space $\mathbb E^k$, the $k$-sphere $S^k(c)$ and the hyperbolic $k$-space $H^k(c)$ are complete simply-connected Riemannian manifolds of constant curvature $0, c>0$ and $c<0$, respectively. 
A complete simply-connected pseudo-Riemannian $k$-manifold, $k\geq 3$, of constant curvature $c$ and with index $s$ is isometric to $\mathbb E^k_s$, or $ S^k_s(c)$ or $H^k_s(c)$ according to $c=0,$ or $ c>0$ or $c<0$, respectively.}
\end{Example}

In the following, we denote a $k$-dimensional  indefinite space form of constant curvature $c$ and index $s$ by $R^k_s(c)$. Also we denote an indefinite space form $R^k_0(c)$ (with index $s=0$) simply by $R^k(c)$.

For a pseudo-Riemannian submanifold $M$ of a  pseudo-Riemannian manifold $R^k_s(c)$ of constant curvature $c$ with index $s$, the equations of Gauss, Codazzi and Ricci reduce to (see, e.g., \cite{c2011})
\begin{align}\label{16} &\<R(X,Y)Z,W\> =c\, (\<X,W\>\<Y,Z\>-\<X,Z\>\<Y,W\>)\\&\notag \hskip1.1in +  \<h(X,W),h(Y,Z)\> - \<h(X,Z),h(Y,W)\>, \\ &  \label{17} (\bar\nabla_X h)(Y,Z) = (\bar\nabla_Y h)(X,Z),
\\&\label{18} \<R^D(X,Y)\xi,\eta\>=\<[A_\xi,A_\eta]X,Y\>
\end{align}
for  vectors $X,Y,Z,W$ tangent to $M$ and $\xi,\eta$ normal to $M$.

\subsection{Gauss image}

The classical Gauss map of a surface in $\mathbb E^3$ was introduced by C. F. Gauss in his fundamental paper on the theory of surfaces \cite{G1827}. He used it to define the Gauss curvature. Since then Gauss maps became one of the important tools in differential geometry.
The classical Gauss map can be extended to arbitrary Euclidean submanifolds as follows:

Let $G(n,m-n)$ denote the Grassmann manifold consisting of linear $n$-subspaces of $\mathbb E^m$.  Then the Grassmann manifold $G(n,m-n)$ admits a canonical Riemannian metric via Pl\"ucker embedding which makes $G(n,m-n)$ into a symmetric space. 
For an $n$-dimensional submanifold $M$ of $\mathbb E^m$, the Gauss map $\Gamma$ of $M$ in $\mathbb E^m$ is defined to be the mapping 
$$\Gamma: M\to G(n,m-n)$$
 which carries a point $p\in M$ into the  linear $n$-subspace of $\mathbb E^m$ obtained via the parallel displacement of the tangent space $T_pM$ of $M$ at $p$. The image $\Gamma(M)$ of $M$ in $G(n,m-n)$  via $\Gamma$ is called the {\it Gauss image} of $M$ (cf. \cite{CY83,CY84}). In the following, we shall assume that the Gauss maps are regular maps.
 
 The following result of B.-Y. Chen and S. Yamaguchi in \cite{CY83} provides a simple characterization of Euclidean submanifolds having totally geodesic Gauss image.
 
 \begin{Theorem}\label{T1} A submanifold $M$ of a Euclidean space has totally geodesic Gauss image if and only if its second fundamental form $h$ satisfies
$$(\bar\nabla_Xh)(Y,Z)=h(\nabla^G_XY,Z)- h(\nabla_XY,Z)$$
for any vector fields $X,Y,Z$ tangent to $M$, where $\nabla$ is the Levi-Civita connection of $M$ and $\nabla^G$ is the Levi-Civita connection of the Gauss image with the induced metric via $\Gamma$.
  \end{Theorem}

\section{Some general properties of parallel submanifolds}\label{sec3}

In this section, we present some basic properties of parallel submanifolds.

\begin{Definition} {\rm A pseudo-Riemannian submanifold $M$ of a pseudo-Riemannian submanifold $(N,\tilde g)$ is called a {\it curvature-invariant submanifold} if each tangent space of $M$ is invariant under the curvature transformation, i.e.,
$ \tilde R(X,Y)(T_pM)\subset T_p(M)$ for any vector fields $X,Y$ tangent to $M$.}
\end{Definition}

The following is an immediate consequence of   equation \e{5} of Codazzi.

\begin{Lemma}\label{L1} Any parallel pseudo-Riemannian submanifold $M$ of a pseudo-Riemannian manifold $(N,\tilde g)$ is curvature-invariant.\end{Lemma}

The following properties of parallel submanifolds are also well-known.

\begin{Lemma}\label{L2}  Every parallel  submanifold $M$ of a Riemannian manifold $(N,\tilde g)$ has constant relative nullity, i.e., the dimension of the relative null subspace is constant.
\end{Lemma}

\begin{Lemma}\label{L3} Every parallel submanifold $M$ of a Riemannian symmetric space $(N,\tilde g)$ is locally symmetric, i.e., the Riemannian curvature tensor $R$ of $M$ satisfies $\nabla R=0$.\end{Lemma}

Further, every  parallel submanifold in $\mathbb E^m$ is of finite type in the sense of Chen (cf. e.g., \cite{c84,c96,book15}). Also, if a rank one compact symmetric space $N$ is regarded as a submanifold of a Euclidean space $\mathbb E^m$ via its first standard embedding, then any parallel submanifold of $N$ via its  first standard embedding is of
finite type in $\mathbb E^m$ (see, e.g., \cite{c96,book15}).

\section{Parallel submanifolds of Euclidean spaces}\label{sec4}

In this section, we present basic properties, characterizations and classification of parallel submanifolds of Euclidean spaces.

\subsection{Gauss map and parallel submanifolds}\label{subsec4.1}

As before, let $G(n, m-n)$ denote the Grassmann manifold of $n$-planes through the origin in $\mathbb E^m$ endowed with its natural Riemannian symmetric space metric and  let $G(\mathbb E^m)$ denote the group of Euclidean motions on $\mathbb E^m$. 

The following result was obtained by J. Vilms in \cite{V72}.

\begin{Theorem}\label{T2} Assume that $M$ is an $n$-dimensional parallel submanifold of  $\mathbb E^m$.  If $M$ is complete, then we have:

\begin{itemize}
\item[{\rm (i)}] If the relative nullity $\nu=0$, then $M$ is a complete totally geodesic submanifold of $G(n,m-n)$.

\item[{\rm (ii)}] If $\nu\geq 1$, then there exists a $(G(\mathbb E^m), \mathbb E^m)$-fibration $\pi:M\to B$, where $B$ is a complete totally geodesic submanifold of $G(n,m-n)$ and the fibres are the leaves of the relative nullity foliation. The metric of $M$ is composed from those on base and fibre, and the fibration admits an integrable connection with totally geodesic horizontal leaves (i.e. it is a totally geodesic Riemannian submersion). 

 \item[{\rm (iii)}]  The original Riemannian connection of $M$, or its projection onto $B$, respectively, coincides with the connection induced from  $G(n,m-n)$.
 
 \item[{\rm (iv)}]  $M$ has nonnegative curvature, and is locally symmetric.
 \end{itemize}
\end{Theorem}

As an application of Theorem \ref{T1}, Chen and S. Yamaguchi \cite{CY83} classified  surfaces  with totally geodesic Gauss image as follows. 

\begin{Theorem}\label{T3} Let $M$ be a surface of $\mathbb E^m$. If $M$ has totally geodesic Gauss image in $G(2,m-2)$, then $M$ is one of the following surfaces:

\begin{itemize}
\item[{\rm (a)}]  A surface  in an affine 3-space $\mathbb E^3$ of $\mathbb E^m$.

\item[{\rm (b)}]  A  surface of $\,\mathbb E^m$ with parallel second fundamental form, i.e., $M$ is a parallel surface.

\item[{\rm (c)}]  A surface in an affine 4-space $\mathbb E^4$ of $\mathbb E^m$ which is locally the Riemannian product of two plane curves of nonzero curvature.

\item[{\rm (d)}]  A complex curve lying fully in ${\mathbb C}^2$, where ${\mathbb C}^2$ denotes an affine $\mathbb E^4$ endowed with some orthogonal almost complex structure.
\end{itemize}
\end{Theorem} 

Another application of Theorem \ref{T1} is the following result of  Chen and Yamaguchi obtained in \cite{CY84}.

\begin{Theorem}\label{T4} A submanifold $M$ of $\mathbb E^m$ is locally the product of some hypersurfaces if and only if $M$ has totally geodesic Gauss image and has flat normal connection.
\end{Theorem} 

Yu A,. Nikolaevskij \cite{Ni93} extended Theorem \ref{T3} in 1993 to the following. 

\begin{Theorem}\label{T5} Let $M$ be an $n$-dimensional submanifold of $\,\mathbb E^m$.
Then $M$ has totally geodesic Gauss image in $G(n,m-n)$ if and only if $M$ is the product of submanifolds such that each of the factors is either 
\begin{itemize}
\item[{\rm (a)}] a real hypersurface, or

\item[{\rm (b)}]  a parallel submanifold, or

\item[{\rm (b)}]  a complex hypersurface.

\end{itemize}
\end{Theorem}

\subsection{Normal sections and parallel submanifolds}\label{subsec4.2}

Let $M$ be an $n$-dimensional submanifold in a Euclidean $m$-space $\mathbb E^m$. For a given  point $p\in M$ and a given unit vector $t$ at $p$ tangent to $M$, the vector $t$ and the normal space $T_p^\perp M$ of $M$  determine an $(m-n+1)$-dimensional subspace $E(p,t)$ in $\mathbb E^m$.
The intersection of $M$ and $E(p,t)$ gives a curve $\gamma_t$ (in a neighborhood of $p$) which is called the {\it normal section} of $M$ at $p$ in the direction $t$ (cf. \cite{c1981,c81,c82}). In general, the normal section $\gamma_t$ is a space curve in $E(p,t)$.

For normal sections, Chen proved the following result in \cite{c1981,c81}.

\begin{Theorem}\label{T6} Let $M$ be an $n$-dimensional $(n>2)$ submanifold of a Euclidean $m$-space $\mathbb E^m$. Then $M$ has planar normal sections if and only if the second fundamental for $h$  and its covariant derivative $\bar \nabla h$ satisfy
\begin{align} h(t,t)\wedge (\bar \nabla_t h)(t,t)=0\end{align}
for any unit vector $t$ tangent to $M$.
\end{Theorem}

An immediate consequence of this theorem is the following.

\begin{Theorem}\label{T7} Every  parallel submanifold $M$ of $\mathbb E^m$ with $n=\dim M>2$ has planar normal sections.
\end{Theorem}

By a {\it vertex} of a planar curve $\gamma(s)$ we mean a point $x$ on the curve such that the curvature function $\kappa(s)$ of $\gamma$ satisfies $\frac{d\kappa^2}{ds}=0$ at $x$.

Another application of Theorem \ref{T6} is the following simple geometric
characterization of parallel submanifolds obtained by Chen in \cite{c1981,c81}.

\begin{Theorem} An  $n$-dimensional  $(n>2)$ submanifold $M$ of a Euclidean space is a parallel submanifold if and only if, for each $p\in M$, each normal section of $M$ at any point $p$  is a planar curve with $p$ as one of its vertices.
\end{Theorem}

For further applications of normal sections, see e.g., \cite{c82,c83,CV81,CV84,CL86,DV86,CL04,c06,Sanchez94,Sanchez99,Sanchez07}.

\subsection{Symmetric submanifolds and parallel submanifolds}\label{subsec4.3}

The notion of {\it extrinsic symmetric submanifolds} was defined by D. Ferus in \cite{F80}. More precisely,  an isometric immersion $\psi:M\to \mathbb E^m$ is called {\it extrinsic symmetric} if for each $p\in M$ there exists an isometry $\phi$ of $M$ into itself such that $\phi(p)=p$ and $\psi\circ \phi=\sigma_p\circ \psi$, where $\sigma_p$ denotes the reflection at the normal space
$T_p^\perp M$ at $p$, i.e., the motion of $\mathbb E^m$ which fixes the space through
$\psi(p)$ normal to $\psi_*(T_pM)$ and reflects $\psi(p)+\psi_*(T_pM)$ at $\psi(p)$. The immersed submanifold $\psi :M\to \mathbb E^m$ is said to be {\it extrinsic locally symmetric} if each point $p\in M$ has a neighborhood $U$ and an isometry $\phi$ of $U$ into itself such that $\phi(p)=p$ and $\psi\circ \phi=\sigma_p\circ \psi$ on $U$. In other words, a submanifold $M$ of $\mathbb E^m$ is extrinsic locally symmetric if each point $p\in M$ has a neighborhood which is invariant under the reflection of $\mathbb E^m$ with respect to the normal space at $p$.

D. Ferus \cite{F80} proved the following result.

\begin{Theorem}\label{T9}  Extrinsic locally symmetric submanifolds of Euclidean spaces have parallel second fundamental form and vice versa.
\end{Theorem}

 Symmetric submanifolds were classified completely by D. Ferus in \cite{F74} as being a very special class of orbits of isotropy representations of semisimple symmetric spaces. For some symmetric spaces $N$, a distinguished class of isotropy orbits (the symmetric $R$-spaces) are symmetric spaces. They are symmetric submanifolds in the corresponding tangent space $T_oN$ of $N$. If $N$ is non-compact, the projection of these symmetric submanifolds from $T_oN$ into $N$ via the exponential map at $o$ provides examples of symmetric submanifolds in $N$. 

In \cite{BENT05}, J. Berndt et. al. extended these symmetric submanifolds to larger one-parameter families of symmetric submanifolds, and proved that if $N$ is irreducible and of rank greater than or equal to 2, then every symmetric submanifold of $N$ arises in this way. This result yields the full classification of symmetric submanifolds in Riemannian symmetric spaces.
For symmetric submanifolds in non-flat Riemannian manifolds of constant curvature, see \cite{St79,BR83,Bl85}.

\subsection{Extrinsic $k$-symmetric submanifolds as $\nabla^c$-parallel submanifolds}\label{subsec4.4}

A canonical connection on a Riemannian manifold $(M,g)$ is defined as any metric connection $\nabla^c$ on $M$ such that the difference tensor $\hat D$ between $\nabla\sp c$ and the Levi-Civita connection $\nabla$ of $(M,g)$ is $\nabla\sp c$-parallel. 
An embedded submanifold $M$ of $\mathbb E^m$ is said to be {\it extrinsic homogeneous} with constant principal curvatures if, for any given $p,q\in M$ and  a given piecewise differentiable curve $\gamma$ from $p$ to $q$, there is an isometry $\phi$ of $\mathbb E^m$ satisfying (1) $\phi(M)=M$, (2) $\phi(p)=q$, and (3) $\phi_*{}_p:T^\perp_p M\to T^\perp_qM$ coincides with the
$\hat D$-parallel transport along $\gamma$.

C. Olmos and C. S\'anchez extended Ferus' result  in \cite{OS91} to the following.

\begin{Theorem}\label{T10}  Let $M$ be a  compact Riemannian submanifold fully in  $\mathbb E^m$ and let $h$ be its second fundamental form. Then the following three statements are equivalent:

\begin{itemize}
\item[{\rm (1)}]  $M$ admits a canonical connection $\nabla^c$ such that $\nabla^c h=0$, 

\item[{\rm (2)}] $M$ is an extrinsic homogeneous submanifold with constant principal curvatures, 

\item[{\rm (3)}] $M$ is an orbit of an $s$-representation, that is, of an isotropy representation of a semisimple Riemannian symmetric space. 
\end{itemize}
\end{Theorem}

Furthermore, C. S\'anchez defined in \cite{San85} the notion of {\it extrinsic $k$-symmetric submanifolds} of $\mathbb E^m$ and classified such submanifolds for odd $k$.
Moreover,  he proved in \cite{San92} that the extrinsic $k$-symmetric submanifolds are essentially characterized by the property of
having parallel second fundamental form with respect to the  canonical connection of $k$-symmetric space. In particular, the above result implies that every extrinsic $k$-symmetric submanifold of a Euclidean space is an orbit of an $s$-representation.

\section{Symmetric $R$-spaces and parallel submanifolds of real space forms}\label{sec5}

Symmetric spaces are the most beautiful and important Riemannian manifolds. Such spaces arise in a wide variety of situations in both mathematics and physics.  
 This class of spaces contains many prominent examples which are of great importance for various branches of mathematics, like compact Lie groups, Grassmannians and bounded symmetric domains.  Symmetric spaces are also important objects of study in representation theory, harmonic analysis as well as in differential geometry.
 
We refer to \cite{H,c2018,C1987,CN78,CN88} for general information on compact symmetric spaces. 

\subsection{Symmetric $R$-spaces and Borel subgroups}\label{subsec5.1}

An isometry $s$ of a Riemannian manifold is called an {\it involutive} if $s^{2}=id$. A Riemannian manifold $M$ is called a {\it symmetric space} if for each $p\in M$ there is an involutive isometry $s_{p}$ such that $p$ is an isolated fixed point of $s_{p}$; the involutive isometry $s_{p}\ne id$ is called the {\it symmetry at $p$}. 

Let $M$ be a symmetric space. Denote by $G=G_M$ the closure of the group of isometries on $M$ generated by $\{s_p:p\in M\}$ in the compact-open topology. Then $G$ is a Lie group which acts transitively on the symmetric space. Thus the typical {\it isotropy subgroup} $K$, say at a point $o\in M$, is compact and $M=G/K$.  Let $I_0(M)$ denote the connected group of isometries of a compact symmetric Riemannian manifold $M$.

A symmetric $R$-space is a special type of compact symmetric space for which several characterizations were known.  Originally in 1965,  T. Nagano  defined in \cite{Na65} a symmetric $R$-space as a compact symmetric space $M$ which admits a Lie transformation group $P$ which is non-compact and contains the identity component of the isometric group $I_0(M)$ of $M$  as a subgroup. 
 
 In the theory of algebraic groups, a {\it Borel subgroup} of an algebraic group $G$ is a maximal Zariski closed and connected solvable algebraic subgroup (cf. \cite{H72,B01}).  Subgroups between a Borel subgroup $B$ and the ambient group $G$ are called {\it parabolic subgroups}.  Working over algebraically closed fields, the Borel subgroups turn out to be the minimal parabolic subgroups in this sense. Thus $B$ is a Borel subgroup when the homogeneous space $G/B$ is a complete variety which is as large as possible.

In 1965, M. Takeuchi used the terminology {\it symmetric $R$-space} for the first time in \cite{Ta65}.  He gave a cell decomposition of an $R$-space in \cite{Ta65}, which is a kind of generalization of a symmetric $R$-space. Here, by an {\it R-space} we mean $M = G/U$ where $G$ is a connected real semisimple Lie group without center and $U$ is a parabolic subgroup of $G$.
    A compact symmetric space $M$ is said to have a {\it cubic lattice} if a maximal torus of $M$ is isometric to the quotient of ${\mathbb E}^r$ by a lattice of ${\mathbb E}^r$ generated by an orthogonal basis of the same length. 
  
  In  1985,  O. Loos \cite{Lo85} provided another intrinsic characterization of symmetric $R$-spaces which states that a compact symmetric space $M$ is a symmetric $R$-spaces if and only if the unit lattice of the maximal torus of  $M$ is a cubic lattice. The proof of Loos is based on the correspondence between the symmetric $R$-spaces and compact Jordan triple systems.

\subsection{Classification of symmetric $R$-spaces}\label{subsec5.2}

An affine subspace of $\mathbb E^m$ or a symmetric $R$-space $M\subset \mathbb  E^m$,
which is minimally embedded in a hypersphere of $\mathbb  E^m$ as described  in \cite{TK68} by M. Takeuchi and S. Kobayashi, is a parallel submanifold of $\mathbb E^m$. The class of
symmetric $R$-spaces includes (see \cite{TK68}):
\begin{itemize}
\item[{\rm (a)}] all Hermitian symmetric spaces of compact
type,

\item[{\rm (b)}] Grassmann manifolds
$O(p+q)/O(p)\times O(q), Sp(p+q)/Sp(p)\times
Sp(q),$

\item[{\rm (c)}]  the classical groups
$SO(m),\,U(m),\,Sp(m)$,

\item[{\rm (d)}]  $U(2m)/Sp(m),\, U(m)/O(m)$,

\item[{\rm (e)}]  $(SO(p+1)\times  SO(q+1))/S(O(p)\times
O(q))$, where $S(O(p)\! \times\! O(q))$ is 
the subgroup of $SO(p+1)\times SO(q+1)$
consisting of matrices of the form
$$\begin{pmatrix} \varepsilon& 0 & & \\ 0 & A&&\\ &&\varepsilon& 0\\
&&0&B
\end{pmatrix},\;\; \varepsilon=\pm1,\quad A\in
O(p),\quad B\in O(q),$$

\item[{\rm (f)}] the Cayley projective plane ${\mathcal O}P^2$, and 

\item[{\rm (g)}] the three exceptional spaces $E_6/Spin(10)\times T, E_7/E_6\times
T,$ and $E_6/F_4.$
\end{itemize}

\subsection{Ferus'  theorem}\label{subsec5.3}

A  classification theorem of parallel submanifolds in Euclidean spaces was obtained in 1974 by D. Ferus \cite{F74}. He proved that essentially these submanifolds mentioned above exhaust all parallel submanifolds of $\mathbb  E^m$ in the following sense.

\begin{Theorem}\label{T11} A complete full parallel submanifold of the Euclidean $m$-space $\mathbb  E^m$ is congruent to 
\begin{itemize}
\item[{\rm (1)}] $M=\mathbb  E^{m_0}\times M_1\times \cdots\times M_s\subset \mathbb E^{m_0}\times \mathbb E^{m_1}\times\cdots\times \mathbb  E^{m_s}=\mathbb  E^m$, $s\geq 0$, 
or to
\item[{\rm (2)}] $M= M_1\times \cdots\times M_s\subset  \mathbb E^{m_1}\times\cdots\times  
\mathbb E^{m_s}=\mathbb  E^m$, $s\geq 1$,
\end{itemize}
where each $M_i\subset \mathbb E^{m_i}$ is an irreducible symmetric $R$-space.
Notice that in case $(1)$ $M$ is not contained in any hypersphere of $\mathbb  E^m$, but in case $(2)$ $M$ is contained in a hypersphere of $\mathbb  E^m$. 
\end{Theorem}

\subsection{Parallel submanifolds in spheres}\label{subsec5.4}

For the standard inclusion of a unit hypersphere $S^{m-1}$ in a Euclidean $m$-space $ \mathbb E^m$, a submanifold $M\subset S^{m-1}$ is a parallel submanifold if and only if $M\subset S^{m-1}\subset  \mathbb E^{m}$ is a parallel submanifold of $\mathbb E^m$. 
Hence, Ferus' classification theorem given in \S5.3 implies that $M$ is a parallel submanifold of $S^{m-1}$ if
and only if $M$ is obtained by a submanifold of type (2).

For parallel submanifolds of spaces of constant curvature, see also \cite{Mir78,T81}.

\subsection{Parallel submanifolds in hyperbolic spaces}\label{subsec5.5}

Parallel submanifolds of a hyperbolic space were classified by M. Takeuchi
\cite{T81}  in 1981  as follows.
 
\begin{Theorem} Let $H^m(\bar c)$ be the hyperbolic $m$-space  defined by
$$H^m(\bar c)=\{(x_0,\ldots,x_{m})\in \mathbb E^{m+1}:-x_0^2+x_1^2+\cdots+x_m^2=\bar c^{-1}, x_0>0\},\;\; \bar c<0.$$ 
If $M$ is a parallel submanifold of $H^m(\bar c)$, then we have:
\begin{itemize}
\item[{\rm (1)}]  If $M$ is not contained in any complete totally geodesic hypersurface of
$H^m(\bar c)$, then $M$ is congruent to the product 
$$H^{m_0}(c_0)\times M_1\times \cdots\times
M_s\subset H^{m_0}(c_0)\times S^{m-m_0-1}(c')\subset H^{m_0}(\bar c)$$ 
with $c_0<0,\, c'>0, 1/c_0+1/c'=1/\bar c,\, s\geq 0$, where  $M_1\times \cdots\times
M_s\subset$ $ S^{m-m_0-1}(c')$ is a parallel
submanifold as described  in Ferus' result.

\item[{\rm (2)}]  If $M$ is contained in a complete totally geodesic hypersurface $N$ of
$H^m(\bar c)$, then $N$ is  isometric to an  $(m-1)$-sphere, or to a Euclidean
$(m-1)$-space, or to a hyperbolic $(m-1)$-space. Consequently, such parallel
submanifolds reduce to the parallel submanifolds  described before.
\end{itemize}\end{Theorem}

\section{Parallel Kaehler submanifolds}\label{sec6}

By a {\it complex space form} $\tilde M^m(4c)$, we mean a complex $m$-dimensional Kaehler manifold of constant holomorphic sectional curvature $4c$. It is well known that a complete simply-connected complex space form $\tilde M^m(4c)$ is holomorphically isometric to a complex projective $m$-space $CP^m(4c)$, a complex Euclidean $m$-space $\mathbb C^m$, or a complex hyperbolic $m$-space $CH^m(4c)$ depending on $c>0$, $c=0$ or $c<0$, respectively.

\subsection{Segre and Veronese maps}\label{subsec6.1}

Let $(z_0^i,\ldots,z_{n_i}^i)$ $(1\leq i\leq s)$ denote the homogeneous coordinates of
$CP^{n_i}$. Define a map:
$$S_{n_1\cdots n_s}:CP^{n_1}\times\cdots\times CP^{n_s}\to CP^n,\quad n=\prod_{i=1}^s (n_i+1)-1,$$
which maps a point $((z_0^1,\ldots,z_{n_1}^1),\ldots,(z_0^s,\ldots,z_{n_s}^s))$ of the product Kaehler manifold $CP^{n_1}\times\cdots\times CP^{n_s}$ to the point 
$(z^1_{i_1}\cdots z^s_{i_s})_{0\leq i_1\leq n_1,\,\ldots\,,0\leq i_s\leq n_s}$ in
$CP^n$. Is it well known that the map
$S_{n_1\cdots n_s}$ is a Kaehler embedding, known as the {\it Segre embedding}.

B.-Y. Chen \cite{C81.1} and Chen and W. E. Kuan \cite{CK83,CK85}  proved the following simple characterization for Segre embeddings for   $n=2$  and for $n\geq 3$, respectively (see, also \cite{C81.2,C81.,C81.3,Chen02}).

\begin{Theorem} Let $M_1,\ldots,M_s$ be Kaehler manifolds of complex dimensions $n_1,\ldots,n_s$, respectively. Then every Kaehler immersion
$$\phi:M_1\times\cdots\times M_s\to CP^n,\quad n=\prod_{i=1}^s (n_i+1)-1,$$ 
of $M_1\times\cdots\times M_s$ into $CP^n$   is locally  the Segre embedding, i.e.,
$M_1,\ldots,M_s$ are open portions of $CP^{n_1},\ldots, CP^{n_s}$, respectively, and
moreover, the Kaehler immersion $\phi$ is congruent to the Segre embedding.
\end{Theorem}

A complex projective $n$-space $CP^n(c)$ of constant holomorphic sectional
curvature $c$ can be holomorphically isometrically embedded into an
$\big({{n+\nu} \choose \nu}-1\big)$-dimensional complex projective space of constant holomorphic sectional curvature $\mu c$ as 
$$(z_0,\ldots,z_n)\to \left(z_0^\nu,\sqrt{\nu}z_0^{\nu-1} z_1,\ldots,\sqrt{\frac{\nu !}{\alpha_0!\cdots \alpha_n !}}z_0^{\alpha_0} \cdots z_n^{\alpha_n},\ldots,z_n^\nu\right),\;\; \sum_{i=0}^n \alpha_i=\nu,$$
 which is called  the $\nu$-th {\it Veronese embedding} of $CP^n(c)$. The degree of the 
 $\nu$-th Veronese embedding is $\nu$ (cf. e.g., page 83 of \cite{Og74}).

The Veronese embeddings were characterized  by A. Ros \cite{Ros86} in terms of holomorphic sectional curvature $H$ in the following result.

\begin{Theorem} If a compact $n$-dimensional Kaehler submanifold $M$ immersed in $CP^m(c)$ satisfies 
$${c\over{\nu+1}}< H\leq {c\over\nu},$$ then $M=CP^n({c\over\nu})$
and the immersion is given by the $\nu$-th Veronese embedding.
\end{Theorem}

\subsection{Classification of parallel Kaehler submanifolds of $CP^m$ and $CH^m$}\label{subsec6.2}

In 1972, K. Ogiue classified parallel complex space forms in complex space forms  in \cite{Og72}. More precisely, he proved the following.

\begin{Theorem}  Let $M^n(c)$ be a complex space form holomorphically isometrically immersed in  another complex space form $M^m(\bar c)$. If the second fundamental form of the immersion is parallel, then either the immersion is totally geodesic or $\bar c>0$ and the immersion is given by the second Veronese embedding.
\end{Theorem}

All complete parallel Kaehler submanifolds of a complex projective space were classified by H. Nakagawa and R. Tagaki \cite{NT76}  in 1976 (also \cite{Ta78} by M. Takeuchi in 1978).

\begin{Theorem} Let $M$ be a complete parallel Kaehler submanifold in $CP^m(c)$.  If $M$ is irreducible, then $M$ is congruent to one of the following
six kinds of Kaehler submanifolds:
\begin{equation}\begin{aligned} &CP^n(c),\;\; CP^{n}\left({c\over 2}\right),\; Q_n=SO(n+2)/SO(n)\times SO(2),
\\&\hskip.13in SU(r+2)/S(U(r)\times U(2)),\; r\geq 3,\;\; SO(10)/U(5),
\\ &\hskip1.1in  E_6/Spin(10)\times SO(2).\end{aligned}\end{equation}
If $M$ is reducible, then $M$ is congruent to $CP^{n_1}\times CP^{n_2}$ with $n=n_1+n_2$ and the embedding is given by the Segre embedding.
\end{Theorem}

On the other hand,  M. Kon \cite{K74} proved in 1974 the following result for parallel Kaehler submanifolds in complex hyperbolic spaces.

\begin{Theorem} Every parallel Kaehler submanifold of $CH^m(-4)$ is totally geodesic.
\end{Theorem}

\subsection{Parallel Kaehler submanifolds of Hermitian symmetric spaces}\label{subsec6.3}

Parallel submanifolds of Hermitian symmetric spaces were studied in 1985 by K. Tsukada \cite{Ts85.1} as follows.

\begin{Theorem} Let $\phi: M\to \tilde M$   be a parallel Kaehler immersion of a connected complete Kaehler manifold M into a simply connected Hermitian symmetric space $\tilde M$. Then M is the direct product of a complex Euclidean space and semisimple Hermitian symmetric spaces. Moreover, $\phi=\phi_2\circ \phi_1$, where $\phi_1$ is a direct product of identity maps and (not totally geodesic) parallel Kaehler embeddings into complex projective spaces, and $\phi_2$ is a totally geodesic Kaehler embedding.
\end{Theorem}

All non-totally geodesic parallel Kaehler embeddings into complex projective spaces have been classified earlier by H. Nakagawa and R. Takagi \cite{NT76}  in 1976.
More precisely, these are Veronese maps and Segre maps applied to complex projective spaces, and the first standard embeddings applied to rank two compact irreducible Hermitian symmetric spaces.

\subsection{Parallel Kaehler manifolds in complex Grassmannian manifolds}\label{sec6.4}

   Let $G^{\bf C}(n,p)$ denote the complex Grassmannian manifold of complex $p$-planes in $\mathbb C^n$. We denote by $S\to G^{\bf C}(n,p)$ the tautological vector bundle over $G^{\bf C}(n,p)$ (cf. e.g., \cite{PS11}).  Since the taulogical bundle $S\to G^{\bf C}(n,p)$ is a subbundle of a trivial bundle $G^{\bf C}(n,p)\times {\mathbb C}^n\to G^{\bf C}(n,p)$, one has the quotient bundle $Q\to G^{\bf C}(n,p)$, which is called the {\it universal quotient bundle}.
   
The holomorphic tangent bundle $T_{1,0}(G^{\bf C}(n,p))$ over $G^{\bf C}(n,p)$ can be identified with the tensor product of holomorphic vector bundles $S^*$ and $Q$, where $S^* \to  G^{\bf C}(n,p)$ is the dual bundle of $S \to G^{\bf C}(n,p)$. If ${\mathbb C}^n$ has a Hermitian inner product, $S$, $Q$ have Hermitian metrics and Hermitian connections and so $G^{\bf C}(n,p)$ has a Hermitian metric induced by the identification of $T_{1,0}(G^{\bf C}(n,p))$ and $S^*\otimes Q$ is called the standard metric on $G^{\bf C}(n,p)$. 

In \cite{KN16}, I. Koga and Y. Nagatomo proved the following result for parallel Kaehler manifolds in a complex Grassmannian manifold.

 \begin{Theorem}\label{T19}  Let   $G^{\bf C}(n,p)$ be the complex Grassmannian manifold of complex $p$-planes in ${\mathbb C}^n$ with the standard metric $h_{Gr}$ induced from a Hermitian inner product on ${\mathbb C}^n$ and $\phi$ be a holomorphic isometric immersion of a compact Kaehler manifold $(M,h_M)$ with a Hermitian metric $h_M$ into $G^{\bf C}(n,p)$. We denote by $Q\to G^{\bf C}(n,p)$ the universal quotient bundle over $G^{\bf C}(n,p)$ of rank $n-p$. Assume that the pull-back bundle of $Q\to G^{\bf C}(n,p)$  is projectively flat. Then $\phi$ has parallel second fundamental form if and only if the holomorphic sectional curvature of $M$ is greater than or equal to 1.
 \end{Theorem}

\section{Parallel totally real submanifolds}\label{sec7}

\subsection{Basics on totally real submanifolds}\label{subsec7.1}

  A {\it totally real submanifold} $M$ of an almost Hermitian manifold $\tilde M$ is a submanifold such that the almost Hermitian structure $J$ of $\tilde M$ carries each tangent vector of
$M$ into the corresponding normal space of $M$ in $\tilde M$, i.e.,
$J(T_{p}M)\subseteq T_{p}^{\perp}M$
for any point $p\in M$ (cf. \cite{CO}).  When $\dim_{\bf R} M=\dim_{\bf C} \tilde M$, the totally real submanifold $N$ in $M$ is also known as a {\it Lagrangian submanifold\/}.

The following  result of Chen and K. Ogiue in \cite{CO} is well-known.

\begin{Theorem}\label{T20} A parallel submanifold $M$ of dimension $\geq 2$ of a non-flat complex space form is either a Kaehler submanifold or a totally real submanifold.
\end{Theorem}

H. Naitoh \cite{Na81} proved in 1981 that the classification of complete totally real parallel
submanifolds in complex projective spaces is reduced to that of certain cubic forms of
$n$-variables. Further, H. Naitoh and M. Takeuchi \cite{NT82} classified in 1982 these submanifolds by the theory of symmetric bounded domains of tube type. 

In 1983,  H. Naitoh \cite{Na83.1,Na83.2} proved the following reduction theorem.

\begin{Theorem}\label{T21} A parallel totally real submanifold of a complex space form $\tilde M^n(c)$ with $c\ne 0$ is either a totally real submanifold which is contained in a totally real totally geodesic submanifold, or a totally real submanifold which is contained in a totally geodesic Kaehler submanifold whose dimension is twice of the dimension of the submanifold.
\end{Theorem}

The classifications of Naitoh and Naitoh-Takeuchi given above rely heavily on the theory of Lie groups and symmetric spaces.

\begin{Remark} Theorem \ref{T21} implies that the classification of complete parallel submanifolds of complex projective space $CP^m(c)$ is reduced to those of D. Ferus \cite{F80} and H. Naitoh and  M. Takeuchi \cite{NT82}.\end{Remark}

\begin{Remark} For parallel totally real submanifolds in a complex hyperbolic space $CH^m$, Theorem \ref{T21}  implies that the classification reduces to those of M. Takeuchi \cite{T81}.
\end{Remark}

\subsection{Parallel Lagrangian submanifolds of $CP^n$}\label{subsec7.2}

  F. Dillen,  H. Li, L. Vrancken and X. Wang gave in \cite{DLVW12}  explicitly and geometrically  classification of parallel Lagrangian submanifolds in $CP^n(4)$  using a different method, which applies the warped products of Lagrangian immersions, called {\it Calabi products}, and the characterization of parallel Lagrangian submanifolds by Calabi products.  For the definition of Calabi products and their characterization, see, e.g., \cite{Bo09,LW11}.
    
 The advantage of this classification given by Dillen et. al. is that it allows the study of details for these submanifolds, in particular, for their reduced cases. The classification theorem they obtained is as follows: 
 
\begin{Theorem}\label{T22}  Let $M$ be a parallel Lagrangian submanifold in $CP^n(4)$. Then either $M$ is totally geodesic, or

\begin{itemize}
\item[{\rm (1)}] $M$ is locally the Calabi product of a point with a lower-dimensional parallel Lagrangian submanifold;
\item[{\rm (2)}] $M$  is locally the Calabi product of two lower-dimensional parallel Lagrangian submanifolds; or
\item[{\rm (3)}] $M$  is congruent to one of the following symmetric spaces: {\rm (a)} $SU(k)/SO(k)$ with $n=k(k+1)/2-1$ and $k\geq 3$, {\rm (b)} $SU(k)$ with $n=k^2-1$ and $k\geq 3$, $SU(2k)/Sp(k)$ with $n=2k^2-k-1$ and $k\geq 3$, or {\rm (c)} $E_6/F_4$ with $n=26$.
\end{itemize}\end{Theorem}

\subsection{Parallel  surfaces of  $CP^2$ and $CH^2$}\label{subsec7.3}

For the explicit classification of parallel surfaces in  $CP^2$ (see  \cite{CDV10}).

\begin{Theorem}\label{T23} If $M$ is a parallel surface in the complex projective plane $CP^2(4)$, then it is either holomorphic or Lagrangian in $CP^2(4)$.
\begin{itemize}
\item[{\rm (a)}]   If $M$ is holomorphic, then locally either
\begin{itemize}
\item[{\rm (a.1)}]  $M$ is a totally geodesic complex projective line $CP^1(4)$ in $CP^2(4)$, or

\item[{\rm (a.2)}]  $M$ is the complex quadric $Q^1$ embedded in $CP^2(4)$ as
$\big\{(z_0,z_1,z_2)\in C P^2(4)\ |\ z_0^2+z_1^2+z_2^2=0\big\},$
\indent where $z_0,z_1,z_2$ are complex homogeneous coordinates on $CP^2(4)$.
\end{itemize}
\item[{\rm (b)}]  If $M$ is Lagrangian, then locally either 
\begin{itemize}
\item[{\rm (b.1)}] $M$ is a totally geodesic real projective plane $RP^2(1)$ in $CP^2(4)$, or

\item[{\rm (b.2)}]  $M$ is a flat surface and the immersion is congruent to $\pi\circ  L$, where $\pi:S^5(1)\to CP^2(4)$ is the Hopf-fibration and $L:M\to S^5(1)\subseteq {\mathbb C}^3$ is
given by
\begin{equation}\begin{aligned}\notag &\hskip.5in L(x,y)=\Bigg(\frac{a\,e^{-ix/a}}{\sqrt{1+a^2}},
\frac{e^{i (ax+by)}}{\sqrt{1+a^2+b^2}}\sin \(\sqrt{1+a^2+b^2}\,y\),
\\& \hskip.3in \frac{e^{i (ax+by)}} {\sqrt{1+a^2}} \left(\cos \(\sqrt{1+a^2+b^2}\,y\)-\frac{i b}{\sqrt{1+a^2+b^2}} \sin\(\sqrt{1+a^2+b^2}\,y\)\right)\Bigg), \end{aligned}\end{equation}
where $a$ and $b$ are  real numbers with $a\ne 0$.
\end{itemize}\end{itemize}
\end{Theorem}

For parallel surfaces in  $CH^2$, we have the following result from \cite{CDV10}.
 
\begin{Theorem}\label{T24} If $M$ is a parallel surface in the complex hyperbolic plane $C H^2(-4)$, then it is either holomorphic or Lagrangian in $C H^2(-4)$. 
\begin{itemize}
\item[{\rm (a)}]  If $M^2$ is holomorphic, then it is an open part of a totally geodesic complex submanifold $CH^1(-4)$ in $CH^2(-4)$.

\item[{\rm (b)}] If $M$ is Lagrangian, then locally either
\begin{itemize}
\item[{\rm (b.1)}]  $M$ is a totally geodesic real hyperbolic plane $RH^2(-1)$ in $CH^2(-4)$, or

\item[{\rm (b.2)}] $M$ is flat and the immersion is congruent to $\pi\circ L$, where $\pi:H^5_1(-1)\to \indent C H^2(-4)$ is the Hopf fibration and $L:M^2\to H^5_1(-1)\subseteq {\mathbb C}_1^3$ is one of the \indent following six maps:\\

\item[{\rm (1)}]  $ L= \displaystyle{ \Bigg( \frac{e^{i(ax+by)}}{\sqrt{1-a^2}} \left( \cosh \(\sqrt{1-a^2-b^2}\, y\)- \frac{i b\,\sinh \(\sqrt{1-a^2-b^2}\, y\)}{\sqrt{1-a^2-b^2}} \right),  }$ 
$$ \displaystyle{  \frac{e^{i(ax+by)}}{\sqrt{1-a^2-b^2}} \sinh \(\sqrt{1-a^2-b^2}\, y\), \frac{a\,e^{ix/a}}{\sqrt{1-a^2}} \Bigg) }, \; a,b\in{\bf R},\; a\neq 0,\;a^2+b^2<1;$$

\item[{\rm (2)}]  $ L(x,y)=  \displaystyle{ \Bigg(\! \(\frac{i}{b}+y\) e^{i(\sqrt{1-b^2}x+by)}, y e^{i(\sqrt{1-b^2}x+by)}, \frac{\sqrt{1-b^2}}{b} e^{ix/\sqrt{1-b^2}} \Bigg) },$
              $\; b\in {\bf R},\ 0<b^2<1;$

\item[{\rm (3)}]  $ L(x,y)= \displaystyle{ \Bigg( \frac{e^{i(ax+by)}}{\sqrt{1-a^2}} \left( \cos \(\sqrt{a^2+b^2-1}\, y\)- \frac{i b\,\sin \(\sqrt{a^2+b^2-1}\, y\)}{\sqrt{a^2+b^2-1}} \right), } $
             $$ \displaystyle{ \frac{e^{i(ax+by)}}{\sqrt{a^2+b^2-1}} \sin \(\sqrt{a^2+b^2-1}\, y\Bigg), \frac{a\,e^{ix/a}}{\sqrt{1-a^2}} \) },\;\; a,b\in{\bf R},\;0<a^2<1,\; a^2+b^2>1; $$

\item[{\rm (4)}]  $ L(x,y)= \displaystyle{ \( \frac{a\,e^{ix/a}}{\sqrt{a^2-1}}, \frac{e^{i(ax+by)}}{\sqrt{a^2+b^2-1}} \sin\(\sqrt{a^2+b^2-1}\,y\), \right. } $
             $$ \displaystyle{ \left. \frac{e^{i(ax+by)}}{\sqrt{a^2-1}} \left( \cos\(\sqrt{a^2+b^2-1}\,y\)- \frac{ib\,\sin\(\sqrt{a^2+b^2-1}\, y\)}{\sqrt{a^2+b^2-1}}  \right) \) },\;\; a,b\in {\bf R},\ a^2>1;$$
             
\item[{\rm (5)}]  $ L(x,y)=\displaystyle{ \frac{e^{i x}}{8b^2} \(i+8b^2(i+x)-4by, i+8b^2x-4by, 4b e^{2i by}\) }$, $\; {\bf R}\ni b\neq 0;$

\item[{\rm (6)}]  $ L(x,y)= \displaystyle{ e^{ix} \( 1+\frac{y^2}{2}-ix, y,
               \frac{y^2}{2}-ix\) }.$
\end{itemize}\end{itemize}\end{Theorem}

\subsection{Parallel totally real submanifolds in nearly Kaehler $S^6$}\label{sec7.4}

Let $\mathcal O$ denote the Cayley numbers.
E. Calabi \cite{Cal58} showed in 1958  that any oriented submanifold $M^6$ of the hyperplane  $\hbox{Im}\,\mathcal O$ of the imaginary octonions  carries a $U(3)$-structure,  i.e., an almost Hermitian structure $J$. 

The almost Hermitian structure $J$ on $S^6(1)\subset \hbox{Im}\,\mathcal O$ is a  nearly Kaehler structure in the sense that the (2,1)-tensor field $G$ on $S^6(1)$, defined by
$G(X,Y) = (\widetilde \nabla_XJ)(Y),$ is skew-symmetric, where $\widetilde \nabla$ is the Riemannian connection on $S^6(1)$. 
The group of automorphisms of this nearly K\"ahler structure is the exceptional simple Lie group $G_2$ which acts transitively on $S^6$ as a group of isometries.

In 1969, A. Gray proved in \cite{Gray69}  the following.

\begin{Theorem}\label{T25}  {\rm (1)} Every almost complex submanifold of the nearly Kaehler $S^6(1)$ is a minimal submanifold, and 
{\rm (2)} the nearly Kaehler $S^6(1)$ has no 4-dimensional almost complex submanifolds. 
\end{Theorem}

N. Ejiri proved in \cite{Ej81} that a 3-dimensional totally real submanifold of  the nearly Kaehler  $S^6(1)$ is minimal and orientable

It was proved by B. Opozda in \cite{Op87} that every 3-dimensional parallel Lagrangian submanifold (respectively, a 2-dimensional totally real and minimal submanifold) of  the nearly Kaehler  $S^6(1)$  is totally geodesic (see also \cite{ZD16}).
Opozda also proved in \cite{Op87} that a 2-dimensional parallel totally real, minimal surface of  the nearly Kaehler  $S^6(1)$  is also totally geodesic.
The same result holds for Lagrangian submanifolds of the nearly Kaehler $S^3\times S^3$; namely,  a (3-dimensional) parallel Lagrangian submanifold of the nearly Kaehler $S^3\times S^3$ is totally geodesic (see, e.g., B. Dioos's PhD thesis \cite{Dioos15}).

\section{Parallel slant submanifolds of complex space forms}\label{sec8}

\subsection{Basics on slant submanifolds}\label{subsec8.1}
 Besides Kaehler and totally real submanifolds in a Kaehler manifold $\tilde M$, there is another important family of submanifolds, called {\it slant submanifolds} (cf. \cite{c90,book90}).  

Let $N$ be a submanifold of a K\"ahler manifold (or an almost Hermitian manifold) $(M,J,g)$. For any vector $X$ tangent to $M$, we put 
\begin{align}\label{1.3} JX=PX+FX,\end{align}
 where $PX$ and $FX$ denote the tangential and the normal components of $JX$,
respectively.  Then $P$ is an endomorphism of the tangent bundle $TN$. For any nonzero vector $X\in T_pN$ at  $p\in N$, the angle $\theta (X)$  between
$JX$ and the tangent space $T_pN$ is called the {\it Wirtinger angle} of $X$.

\vskip.06in
In 1990, the author \cite{c90} introduced   the notion of slant submanifolds as follows.

\begin{Definition} {\rm A submanifold $N$ of an almost Hermitian manifold $(M,J,g)$ is called {\it a slant submanifold} if the Wirtinger angle $\theta (X)$ is independent of the choice of $X \in T_{p}N$ and of $p \in N$. The Wirtinger angle of a slant submanifold is called the {\it slant angle}. A slant submanifold with slant angle $\theta$ is simply called {\it $\theta$-slant}.
}\end{Definition}

Complex submanifolds and totally real submanifolds are exactly $\theta$-slant
submanifolds with $\theta = 0$ and $\theta = \frac{\pi}{2}$, respectively. A slant submanifold is called {\it proper slant} if it is neither complex nor totally real. 

The following basic result on slant submanifolds was proved in \cite{CT91} by Chen and Y. Tazawa.

\begin{Theorem}\label{T26}  Let $M$ be a slant submanifold in a complex Euclidean $m$-space ${\mathbb C}^m$. If $M$ is not totally real, then $M$ is non-compact.
In particular, there do not exist  compact proper slant submanifolds in any complex Euclidean $m$-space.
\end{Theorem}

The next result on slant surface was proved  in \cite{CT00} by Chen and Y. Tazawa. 

\begin{Theorem}\label{T27}  Every proper slant surface  of  ${CP^2}$  or of $CH^2$ is non-minimal. \end{Theorem}

\subsection{Classification of parallel slant submanifolds}\label{subsec8.2}

For parallel slant surfaces in ${\mathbb C}^m$, we have the following classification result.

\begin{Theorem}\label{T28} Let $M$ be a slant surface of ${\mathbb C}^m$. Then $M$ is a parallel surface if and only if $M$ is one of the following surfaces:

\begin{itemize}

\item[{\rm (a)}]  An open portion of a slant plane in  ${\mathbb C}^2\subset {\mathbb C^m}$;

\item[{\rm (b)}] An open portion of the product surface of two plane circles;

\item[{\rm (c)}]  An open portion of a circular cylinder which is contained in a hyperplane of  ${\mathbb C}^2\subset {\mathbb C^m}$
\end{itemize}

If case {\rm (b)} or case {\rm (c)} occurs, the $M$ is totally real.
\end{Theorem}

Theorem \ref{T28} follows from Theorem 1.2 of \cite{book90} and  that every parallel surface of a Euclidean space lies in affine 4-space of the ambient space. 

For higher dimensional parallel slant submanifolds, we have the following result by applying  Theorem \ref{T19},  the list of symmetric $R$-spaces and Ferus' Theorem.

\begin{Theorem}\label{T29}  A proper slant  submanifold of ${\mathbb C}^m$ is parallel if and only if it is an open part of a slant $n$-plane of $\mathbb C^m$.
\end{Theorem}

For further results on slant submanifolds, see, e.g., \cite{c2011,book90,C2020,CV97,CV01}.

\section{Parallel submanifolds of  quaternionic space forms and Cayley plane}\label{sec9}

\subsection{Parallel submanifolds of  quaternionic space forms}\label{subsec9.1}

K. Tsukada \cite{Ts85.2} classified  in 1985 all parallel submanifolds of a quaternionic projective $m$-space $HP^m$.  Tsukada's results states that such submanifolds are either parallel totally real submanifolds in a totally real totally geodesic submanifold $RP^m$, or  parallel
totally real submanifolds  in a totally complex totally geodesic submanifold $CP^m$, or parallel complex submanifolds in a totally complex totally geodesic submanifold $CP^m$, or parallel totally complex submanifolds in a totally geodesic quaternionic submanifold $HP^k$ whose dimension is twice the dimension of the parallel submanifold. 
In \cite{Ts85.2}, K. Tsukada also classified parallel submanifolds of the  non-compact dual of $HP^m$.

\subsection{Parallel submanifolds of the Cayley plane}\label{subsec9.2}

A result of K. Tsukada \cite{Ts85.3} in 1985 states that parallel submanifolds of the Cayley plane $\mathcal OP^2$ are contained either in a totally geodesic quaternion projective  plane $HP^2$ as  parallel submanifolds or in a totally geodesic $8$-sphere as parallel submanifolds. Hence, all these immersions are completely known. 

The non-compact case is treated in a similar way.

\section{Parallel spatial submanifolds in pseudo-Euclidean spaces}\label{sec10}

The first classification result of parallel submanifolds in indefinite real space forms was given by M. A. Magid \cite{Ma84}   in 1984 in which he classified parallel immersions of 
$\mathbb E^n\to \mathbb E^{n+k}_1$, $\mathbb E^n_1\to \mathbb E^{n+2}_1$ and $\mathbb E^n_1\to \mathbb E^{n+k}_2$. He showed that such immersions are either quadratic in nature, like the  flat umbilical immersion with light-like mean curvature vector, or the product of the identity map and previously determined low dimensional maps.
In this section we survey known results on parallel pseudo-Riemannian submanifolds in  indefinite real space forms.

First we recall the next lemma  which is an easy consequence of Erbacher-Magid's reduction theorem (see Lemma 3.1 of \cite{c10.1}).

\begin{Lemma}\label{L:3.1}   Let $\psi: M^n_i\to \mathbb E^m_s$ be an isometric immersion of a pseudo-Riemannian $n$-manifold $M^n_i$ into $\mathbb E^m_s$.  If $M$ is a parallel submanifold, then there exists a complete $(n+k)$-dimensional totally geodesic submanifold $E^*$ such that $\psi(M)\subset E^*$, where $k$ is the dimension of the first normal spaces.
\end{Lemma}

\subsection{Marginally trapped surfaces}\label{subsec10.1}

Now, we recall the notion of marginally trapped surfaces for later use.

 The concept of trapped surfaces, introduced R. Penrose in \cite{penrose65} plays very important role in  the theory of  cosmic black holes. If there is a massive source inside the surface, then close enough to a massive enough source, the outgoing light rays may also be converging; a  trapped surface.  Everything inside is trapped. Nothing can escape, not even light.  It is believed that  there will be a marginally trapped surface, separating the trapped surfaces from the untrapped ones, where  the outgoing light rays are instantaneously parallel.  The surface of a  black hole is  the marginally trapped surface.  As times develops, the marginally trapped surface generates a hypersurface in spacetime,  a {\it trapping horizon.}

Spatial surfaces in pseudo-Riemannian manifolds play important roles in mathematics and physics, in particular in general relativity theory. For instance, a {\it marginally trapped surface} in a spacetime is a spatial surface with light-like mean curvature vector field.

In this article, we also call a Lorentzian surfaces in a pseudo-Riemannian manifold  {\it marginally trapped} (or {\it quasi-minimal\/}) if it has light-like mean curvature vector field. A nondegenerate surface  in a pseudo-Riemannian manifold  is called {\it trapped} (respectively, {\it untrapped\/}) if it has time-like (respectively, space-like) mean curvature vector field.

\subsection{Classification of parallel spatial surfaces in $\mathbb E^m_s$}\label{subsec10.2} 

In this subsection, we provide the classification of parallel spatial surfaces in indefinite space forms with arbitrary index and arbitrary dimension obtained by Chen in \cite{c10.1} as follows.

\begin{Theorem} \label{T30} Let $L: M\to \mathbb E^m_s$ be a parallel isometric immersion of a spatial surface into the pseudo-Euclidean $m$-space $\mathbb E^m_s$. Then, up to dilations and rigid motions of $\mathbb E^m_s$,  we have either
\vskip.05in

\noindent {\rm (A)}  the surface is an open part of one of the following 11 surfaces:

\begin{enumerate}

\item[{\rm (i)}]  a totally geodesic Euclidean 2-plane $\mathbb E^2\subset \mathbb E^m_s$ given by $(0,\ldots,0,u,v);$

\item[{\rm (ii)}] a totally umbilical $S^2(1)$  in a totally geodesic $\mathbb E^3\subset \mathbb E^m_s$ given by
$\big(0,\ldots,0,  \cos u,\sin u\cos v,\sin u\sin v\big);$

\item[{\rm (iii)}] a flat cylinder $\mathbb E^1\times S^1$ lying  in a totally geodesic $\mathbb E^3\subset\mathbb E^m_s$  given by $\big(0,\ldots,0,u, \cos v, \sin v\big);$

\item[{\rm (iv)}]  a flat torus $S^1\times S^1$ in a totally geodesic $ \mathbb E^4\subset \mathbb E^m_s$   given by $\big(0,\ldots,0,a\cos u, a\sin u, b \cos v, b\sin v\big)$ with $a, b>0;$

\item[{\rm (v)}] a real projective plane of curvature $\frac{1}{3}$ lying in a totally geodesic $\mathbb E^5\subset \mathbb E^m_s$ given by
\begin{equation}\begin{aligned}& \hskip.2in\notag  \(0,\ldots,0,\text{\small$\frac{vw}{\sqrt{3}}, \frac{uw}{\sqrt{3}}, \frac{uv}{\sqrt{3}}, \frac{u^2-v^2}{2\sqrt{3}},  \frac{1}{6}\(u^2+v^2-2w^2\)$}\hskip-.02in \),\; u^2+v^2+w^2=3;  \end{aligned}\end{equation} 
\item[{\rm (vi)}] a  hyperbolic 2-plane $H^2$ in a totally geodesic $\mathbb E^3_1\subset \mathbb E^m_s $ as
$\big(\cosh u,0,\ldots,0,\sinh u \cos v,\sinh u\sin v\big);$

\item[{\rm (vii)}]  a flat cylinder  $H^1\times \mathbb E^1$ lying in a totally geodesic $\mathbb E^3_1\subset \mathbb E^4_1$ given by $\big(\cosh u, 0,\ldots,0,\sinh u, v\big);$

\item[{\rm (viii)}] a flat surface $H^1\times S^1$  in a totally geodesic $\mathbb E^4_1\subset \mathbb E^m_s$  given by $\big(a\cosh u, 0,\ldots,0,a\sinh u, b\cos v, b\sin v\big)$ with  $a,b>0;$

\item[{\rm (ix)}] a flat totally umbilical surface of a totally geodesic $\mathbb E^4_1\subset \mathbb E^m_s$ defined by
$$\(u^2+v^2+\text{\small$\frac{1}{4}$},0,\ldots,0,u,v, u^2+v^2-\text{\small$\frac{1}{4}$}\);$$

\item[{\rm (x)}]  a  flat surface $H^1\times H^1$ lying in a totally geodesic $\mathbb E^4_2\subset \mathbb E^m_s$ given by
$$\big(a\cosh u,b\cosh v, 0,\ldots,0,a \sinh u,b\sinh v\big),\; a,b>0;$$

\item[{\rm (xi)}] a surface of  curvature $-\frac{1}{3}$ lying in a totally geodesic $\mathbb E^5_3\subset \mathbb E^m_s$ given by
\begin{equation}\begin{aligned}&\notag  \text{\small$\Bigg(\sinh \Big(\frac{2s}{\sqrt{3}}\Big)- \frac{t^2}{3}-\(\frac{7}{8}+\frac{t^4}{18}\)$}e^{\frac{2s}{\sqrt{3}}}, \, t+\( \text{\small $\frac{t^3}{3}-\frac{t}{4}$}\)e^{\frac{2s}{\sqrt{3}}},\text{\small $\frac{1}{2}$} + \text{\small $\frac{t^2}{2}$} e^{\frac{2s}{\sqrt{3}}},
\\& \hskip.35in 0,\ldots,0,   t+ \(\text{\small $\frac{t^3}{3}$}+ \text{\small $\frac{t}{4}$}\)e^{\frac{2s}{\sqrt{3}}},\,
 \sinh \Big(\text{\small$\frac{2s}{\sqrt{3}}$}\Big)-\frac{t^2}{3}-\text{\small$\(\frac{1}{8}+\frac{t^4}{18}\)$}e^{\frac{2s}{\sqrt{3}}}\text{\small$\Bigg)$}, \; or \end{aligned}\end{equation}
 
\end{enumerate}
 
 \noindent  {\rm (B)} $L=(f_1,\ldots,f_\ell,\phi,f_\ell,\ldots,f_1)$, where  $\phi$ is a surface given by  {\rm (i), (iii), (iv), (vii), (viii), (ix),} or {\rm (x)} from $(A)$ and  $f_1,\ldots,f_\ell$ are polynomials of degree $\leq 2$ in $u,v$.
 \vskip.05in
 \end{Theorem}

 \subsection{Special case: parallel spatial surfaces in $\mathbb E^3_1$}\label{subsec10.3}
 
For parallel surfaces in $\mathbb E^3_1$, Theorem \ref{T30} implies  the following.
 
 \begin{Corollary}\label{C1}   A  parallel spatial surface in $\mathbb E^3_1$
is congruent to an open part of one of the following eight types of surfaces:
\begin{enumerate}

\item[{\rm (1)}]  the Euclidean plane $\mathbb E^2$ given by $(0, u, v)$;
\item[{\rm (2)}]   a hyperbolic plane $H^2$ given by $a(\cosh u \cosh v, \cosh u \sinh v, \sinh u)\,a > 0$;
\item[{\rm (3)}]  a cylinder $H^1\times \mathbb E^1$ defined by $(a \cosh u, a \sinh u, v),\, a > 0$;
\end{enumerate}\end{Corollary}
 
\begin{Remark} The surfaces (1) is totally geodesic, the surfaces (2) is totally umbilical but not totally geodesic and surfaces  (1) and (3) are products of parallel curves in totally geodesic
subspaces.
\end{Remark}

 \section{Parallel spatial surfaces in $S^m_s$}\label{sec11} 

\subsection{Classification of parallel spatial surfaces in $S^m_s$}\label{subsec11.1}

For parallel spatial surfaces in a pseudo-sphere $S^m_s$, we have the following classification theorem proved in \cite{c10.1}.

\begin{Theorem}\label{T31} Let $\psi: M\to S^m_s(1)$ be a parallel immersion of a spatial surface into the unit pseudo-Riemannian $m$-sphere $S^m_s(1)$ and $L=\iota:\psi:M\to \mathbb E^{m+1}_s$ be the composition of $\psi$ and the inclusion $\iota:S^m_s(1)\to \mathbb E^{m+1}_s$. Then   either
\vskip.05in

\noindent {\rm (A)} the surface is congruent to an open part of one of the following 18 surfaces:

\begin{enumerate}

\item[{\rm (1)}]  a totally geodesic 2-sphere $S^2(1)\subset S^m_s(1)$;

\item[{\rm (2)}] a totally umbilical $S^2$ immersed in $S^m_s(1)\subset \mathbb E^{m+1}_s$ as
$$\(0,\ldots,0,r \sin u, r\cos u\cos v, r\cos u\sin v,\sqrt{1-r^2}\),\; \;  0<r<1;$$

\item[{\rm (3)}] a totally umbilical $S^2$ immersed in $S^m_s(1)\subset \mathbb E^{m+1}_s$ as
$$\(\sqrt{r^2-1},0,\ldots,0,r \sin u, r\cos u\cos v, r\cos u\sin v\),\;\; r>1\;\; s\geq 1;$$

\item[{\rm (4)}]  a flat torus $S^1\times S^1$ immersed in $S^m_s(1)\subset \mathbb E^{m+1}_s$ as
$$\(0,\ldots,0,b\cos u, b\sin u,c\cos v,c\sin v,\sqrt{1-b^2-c^2}\),\;\;
b,c>0,\; b^2+c^2\leq 1;$$

\item[{\rm (5)}]  a flat torus $S^1\times S^1$ immersed in $S^m_s(1)\subset \mathbb E^{m+1}_s$ as
$$\(\sqrt{b^2+c^2-1},0,\ldots,0,b\cos u, b\sin u,c\cos v,c\sin v\),\;\;
b,c,s>0,\; b^2+c^2>1;$$

\item[{\rm (6)}] a real projective plane $RP^2$ immersed in $S^m_s(1)\subset \mathbb E^{m+1}_s$ as
\begin{equation}\begin{aligned}& \notag \(0,\ldots,0,\text{\small$\frac{rvw}{\sqrt{3}}, \frac{ruw}{\sqrt{3}}, \frac{ruv}{\sqrt{3}}, \frac{r(u^2-v^2)}{2\sqrt{3}},  \frac{r}{6}(u^2+v^2-2w^2)$},\sqrt{1-r^2}\) \end{aligned}\end{equation}
with $u^2+v^2+w^2=3$ and $0<r\leq 1$;

\item[{\rm (7)}] a real projective plane $RP^2$ immersed in $S^m_s(1)\subset \mathbb E^{m+1}_s$ as
\begin{equation}\begin{aligned}& \notag \(\hskip-.02in \sqrt{r^2-1},0,\ldots,0,\text{\small$\frac{rvw}{\sqrt{3}}, \frac{ruw}{\sqrt{3}}, \frac{ruv}{\sqrt{3}}, \frac{r(u^2-v^2)}{2\sqrt{3}},  \frac{r}{6}(u^2+v^2-2w^2)$}\hskip-.02in\) \end{aligned}\end{equation}
with $u^2+v^2+w^2=3$ and $r> 1,s\geq 1;$

\item[{\rm (8)}]  a hyperbolic 2-plane $H^2$ immersed in $S^m_s(1)\subset \mathbb E^{m+1}_s$ as
$$\(r\cosh u,0,\ldots,0,r\sinh u\cos v, r\sinh u\sin v,\sqrt{1+r^2}\),\;\; r,s>0;$$

\item[{\rm (9)}]  a flat surface $H^1\times H^1$ immersed in $S^m_s(1)\subset \mathbb E^{m+1}_s$ as
$$ \(b \cosh u, c\cosh v,0,\ldots,0,b\sinh u,c\sinh v,\sqrt{1+b^2+c^2}\),\;\; b,c>0,\; s\geq 2;$$

\item[{\rm (10)}] a flat surface $H^1\times S^1$ immersed in $S^m_s(1)\subset \mathbb E^{m+1}_s$ as $$\(b\cosh u,0,\ldots,0,b\sinh u, c\cos v, c\sin v, \sqrt{1+b^2-c^2}\big)\),\;\;
b,c,s>0,\; c^2\leq 1+b^2;$$

\item[{\rm (11)}] a flat surface $H^1\times S^1$ immersed in $S^m_s(1)\subset \mathbb E^{m+1}_s$ as $$\( \sqrt{c^2-b^2-1},b\cosh u,0,\ldots,0,b\sinh u, c\cos v, c\sin v\),\;\; c^2> 1+b^2>1;$$

\item[{\rm (12)}] a flat  surface immersed in $S^m_s(1)\subset \mathbb E^{m+1}_s$ as
$$ r \(u^2+v^2+b+\text{\small$\frac{1}{4}$}, 0,\ldots,0,\text{\small$\frac{\sqrt{1+b r^2}}{r}$},u,v,u^2+v^2+b-\text{\small$\frac{1}{4}$} \),\;\; r,s>0,\; b\geq -r^{-2};$$

\item[{\rm (13)}]   a  flat  surface immersed in $S^m_s(1)\subset \mathbb E^{m+1}_s$ as
$$ r \( u^2+v^2-b+\text{\small$ \frac{1}{4}$},\text{\small$\frac{\sqrt{br^2-1}}{r}$}, 0,\ldots,0,u,v,u^2+v^2-b- \frac{1}{4} \)$$
with $r>0,\; s\geq 2,\;b>r^{-2};$

\item[{\rm (14)}]  a flat surface immersed in $S^m_s(1)\subset \mathbb E^{m+1}_s$ as
\begin{align}&\notag r\(u^2+b- \text{\small$\frac{3}{4}$},0,\ldots,0, \text{\small$\frac{\sqrt{1-(1-b+c^2)r^2}}{r}$},u,c\cos v, c\sin v,  u^2+b-\text{\small$\frac{5}{4}$}\)\end{align} 
with  $r,s>0$ and $b\geq 1+c^2-r^{-2};$

\item[{\rm (15)}]  a flat surface immersed in $S^m_s(1)\subset \mathbb E^{m+1}_s$ as
\begin{align}&\notag  r\(u^2+b- \text{\small$\frac{3}{4}$}, \text{\small$\frac{\sqrt{(1-b+c^2)r^2-1}}{r}$},0,\ldots,0,u,c\cos v, c\sin v,  u^2+b-\text{\small$\frac{5}{4}$}\) \end{align} with $r>0$, $s\geq 2$ and $b< 1+c^2-r^{-2};$

\item[{\rm (16)}]   a flat surface immersed in $S^m_s(1)\subset \mathbb E^{m+1}_s$ as
$$ r \left(v^2-b+\text{$\frac{5}{4}$}, c\cosh u,0,\ldots,0,\text{$\frac{\sqrt{1\hskip-.02in +\hskip-.02in (1\hskip-.02in -\hskip-.02in b\hskip-.02in +\hskip-.02in c^2)r^2}}{r}$},c\sinh u, v, v^2-b+\text{$\frac{3}{4}$}\right)$$
with $c,r>0$, $s\geq 2$ and $b\leq 1+c^2+r^{-2};$

\item[{\rm (17)}]   a flat surface immersed in $S^m_s(1)\subset \mathbb E^{m+1}_s$ as
$$ r \left(v^2-b+\text{$\frac{5}{4}$}, c\cosh u,\text{$\frac{\sqrt{(b\hskip-.02in -\hskip-.02in c^2\hskip-.02in -\hskip-.02in 1)r^2\hskip-.02in -\hskip-.02in 1}}{r}$},0,\ldots,0,c\sinh u, v, v^2-b+\text{$\frac{3}{4}$}\right)$$
with $c,r>0$, $s\geq 3$ and $b> 1+c^2+r^{-2};$

\item[{\rm (18)}]  a  surface of constant negative curvature immersed in $S^m_s(1)\subset \mathbb E^{m+1}_s$ as \begin{equation}\begin{aligned}&\hskip.3in \notag r \(\sinh \Big(\text{$\frac{2s}{\sqrt{3}}$}\Big)-\text{$ \frac{t^2}{3}$}-\(\text{$\frac{7}{8}+\frac{t^4}{18}$}\)e^{\frac{2s}{\sqrt{3}}}, \, t+\( \text{ $\frac{t^3}{3}-\frac{t}{4}$}\)e^{\frac{2s}{\sqrt{3}}},\text{ $\frac{1}{2}$} + \text{$\frac{t^2}{2}$} e^{\frac{2s}{\sqrt{3}}},\right.
\\&\left. \hskip.5in 0,\ldots,0,  t+ \(\text{$\frac{t^3}{3}$}+ \text{$\frac{t}{4}$}\)e^{\frac{2s}{\sqrt{3}}},\, \sinh \Big(\text{$\frac{2s}{\sqrt{3}}$}\Big)-\frac{t^2}{3}-\text{$\(\frac{1}{8}+\frac{t^4}{18}\)$}e^{\frac{2s}{\sqrt{3}}}, \frac{\sqrt{1+r^2}}{r}\)\end{aligned}\end{equation}
 with $r>0$ and $s\geq 3$, or
\end{enumerate}

 \noindent  {\rm (B)}  $L=(f_1,\ldots,f_\ell,\phi,f_\ell,\ldots,f_1)$, where $\phi$ is a surface given by  {\rm (4), (5)} or {\rm (9)-(17)} from $(A)$ and $f_1,\ldots,f_\ell$ are polynomials of degree $\leq 2$ in $u,v$, or
 \vskip.05in
 
 \noindent {\rm (C)} $L=(r,\phi,r)$, where  $r\in {\mathbb R^+}$  and $\phi$ is a surface given by {\rm (1), (2), (3), (6), (7), (8)} or {\rm  (18)} from $(A)$.
 \end{Theorem}

\subsection{Special case: parallel spatial surfaces in $S^3_1$}\label{subsec11.2}

For parallel spatial surfaces in a de Sitter space-time $S^3_1$, Theorem \ref{T31} implies the following.

\begin{Corollary}\label{C2} If $M$ is a  parallel spatial surface in $S^3_1(1)\subset \mathbb E^4_1$, then $M$ is congruent to one of the following ten types of surfaces:
\begin{itemize}

\item[{\rm (1)}] a totally umbilical sphere $S^2$ given locally by
$ (a, b \sin u, b \cos u \cos v, b \cos u\sin v),\;\;  b^2 -a^2= 1;$

\item[{\rm (2)}]  a totally umbilical hyperbolic plane $H^2$ given by
$ (a \cosh u \cosh v, a \cosh u \sinh v, a \sinh u, b) $ with $b^2 -a^2= 1;$

\item[{\rm (3)}] a flat surface $H^1\times S^1$ given by
$(a \cosh u, a \sinh u, b \cos v, b\sin v)$ with  $a^2 +b^2= 1.$

\item[{\rm (4)}]  a totally umbilical  Euclidean $\mathbb E^2$ plane given by
$$\frac{1}{\sqrt{c}}\(u^2+v^2-\frac{3}{4},u^2+v^2-\frac{5}{4},u,v\);$$

\end{itemize}
\end{Corollary}

\begin{Remark} The surfaces (1), (2), and (4) are totally umbilical;  the surfaces (1) with $a = 0$ and (2) with $b = 0$ are totally geodesic; the surfaces (3) and (4) are flat. And the
surface (4) is a totally umbilical isometric immersion of $\mathbb E^2$ into $S^3_1(c)$.
\end{Remark}

\section{Parallel spatial surfaces in $H^m_s$}\label{sec12} 

\subsection{Classification of parallel spatial surfaces in $H^m_s$}

For parallel spatial surfaces in a pseudo-hyperbolic space $H^m_s$, we have the following classification theorem also proved in \cite{c10.1}.

 \begin{Theorem} \label{T32} Let $\psi: M\to H^m_s(-1)$ be a parallel immersion of a spatial surface into the pseudo-hyperbolic $m$-space $H^m_s(-1)$  and let $L=\iota:\psi:M\to \mathbb E^{m+1}_{s+1}$ be the composition of $\psi$ and the inclusion $\iota:H^m_s(-1)\to \mathbb E^{m+1}_{s+1}$. Then   either
\vskip.05in

\noindent {\rm (A)} the surface is congruent to an open part of one of the following 18 surfaces:
\begin{enumerate}

\item[{\rm (1)}]  a totally geodesic  $H^2(-1)$ immersed in $H^m_s(-1)\subset \mathbb E^{m+1}_{s+1}$ as 
$\(\cosh u, 0,\ldots,0, \sinh u\cos v, \sinh u \sin v\)$ with $b>0;$

\item[{\rm (2)}] a totally umbilical $H^2$ immersed in  $H^m_s(-1)\subset \mathbb E^{m+1}_{s+1}$ as 
$$\(r\cosh u,0,\ldots,0, r\sinh u\cos v, r \sinh u \sin v,\text{\small$\sqrt{r^2-1}$}\,\)
\; r>1;$$

\item[{\rm (3)}] a totally umbilical $H^2$ immersed in  $H^m_s(-1)\subset \mathbb E^{m+1}_{s+1}$ as 
$$\qquad \quad\(r\cosh u,\text{\small$\sqrt{1-r^2}$}, 0,\ldots,0, r\sinh u\cos v, r \sinh u \sin v\,\),\; s\geq 1,\; 0<r<1;$$

\item[{\rm (4)}] a totally umbilical $S^2$ immersed in $H^m_s(-1)\subset \mathbb E^{m+1}_{s+1}$ as 
$$\(\text{\small$\sqrt{1+r^2}$}, 0,\ldots,0,r \sin u, r\cos u\cos v, r\cos u\sin v\),\; r>0;$$ 

\item[{\rm (5)}]  a flat torus $S^1\times S^1$ in $H^m_s(-1)\subset \mathbb E^{m+1}_{s+1}$ as 
$\(\text{\small$\sqrt{1+b^2+c^2}$}\,,0,\ldots,0,b\cos u, b\sin u,c\cos v,c\sin v,\),$ with $b,c>0;$

\item[{\rm (6)}] a surface of constant positive curvature immersed in $H^m_s(-1)\subset \mathbb E^{m+1}_{s+1}$ as
\begin{equation}\begin{aligned}&\hskip.5in \notag \(\hskip-.02in \sqrt{1+r^2},0,\ldots,0,\text{\small$\frac{rvw}{\sqrt{3}}, \frac{ruw}{\sqrt{3}}, \frac{ruv}{\sqrt{3}}, \frac{r(u^2-v^2)}{2\sqrt{3}},  \frac{r}{6}(u^2+v^2-2w^2)$}\hskip-.02in \) \end{aligned}\end{equation}
with $u^2+v^2+w^2=3$ and $r>0;$

\item[{\rm (7)}] a  flat surface $H^1\times H^1$ in  $H^m_s(-1)$ as 
$ \(b\cosh u,c\cosh v, 0,\ldots,0, b\sinh u,c\sinh v, \text{\small$\sqrt{b^2+c^2-1}$}\)$
with $b,c,s>0$ and $b^2+c^2\geq 1;$

\item[{\rm (8)}] a  flat surface $H^1\times H^1$ in  $H^m_s(-1)$ as 
$\( \text{\small$\sqrt{1-b^2-c^2}$}\,,b\cosh u,c\cosh v, 0,\ldots,0, b\sinh u,c\sinh v\)$
with $b,c>0$, $s\geq 2$ and $b^2+c^2< 1;$

\item[{\rm (9)}] a flat surface $H^1\times S^1$ in $H^m_s(-1)\subset \mathbb E^{m+1}_{s+1}$ as $\(b\cosh u,0,\ldots,0,b\sinh u, c\cos v, c\sin v,\text{\small$ \sqrt{b^2-c^2-1}$}\,\)$
with $b,c>0$ and $b^2\geq c^2+1$;

\item[{\rm (10)}] a flat surface $H^1\times S^1$ immersed in $H^m_s(-1)$ as $\(\text{\small$ \sqrt{1-b^2+c^2}$}\,, b\cosh u,0,\ldots,0,b\sinh u, c\cos v, c\sin v\)$
with $b,c,s>0$ and $b^2< c^2+1$;

\item[{\rm (11)}] a flat  surface immersed in $H^m_s(-1)\subset \mathbb E^{m+1}_{s+1}$ as
$$ r \(u^2+v^2+b+\text{\small$\frac{1}{4}$}, 0,\ldots,0,\text{\small$\frac{\sqrt{br^2-1}}{r}$},u,v,u^2+v^2+b-\text{\small$\frac{1}{4}$}\,\),\;\; r>1,\; b\geq r^{-2};$$ 

\item[{\rm (12)}] a flat  surface immersed in $H^m_s(-1)\subset \mathbb E^{m+1}_{s+1}$ as
$$ r \(u^2+v^2-b+\text{\small$\frac{1}{4}$}, \text{\small$\frac{\sqrt{br^2+1}}{r}$},0,\ldots,0,u,v,u^2+v^2-b-\text{\small$\frac{1}{4}$}\,\),\;\; r,s>0, \; b\geq -r^{-2};$$ 

\item[{\rm (13)}]  a flat surface immersed in $H^m_s(-1)\subset \mathbb E^{m+1}_{s+1}$ as
\begin{align}&\notag  r \(u^2+b- \text{\small$\frac{3}{4}$},0,\ldots,0, \text{\small$\frac{ \sqrt{(b \hskip-.02in - \hskip-.02in c^2 \hskip-.02in - \hskip-.02in 1)r^2 \hskip-.02in - \hskip-.02in 1}}{r}$},c \cos v, c\sin v, u, u^2+b-\text{\small$\frac{5}{4}$}\)\end{align} with $r>0$,  and $b\geq 1+c^2+r^{-2};$

\item[{\rm (14)}]  a flat surface immersed in $H^m_s(-1)\subset \mathbb E^{m+1}_{s+1}$ as
\begin{align}&\notag  r \(u^2+b- \text{\small$\frac{3}{4}$}, \text{\small$\frac{ \sqrt{1 \hskip-.02in + \hskip-.02in (1 \hskip-.02in - \hskip-.02in b \hskip-.02in + \hskip-.02in c^2)r^2}}{r}$}, 0,\ldots,0,c \cos v, c\sin v, u, u^2+b-\text{\small$\frac{5}{4}$}\)\end{align} with $r,s>0$ and $b<1+c^2+r^{-2};$

\item[{\rm (15)}]   a flat surface immersed in $H^m_s(-1)\subset \mathbb E^{m+1}_{s+1}$ as
$$ r \left(v^2+b+\text{$\frac{5}{4}$}, b\cosh u,0,\ldots,0, \text{$\frac{\sqrt{(1\hskip-.02in +\hskip-.02in b\hskip-.02in +\hskip-.02in c^2)r^2\hskip-.02in -\hskip-.02in 1}}{r}$},b\sinh u, v, v^2+b+\text{$\frac{3}{4}$}\right)$$
with $b,r>0$, $s\geq 1$ and $b\geq r^{-2}-1-c^2;$ 
 
\item[{\rm (16)}]   a flat surface immersed in $H^m_s(-1)\subset \mathbb E^{m+1}_{s+1}$ as
$$ r \left(v^2+b+\text{$\frac{5}{4}$}, b\cosh u,\text{$\frac{\sqrt{1\hskip-.02in -\hskip-.02in (a\hskip-.02in +\hskip-.02in b\hskip-.02in +\hskip-.02in c^2)r^2}}{r}$}, 0,\ldots,0,b\sinh u, v, v^2+b+\text{$\frac{3}{4}$}\right)$$
with $b,r>0$, $s\geq 2$, and $b< r^{-2}-1-c^2;$

\item[{\rm (17)}] a surface of constant negative curvature immersed in  $H^m_s(-1)\subset \mathbb E^{m+1}_{s+1}$ as
 \begin{equation}\begin{aligned}&\notag  r \(\hskip-.03in \sinh \Big(\text{$\frac{2u}{\sqrt{3}}$}\Big)-\text{$ \frac{v^2}{3}$}-\(\text{$\frac{7}{8}+\frac{v^4}{18}$}\)e^{\frac{2u}{\sqrt{3}}}, \, v+\( \text{$\frac{v^3}{3}-\frac{v}{4}$}\)e^{\frac{2u}{\sqrt{3}}},\text{$\frac{1}{2}$} + \text{ $\frac{v^2}{2}$} e^{\frac{2u}{\sqrt{3}}}, \right.
\\&
\left. \hskip.4in 0,\ldots,0,   v+ \(\text{ $\frac{v^3}{3}$}+ \text{$\frac{v}{4}$}\)e^{\frac{2u}{\sqrt{3}}},\,
 \sinh \Big(\text{$\frac{2u}{\sqrt{3}}$}\Big)-\frac{v^2}{3}-\text{$\(\frac{1}{8}+\frac{v^4}{18}\)$}e^{\frac{2u}{\sqrt{3}}},\text{$ \frac{\sqrt{r^2-1}}{r}$}\) \end{aligned}\end{equation} with $r\geq 1$ and $s\geq 2;$
 
 \item[{\rm (18)}] a surface of constant negative curvature  immersed in $H^4_2(-1) \subset H^m_s(-1)\subset \mathbb E^{m+1}_{s+1}$ defined as
 \begin{equation}\begin{aligned}&\notag  r \(\hskip-.03in \sinh \Big(\text{$\frac{2u}{\sqrt{3}}$}\Big)-\text{$ \frac{v^2}{3}$}-\(\text{$\frac{7}{8}+\frac{v^4}{18}$}\)e^{\frac{2u}{\sqrt{3}}}, \, v+\( \text{$\frac{v^3}{3}-\frac{v}{4}$}\)e^{\frac{2u}{\sqrt{3}}},\text{ $\frac{1}{2}$} + \text{$\frac{v^2}{2}$} e^{\frac{2u}{\sqrt{3}}}, \right.
\\&\left. \hskip.4in \text{$ \frac{\sqrt{1-r^2}}{r}$}, 0,\ldots,0,   v+ \(\text{$\frac{v^3}{3}$}+ \text{$\frac{v}{4}$}\)e^{\frac{2u}{\sqrt{3}}},\,
 \sinh \Big(\text{$\frac{2u}{\sqrt{3}}$}\Big)-\frac{v^2}{3}-\text{$\(\frac{1}{8}+\frac{v^4}{18}\)$}e^{\frac{2u}{\sqrt{3}}}\) \end{aligned}\end{equation} with $r< 1$ and $s\geq 3$, or
\end{enumerate}

\noindent  {\rm (B)}  $L=(f_1,\ldots,f_\ell,\phi,f_\ell,\ldots,f_1)$, where $f_1,\ldots,f_\ell$ are polynomials of degree $\leq 2$ in $u,v$ and $\phi$ is a surface given by  {\rm  (5), (7), (8)} or {\rm (11)-(18)} from $(A)$, or
 \vskip.05in
 
 \noindent {\rm (C)} $L=(r,\phi,r)$, where  $r$ is a positive number and $\phi$ is a surface given by {\rm (1)-(4), (6), (9)} or {\rm  (10)} from $(A)$.
\end{Theorem}

\subsection{A parallel spatial surfaces in $H^4_2$}\label{subsec12.2}

There is a famous minimal immersion of the 2-sphere $S^2(\frac{1}{3})$ of curvature $\frac{1}{3}$ into the unit 4-sphere $S^4(1)$, known as the Veronese surface, which is constructed by using spherical harmonic homogeneous polynomials of degree two defined as
\begin{equation}\begin{aligned}& \label{1.3} \text{\small$\(\frac{vw}{\sqrt{3}}, \frac{uw}{\sqrt{3}}, \frac{uv}{\sqrt{3}}, \frac{u^2-v^2}{2\sqrt{3}}, 
\frac{u^2+v^2-2w^2}{6}\)$}, \; \; u^2+v^2+w^2=3.
\end{aligned}\end{equation} 
It is well known that the Veronese surface is the only minimal parallel surface lying fully in $S^4(1)$ (see, e.g.,  \cite{chern70,handbook,K}). 
On the other hand, it was also known that there does not exist minimal  surface of constant Gauss curvature lying fully in the hyperbolic 4-space $H^4(-1)$ (cf. \cite{handbook,K,chen72}).
Furthermore, it was known from \cite{CV09} that there exist no minimal spatial parallel surfaces lying fully in  $H^4_1(-1)$.

 B.-Y.  Chen discovered in \cite{chen2010} a minimal immersion of the hyperbolic plane $H^2(-\frac{1}{3})$ of Gauss curvature $-\frac{1}{3}$  into the unit neutral pseudo-hyperbolic 4-space $H^4_2(-1)$ as follows:
 
The following map ${\mathcal B}:{\bf R}^2\to \mathbb E^5_3$:
\begin{equation}\begin{aligned}&\label{24} {\mathcal B}(s,t)=\Bigg(\sinh \Big(\text{\small$\frac{2s}{\sqrt{3}}$}\Big)-\text{\small$ \frac{t^2}{3}$}-\(\text{\small$\frac{7}{8}+\frac{t^4}{18}$}\)e^{\frac{2s}{\sqrt{3}}}, t+\( \text{\small $\frac{t^3}{3}-\frac{t}{4}$}\)e^{\frac{2s}{\sqrt{3}}},
\\& \hskip.1in
 \text{\small $\frac{1}{2}$} + \text{\small $\frac{t^2}{2}$} e^{\frac{2s}{\sqrt{3}}},  t+ \(\text{\small $\frac{t^3}{3}$}+ \text{\small $\frac{t}{4}$}\)e^{\frac{2s}{\sqrt{3}}},
 \sinh \Big(\text{\small$\frac{2s}{\sqrt{3}}$}\Big)-\frac{t^2}{3}-\text{\small$\(\frac{1}{8}+\frac{t^4}{18}\)$}e^{\frac{2s}{\sqrt{3}}}\Bigg) \end{aligned}\end{equation}
 was introduced in \cite{chen2010}.
 It is direct to verify that the position vector field $x$ of $\mathcal B$ satisfies $\<x,x\>=-1$ and the induced metric  is given by $g=ds^2+e^{\frac{2s}{\sqrt{3}}}dt^2.$ Thus ${\mathcal B}$ defines an isometric immersion  
$\psi_{\mathcal B}:H^2(-\tfrac{1}{3})\to H^4_2(-1) $ 
  of the hyperbolic plane $H^2(-\frac{1}{3})$ of curvature $-\frac{1}{3}$ into $H^4_2(-1)$.

In \cite{chen2010}, Chen characterized this parallel immersion $\psi_{\mathcal B}:H^2(-\tfrac{1}{3})\to H^4_2(-1) $ as the following.

 \begin{Theorem}\label{T33} Up to rigid motions, the isometric immersion $\psi_{\mathcal B}:H^2(-\tfrac{1}{3})\to H^4_2(-1)$ defined by \eqref{24}  is the only minimal  parallel spatial  surface lying fully in $ H^4_2(-1)$.
\end{Theorem}

\begin{Remark}  Although our construction of this minimal surface in $H^4_2(-1)$ is  quite different  from the Veronese surface  given by \e{1.3}, we show in \cite{chen2010} that  this parallel surface defined by \e{24} does share several important geometric properties with Veronese surface.  
\end{Remark}

\subsection{Special case: parallel surfaces in $H^3_1$}\label{subsec12.3}

Theorem \ref{T32} implies the following classification of parallel surfaces in $H^3_1$. 

\begin{Corollary}\label{C3} A  parallel spatial surface in $H^3_1(-1)\subset \mathbb E^4_2$ is
congruent to an open part of one of the following ten types of surfaces:
\begin{itemize}
\item[{\rm (1)}]  a hyperbolic plane $H^2$ defined by $(a, b \cosh u \cosh v, b \cosh u \sinh v, b \sinh u)$, $a^2 + b^2 = 1$;
\item[{\rm (2)}] a surface $H^1\times H^1$ defined by $(a \cosh u, b \cosh v, a \sinh u, b \sinh v)$,
 $a^2 + b^2 = 1$.
\end{itemize}
\end{Corollary}

\begin{Remark} The surfaces (1) of Corollary \ref{C3} are totally umbilical and (1) with $a = 0$ is  totally geodesic,  Further, the surfaces (2) are flat and surface (2) with $a^2 = b^2 =\frac{1}{2}$ is minimal.
\end{Remark}

 \section{Parallel Lorentz surfaces in pseudo-Euclidean spaces}\label{sec13} 

 Lorentzian geometry is a vivid field that represents the mathematical foundation of the general theory of relativity, which is probably one of the most successful and beautiful
theories of physics. An interesting phenomenon for Lorentzian surfaces in Lorentzian Kaehler surfaces states that Ricci equation is a consequence of Gauss and Codazzi equations (see \cite{c09.a}). This indicates that  Lorentzian surfaces have many interesting properties which are different from surfaces in Riemannian manifolds.
In particular, Lorentzian surfaces in indefinite real space forms behaved differently from surfaces in Riemannian space forms. For instance, the family of minimal surfaces in Euclidean spaces is huge (see, e.g., chapter 5 of  \cite{handbook}). In contrast, all Lorentzian minimal surfaces in  a pseudo-Euclidean $m$-space $\mathbb E^m_s$  was completely classified  in \cite{c11.1}  (see also \cite{An11}) as the following.

\begin{Theorem} A Lorentzian surface in a pseudo-Euclidean $m$-space $\mathbb E^m_s$  is  minimal if and only if the immersion takes the form 
$$L(x,y)=z(x)+w(y),$$ where $z$ and  $w$ are null curves satisfying $\<z'(x),w'(y)\>\ne 0$. \end{Theorem}

\subsection{Classification of parallel Lorentzian surfaces in $\mathbb E^m_s$}\label{subsec13.1}
 
In  \cite{c10.2}, we have the following classification theorem for parallel Lorentzian surfaces in an arbitrary pseudo-Euclidean space.

\begin{Theorem}\label{T35} Let $M$ be a parallel Lorentzian surface into the pseudo-Euclidean $m$-space $\mathbb E^m_s, s\geq 1$. Then, up to dilations and rigid motions of $\mathbb E^m_s$, we have either 

\vskip.04in
\noindent $\rm (A)$ the surface is an open portion of one of the following fifteen types of surfaces:

\begin{enumerate}
\item[{\rm (1)}]  a totally geodesic plane $\mathbb E^2_1\subset \mathbb E^m_s$ given by $(x,y)\in \mathbb E^2_1\subset \mathbb E^m_s$;

\item[{\rm (2)}]  a totally umbilical de Sitter space $S^2_1$  in a totally geodesic $\mathbb E^3_1\subset \mathbb E^m_s$ given by $$(\sinh x,\cosh x\cos y,\cosh x\sin y);$$

\item[{\rm (3)}]  a flat cylinder $\mathbb E^1_1\times S^1$ in a totally geodesic $\mathbb E^3_1\subset \mathbb E^m_s$ given by $\big(x,\cos y, \sin y\big)$;

\item[{\rm (4)}]   a flat cylinder $S^1_1\times \mathbb E^1$ in a totally geodesic $\mathbb E^3_1\subset \mathbb E^m_s$ given by $\big(\sinh x,\cosh x, y\big)$;

\item[{\rm (5)}] a flat minimal  surface  in a  totally geodesic $\mathbb E^3_1\subset \mathbb E^m_s$ given by
$$\left(\text{$\frac{1}{6}$}(x-y)^3+x,\text{$\frac{1}{6}$}(x-y)^3+y, \text{$\frac{1}{2}$}(x-y)^2\right);$$

\item[{\rm (6)}]  a flat  surface $S^1_1\times S^1$  in a totally geodesic  $\mathbb E^4_1\subset \mathbb E^m_s$ given by $\big(a\sinh
x, a\cosh x, b\cos y, b\sin y\big),$ with $a,b>0;$

\item[{\rm (7)}]  an anti-de Sitter space $H^2_1$ in  a  totally geodesic  ${\mathbb E}^3_2\subseteq \mathbb E^m_s$ given by
$(\sin x, \cos x \cosh y, \cos x \sinh y);$

\item[{\rm (8)}] a flat minimal surface in a totally geodesic ${\mathbb E}^3_2\subseteq \mathbb E^m_s$ defined by
$$\(\frac{a^2x^2}{2}, \frac{x}{2}-\frac{a^4 x^2}{6}+y, \frac{x}{2}+\frac{a^4 x^2}{6}-y \), \; a>0;$$

\item[{\rm (9)}] a non-minimal flat surface in a totally geodesic ${\mathbb E}^3_2\subseteq \mathbb E^m_s$ defined by
\begin{equation}\begin{aligned}\notag & \text{$ \Bigg($}  \text{$\frac{1}{2b}$}\cos \(\text{$\frac{\sqrt{2b}}{a}$}(a^2x+by)\),  \text{$\frac{1}{2b}$}\sin \(\text{$\frac{\sqrt{2b}}{a}$}(a^2x+by)\),  \text{$\frac{a^2 x-by}{a\sqrt{2b}}$} \text{$\Bigg),$}\;\; a, b>0;
\end{aligned} \end{equation}

\item[{\rm (10)}]  a non-minimal flat surface in a totally geodesic ${\mathbb E}^3_2\subseteq \mathbb E^m_s$ defined by
\begin{equation}\begin{aligned}\notag & \text{$ \Bigg($}  \text{$\frac{a^2 x+by}{a\sqrt{2b}}$}, \text{$\frac{1}{2b}$}\cosh \(\!\text{$\frac{\sqrt{2b}}{a}$}(a^2x-by)\),  \text{$\frac{1}{2b}$}\sinh \(\text{$\frac{\sqrt{2b}}{a}$}(a^2x-by)\) \text{$\Bigg),$}\; a, b>0;\end{aligned} \end{equation}

\item[{\rm (11)}] a flat  surface $H^1_1\times H^1$ in a totally geodesic  $\mathbb E^4_2\subset \mathbb E^m_s$ given by $\big(a\sinh
x, b\cosh v, a\cosh x, b\sinh y\big)$ with $a,b>0;$

\item[{\rm (12)}] a  marginally trapped flat surface in a totally geodesic ${\mathbb E}^4_2\subseteq \mathbb E^m_s$ defined by
\begin{align} \notag &\big( a\cos x \cosh y+b \sin x\sinh y, a\sin x \cosh y-b \cos x\sinh y, \\& \notag \hskip-.1in  b\cos x \cosh y-a \sin x \sinh y,b\sin x \cosh y+a \cos x \sinh y \big),\; a,b\in {\bf R};\end{align}

\vskip.04in
\item[{\rm (13)}]  a marginally trapped flat surface in a totally geodesic ${\mathbb E}^4_2\subseteq \mathbb E^m_s$ given by
\begin{align} \notag &\big((1+a)\sin y-(x+ay)\cos y, (1+a)\cos y+(x+ay)\sin y, \\& \notag \hskip-.1in  (1-a)\sin y  +(x+ay)\cos y,(1-a)\cos y-(x+ay)\sin y\big),\; a\in {\bf R};\end{align}

\item[{\rm (14)}]  a non-minimal flat surface  in a totally geodesic ${\mathbb E}^4_3\subseteq \mathbb E^m_s$ defined by
\begin{equation}\begin{aligned}\notag & \(  \cos \( \frac{\sqrt{b}(a^3 x+by)}{a^{5/2}}\)\!, \sin 
\( \frac{\sqrt{b}(a^3 x+by)}{a^{5/2}} \)\!, \cosh \( \frac{\sqrt{b}(a^3 x-by)}{a^{5/2}} \)\!, \sinh 
\( \frac{\sqrt{b}(a^3 x-by)}{a^{5/2}}\)\),
\end{aligned} \end{equation} with $a,b>0;$

\item[{\rm (15)}] a non-minimal flat surface  in a totally geodesic ${\mathbb E}^4_3\subseteq \mathbb E^m_s$ defined by
\begin{align}\notag &   \Bigg(\frac{\sqrt[\uproot{2} 4]{\delta^2+\varphi^2}\cos \big(\lambda (bx+\sqrt{\delta^2+\varphi^2}y\big)}{\sqrt{2}b\sqrt{\sqrt{\delta^2+\varphi^2}+\d}} ,\frac{\sqrt[\uproot{2} 4]{\delta^2+\varphi^2}\sin \big(\lambda (bx+\sqrt{\delta^2+\varphi^2}y\big)}{\sqrt{2}b\sqrt{\sqrt{\delta^2+\varphi^2}+\delta}},\\&\hskip.2in \notag \:  \frac{\sqrt[\uproot{2} 4]{\delta^2+\varphi^2}\cosh \big(\mu(bx-sqrt{\delta^2+\varphi^2}y\big)}{\sqrt{2}b\sqrt{\sqrt{\delta^2+\varphi^2}-\delta}} ,\frac{\!\sqrt[\uproot{2} 4]{\delta^2+\varphi^2}\sin \big(\mu (bx-\sqrt{\delta^2+\varphi^2}y\big)}{\sqrt{2}b\sqrt{\!\sqrt{\delta^2+\varphi^2}-\delta}}\Bigg)
\end{align}  
with $\delta,\varphi\ne 0, b>0$ and
\begin{align}\notag &\lambda =\text{$\frac{\sqrt{b\sqrt{\d^2+\v^2}+b \d}}{\sqrt{\d^2+\v^2}}$},\, \;\; \mu= \text{$\frac{\sqrt{b\sqrt{\d^2+\v^2}-b \d}}{\sqrt{\d^2+\v^2}}$}, \end{align}
\end{enumerate}
or

\vskip.08in 
\noindent \noindent $\rm (B)$ $M^2_1$ is a flat surface and the immersion takes the form: $(f_1,\ldots,f_\ell,\phi(x,y),f_\ell,\ldots,f_1),$
where $\phi=\phi(x,y)$ is given by one of  {\rm (1), (3)--(6), (8)--(15)} 
 and  $f_1,\ldots,f_\ell\, ( \ell \geq 1)$ are polynomials of degree $\leq 2$ in $x,y$.
\end{Theorem}

\subsection{Classification of parallel Lorentzian surfaces in $\mathbb E^3_1$}\label{subsec13.2}

Theorem \ref{T35} implies the following.

 \begin{Corollary} \label{C4}  A  parallel Lorentzian  surface in the Minkowski 3-space $\mathbb E^3_1$
is congruent to an open part of one of the following eight types of surfaces:
\begin{enumerate}

\item[{\rm (1)}]   the Lorentzian plane $E^2_1:L(u,v) = (u, v, 0)$;
\item[{\rm (2)}]   a de Sitter space $S^2_1:L(u,v)= a(\sinh u, \cosh u \cos v, \cosh u \sin v), a > 0$;
\item[{\rm (3)}]   a cylinder $\mathbb E^1_1 \times S^1:L(u,v) = (u, a \cos v, a \sin v),\, a > 0$;
\item[{\rm (4)}]  a cylinder $S^1_1 \times \mathbb E^1:L(u,v) = (a \sinh u, a \cosh u, v),\, a > 0$;
\item[{\rm (5)}]   the null scroll $\mathbb N^2_1$ with rulings in the direction of $(1, 1, 0)$ of the null cubic given by $\alpha(u)=\(\frac{4}{3}u^3+u, \frac{4}{3}u^3-u,2u^2\)$.
\end{enumerate}\end{Corollary}
 
\begin{Remark} The surface (1) is totally geodesic; the surface (s) is totally umbilical but not totally geodesic, all others are flat; the surfaces (1), (3) and (4) are products of parallel curves in totally geodesic subspaces; the surface (5) is flat and minimal, but not totally geodesic.
\end{Remark}

 \section{Parallel  surfaces in a light cone ${\mathcal L}{\mathcal C}$}\label{sec14}

 The {\it light cone} $\mathcal LC$ of a pseudo-Euclidean $(n+1)$-space $\mathbb E^{n+1}_s$ is defined by $${\mathcal LC}^{n}_{s}=\{x\in \mathbb E^{n+1}_s:\left<x,x\right>=0\}.$$
 A curve in a pseudo-Riemannian manifold is called a {\it null curve} if its velocity vector is a lightlike at each point.

\subsection{Light cones in general relativity}\label{subsec14.1} 

In physics, a {\it space-time} is a time-oriented 4-dimensional Lorentz manifold. 
As with any time-oriented spacetime, the time-orientation is called the {\it future}, and its negative is called the {\it past}. A tangent vector in a future time-cone is called {\it future-pointing}.
Similarly, a tangent vector in the past time-cone is called {\it past-pointing}.

Light cones play a very important role in general relativity. Since signals and other causal influences cannot travel faster than light, the light cone plays an essential role in defining the concept of causality: for a given event $E$, the set of events that lie on or inside the past light cone of $E$ would also be the set of all events that could send a signal that would have time to reach $E$ and influence it in some way. Likewise, the set of events that lie on or inside the future light cone of $E$ would also be the set of events that could receive a signal sent out from the position and time of $E$, so the future light cone contains all the events that could potentially be causally influenced by $E$. Events which lie neither in the past or future light cone of $E$ cannot influence or be influenced by $E$ in relativity.

\subsection{Parallel surfaces in ${\mathcal LC}^3_1\subset \mathbb E^4_1$}\label{subsec14.2}

Parallel surfaces in the light cone ${\mathcal LC}^3_1\subset \mathbb E^4_1$ were classified by Chen and J. Van der Veken in \cite{CV09} as follows.

\begin{Theorem}\label{T36} Let $M$ be a parallel surface of $\mathbb E^4_1$. If $M$ lies in the light cone ${\mathcal LC}^3_1\subset \mathbb E^4_1$, then $M$ is  congruent to an open part of one of the following four types of surfaces:

\begin{itemize}
\item[{\rm (1)}] a totally umbilical surface of positive curvature given by
$ a(1, \cos u \cos v, \cos u \sin v, \sin u) ,\; a > 0 ;$

\item[{\rm (2)}]  totally umbilical surface of negative curvature given by
$a(\cosh u \cosh v, \cosh u \sinh v, \sinh u, 1) ,\; a > 0 ;$

\item[{\rm (3)}] a flat totally umbilical surface given by
$ \(u^2+v^2+\frac{1}{4},u^2+v^2-\frac{1}{4},u,v \);$

\item[{\rm (4)}] a flat surface given by $a(\cosh u, \sinh u, \cos v, \sin v), a > 0.$
\end{itemize}
\end{Theorem}

\subsection{Parallel surfaces in ${\mathcal LC}^3_2\subset \mathbb E^4_2$}\label{subsec14.3}

For parallel surfaces in in the light cone ${\mathcal LC}^3_2\subset \mathbb E^4_2$, we have the following result from \cite{CV09} as well.

\begin{Theorem}\label{T37} Let $M$ be a parallel surface of $\mathbb E^4_2$. If $M$ lies in the light cone ${\mathcal LC}^3_2\subset \mathbb E^4_2$, then $M$ is  congruent to an open part of one of the following eight types of surfaces:

\begin{itemize}
\item[{\rm (1)}] a totally umbilical surface of positive curvature given by
$a(\sinh u, 1, \cosh u \cos v, \cosh u \sin v) ,\; a > 0 ;$

\item[{\rm (2)}] a totally umbilical surface of negative curvature given by
$a(\sin u, \cos u \cosh v, 1, \cos u \sinh v), a > 0 ;$

\item[{\rm (3)}] a totally umbilical flat surface defined by
$$ \(u, u^2+v^2-\frac{1}{4},u^2+v^2+\frac{1}{4}, v \);$$

\item[{\rm (4)}] a flat surface defined by $a(\sinh u, \cosh v, \cosh u, \sinh v), a > 0;$

\item[{\rm (5)}]  a flat surface defined by $a(\sin u, \cos u, \cos v, \sin v), a > 0;$

\item[{\rm (6)}]  a flat surface defined by
\begin{align*} &a(\sinh u \cos v + \sinh u \sin v, \cosh u \sin v -\sinh u \cos v ,
\\&\hskip.2in \cosh u \cos v - \sinh u \sin v, \cosh u \sin v +\sinh u \cos v),\; a > 0;\end{align*}

\item[{\rm (7)}]  a flat surface defined by
$a(cos v- u sin v, sin v + u cos v, cos v + u sin v, sin v- u cos v), a > 0;$

\item[{\rm (8)}]  a flat surface defined by
$a(\cosh u - v \sinh u, \sinh u + v \cosh u, \cosh u + v \sinh u, \sinh u - v \cosh u)$ with $a > 0.$
\end{itemize}
\end{Theorem}

\section{Parallel surfaces in de Sitter space-time $S^4_1$}\label{sec15}

The geometry of 4-dimensional space-time is much more complex than that of 3-dimensional space, due to the extra degree of freedom.
Four-dimensional space-times play extremely important roles in the theory of relativity. In physics, space-time is a mathematical model that combines space and time into a single continuum. Space-time is usually interpreted with space being three-dimensional and time playing the role of a fourth dimension. By combining space and time into a single manifold, physicists have significantly simplified a large number of physical theories, as well as described in a more uniform way the workings of the universe at both the super-galactic and subatomic levels. 

In recent times, physics and astrophysics have played a central role in shaping the understanding of the universe through scientific observation and experiment. 
After Kaluza-Klein's theory, the term space-time has taken on a generalized meaning beyond treating space-time events with the normal 3+1 dimensions. It becomes the combination of space and time.  Some proposed space-time theories include additional dimensions, normally spatial, but there exist some speculative theories that include additional temporal dimensions and even some that include dimensions that are neither temporal nor spatial. How many dimensions are needed to describe the Universe is still a big open question.

\subsection{Classification of parallel spatial surfaces in de Sitter space-time $S^4_1$}\label{subsec15.1}

For parallel spatial surfaces in the  de Sitter space-time $S^4_1(1)$, we have the following classification theorem  proved by Chen and Van der Veken in \cite{CV09}.

\begin{Theorem}\label{T38} If $M$ is a  parallel spatial surface in $S^4_1(1)\subset \mathbb E^5_1$, then $M$ is congruent to one of the following ten types of surfaces:

\begin{itemize}
\item[{\rm (1)}] a totally umbilical sphere $S^2$ given locally by
$(c, b \cos u \cos v, b \cos u \sin v, b \sin u, a),\;  a^2 + b^2- c^2 = 1;$

\item[{\rm (2)}]  a totally umbilical hyperbolic plane $H^2$ given by
$(a \cosh u \cosh v, a \cosh u \sinh v, a \sinh u, b, c)$ with $b^2 + c^2 - a^2 = 1 ;$

\item[{\rm (3)}] a torus $S^1\times S^1$ given by
$(a, b \cos u, b \sin u, c \cos v, c \sin v)$ with $b^2 + c^2 - a^2 = 1 ;$

\item[{\rm (4)}] a flat surface $H^1\times S^1$ given by
$(b \cosh u, b \sinh u, c \cos v, c \sin v, a)$ with $a^2 + c^2 - b^2 = 1 ;$

\item[{\rm (5)}] a totally umbilical flat surface defined by $$\(u^2+v^2+a^2+\frac{1}{4},u^2+v^2+a^2-\frac{1}{4},u,v,\sqrt{1+a^2}\);$$

\item[{\rm (6)}] a flat surface defined by
$$\(v^2-\frac{3}{4}+a^2, a\cos u, a\sin u, v, v^2-\frac{5}{4}+a^2\),\;a>0;$$

\item[{\rm (7)}]  a flat surface defined by
$$\frac{1}{\sqrt{1+a^2}}\(u^2+v^2-\frac{3}{4},u^2+v^2-\frac{5}{4},u,v,a\),\;\; a\in {\bf R};$$

\item[{\rm (8)}] a marginally trapped flat surface defined by
$\frac{1}{2}\(2u^2-1,2u^2-2, 2u,\sin v,\cos v\);$

\item[{\rm (9)}] a marginally trapped flat surface defined by
$$\(\frac{b}{\sqrt{4-b^2}},\frac{\cos u}{\sqrt{2-b}}, \frac{\sin u}{\sqrt{2-b}},\frac{\cos v}{\sqrt{2+b}},\frac{\sin v}{\sqrt{2+b}}\);\; |b|<2;$$

\item[{\rm (10)}] a marginally trapped flat surface defined by
$$\(\frac{\cosh u}{\sqrt{b-2}},\frac{\sinh u}{\sqrt{b-2}}, \frac{\cos v}{\sqrt{2+b}},\frac{\sin  v}{\sqrt{2+b}},\frac{b}{\sqrt{b^2-4}}\);\; b>2.$$
\end{itemize}
\end{Theorem}

For parallel spatial surface in $S^3_1(1)\subset \mathbb E^4_1$, Theorem \ref{T38} implies the following.

\begin{Corollary}\label{C5} If $M$ is a  parallel spatial surface in $S^3_1(c)\subset \mathbb E^4_1,\, c>0$, then $M$ is congruent to one of the following four types of surfaces:
\begin{itemize}
\item[{\rm (1)}] a totally umbilical sphere $S^2$ given locally by
$ (a, b \sin u, b \cos u \cos v, b \cos u\sin v)$ with $  b^2 -a^2= c^{-1};$

\item[{\rm (2)}]  a totally umbilical  Euclidean $\mathbb E^2$ plane given by
$\frac{1}{\sqrt{c}}\(u^2+v^2-\frac{3}{4},u^2+v^2-\frac{5}{4},u,v\);$

\item[{\rm (3)}]  a totally umbilical hyperbolic plane $H^2$ given by
$ (a \cosh u \cosh v, a \cosh u \sinh v, a \sinh u, b) ,$ with $b^2 -a^2= c^{-1};$

\item[{\rm (4)}] a flat surface $H^1\times S^1$ given by
$(a \cosh u, a \sinh u, b \cos v, b\sin v)$ with $ a^2 +b^2= c^{-1}.$
\end{itemize}
\end{Corollary}

\subsection{Classification of parallel Lorentzian surfaces in de Sitter space-time $S^4_1$}\label{subsec15.2}

For parallel Lorentzian surfaces in $S^4_1(1)$, we have the following result also from  \cite{CV09}.

\begin{Theorem}\label{T39} If $M$ is a parallel Lorentzian  surface in $S^4_1(1)\subset \mathbb E^5_1$, then $M$ is congruent to an open part of one of the following two types of surfaces:
\begin{itemize}
\item[{\rm (1)}]  a totally umbilical de Sitter space $S^2_1$ in $S^4_1(1)$ given by
$ (a \sinh u, a \cosh u \cos v, a \cosh u \sin v, b, 0)$ with $ a^2 + b^2 = 1 ;$
\item[{\rm (2)}]  a flat surface $S^1_1\times S^1$ given by
$ (a \sinh u, a \cosh u, b \cos v, b \sin v, 0) , a^2 + b^2 = 1 .$
\end{itemize}
Conversely, each surface defined above is a Lorentzian parallel surface in $S^4_1(1)$\end{Theorem}

\section{Parallel surfaces in anti de Sitter space-time  $H^4_1$}\label{sec16}

\subsection{Classification of parallel spatial  surfaces in  $H^4_1$}\label{subsec16.1}

Parallel  surfaces in the anti de Sitter space-time  $H^4_1(-1)$ were also classified by Chen and Van der Veken in \cite{CV09}.

\begin{Theorem}\label{T40} If $M$ is a  parallel spatial surface in $H^4_1(-1)\subset \mathbb E^5_2$, then $M$ is congruent to one of the following ten types of surfaces:
\begin{itemize}
\item[{\rm (1)}] a totally umbilical sphere $S^2$ given locally by
$(a, c, b \sin u, b \cos u \cos v, b \cos u \sin v)$, $a^2-b^2 + c^2 = 1 ;$

\item[{\rm (2)}] a totally umbilical hyperbolic plane $H^2$ given locally by
$(a, b \cosh u \cosh v, b \cosh u \sinh v, b \sinh u, c)$ with $a^2 + b^2-c^2 = 1 ;$

\item[{\rm (3)}] flat surface $H^1\times S^1$ given by
$(a, b \cosh u, b \sinh u, c \cos v, c \sin v)$ with $a^2 + b^2 -c^2 = 1 ;$

\item[{\rm (4)}] a flat surface $H^1\times H^1$ given by
$ (b \cosh u, c \cosh v, b \sinh u, c \sinh v, a)$ with $ b^2 + c^2-a^2 = 1 ;$

\item[{\rm (5)}] a totally umbilical flat surface defined by
$$\(\sqrt{1-a^2},u^2+v^2+a^2+\frac{1}{4}, u^2+v^2+a^2-\frac{1}{4},u,v\),\; a\in (0,1);$$

\item[{\rm (6)}] a flat surface defined by
$$\(a,b\(u^2+v^2-\frac{3}{4}\),b\(u^2+v^2-\frac{5}{4}\),bu,bv\),\; a^2=1+b^2>1;$$

\item[{\rm (7)}] a flat surface defined by
$$\(v^2+\frac{5}{4}-a^2, a \cosh u, a \sinh u, v, v^2+\frac{3}{4}-a^2\), \; a\ne 0;$$

\item[{\rm (8)}] the marginally trapped flat surface defined by
$$\(u^2+1, \frac{1}{2}\cosh v, u, \frac{1}{2}\sinh v,u^2+\frac{1}{2}\);$$

\item[{\rm (9)}] a marginally trapped flat surface defined by
$$\(\frac{\cosh u}{\sqrt{2-b}},\frac{\cosh v}{\sqrt{2+b}},\frac{\sinh u}{\sqrt{2-b}},\frac{\sinh v}{\sqrt{2+b}},\frac{b}{\sqrt{4-b^2}}\),\;|b|<2;$$
\item[{\rm (10)}]  a flat marginally trapped surface defined by
$$\(\frac{b}{\sqrt{b^2-4}},\frac{\cosh v}{\sqrt{b+2}},\frac{\sinh u}{\sqrt{b+2}},\frac{\cos u}{\sqrt{b-1b}},\frac{\sin u}{\sqrt{b-2}}\),\; b>2.$$
\end{itemize}
Conversely, each surface of the ten types given above is spatial and parallel.
\end{Theorem}

For parallel spatial surfaces in $H^3_1(-1)$,  Theorem \ref{T40} implies the following.

\begin{Corollary} If $M$ is a  parallel spatial surface in $H^3_1(-1)\subset \mathbb E^4_1,\, c>0$, then $M$ is congruent to one of the following two types of surfaces:
\begin{itemize}

\item[{\rm (1)}] a hyperbolic plane $H^2$ defined by
$(a, b \cosh u \cosh v, b \cosh u \sinh v, b \sinh u),\; a^2 +b^2 = 1;$

\item[{\rm (2)}] a surface $H^1\times H^1$ defined by $(a \cosh u, b \cosh v, a \sinh u, b \sinh v),
a^2 +b^2 = 1.$
\end{itemize}
\end{Corollary}

\subsection{Classification of parallel Lorentzian  surfaces in anti de Sitter space-time $H^4_1$}\label{subsec16.2}

Parallel  Lorentzian surfaces in $H^4_1(-1)$ were classified by Chen and J. Van der Veken in \cite{CV09} as follows.

\begin{Theorem}\label{T41} If $M$ is a  parallel Lorentzian surface in $H^4_1(-1)\subset \mathbb E^5_2$, then $M$ is congruent to one of the following twelve types of surfaces:
\begin{itemize}
\item[{\rm (1)}] a totally umbilical de Sitter space $S_1^2$ given by
$ (c, a \sinh u \cos v, a \cosh u \cos v, a \cosh u\sin b, b)$ with  $c^2-a^2 - b^2 = 1;$

\item[{\rm (2)}]  a totally umbilical anti-de Sitter space $H^2_1$ given by
$ (a \sin u , a \cos u \cosh v, a\cos u \sinh v, 0,b)$ with $a^2 - b^2 = 1 ;$

\item[{\rm (3)}] a flat surface $S^1_1\times H^1$ given by
$ (c, a \sinh u, a \cosh u \cos v, a \cosh u \sin v, b)$ with $ c^2-a^2 - b^2 = 1 ;$

\item[{\rm (4)}] a flat surface $H_1^1\times S^1$ given by
$ (a \cos u, a \sin u, b \cos v, b \sin v, c)$ with $a^2 + b^2 - c^2 = 1 ;$

\item[{\rm (5)}] a flat surface $S_1^1\times S^1$ given by 
$(a, b \sinh u, b \cosh u, c \cos v, c \sin v)$ with $a^2 - b^2 - c^2 = 1;$

\item[{\rm (6)}] a totally umbilical flat surface defined by
$\(u^2-v^2-\frac{5}{4}, au, av, a\(u^2-v^2-\frac{3}{4}\),b\)$ with $a^2-b^2=1;$

\item[{\rm (7)}]  a flat surface defined by
$$\Bigg(a \cos v-\frac{a(u-v)}{2}\sin v, a\sin v+\frac{a(u-v)}{2}\cos v,\frac{a(u-v)}{2} \sin v,\frac{a(u-v)}{2}\cos v, b\Bigg),\; a\in {\bf R};$$ 

\item[{\rm (8)}]  a flat surface defined by
$$\Bigg(a \cosh v-\frac{a(u+v)}{2}\sinh v, \frac{a(u+v)}{2}\cosh v, a\sinh v-\frac{a(u+v)}{2} \cosh v,\frac{a(u+v)}{2}\sinh v, b\Bigg)$$ with $a^2-b^2=1;$

\item[{\rm (9)}]  a surface defined by
\begin{align*}& (a \cos u \cosh v-a \tan k \sin u \sinh v, a \sec k \sin u \cosh v,
\\&\hskip.2in a \cos u \sinh v-a \tan k \sin u \cosh v, a \sec k \sin u \sinh v, b) ,\end{align*}
with $ a^2-b^2 = 1, \cos k\ne 0;$

\item[{\rm (10)}] a surface defined by
$$\(\frac{b^2(u^2-k^2-1)-1}{2b^2k}, u, \frac{\cos bv}{b}, \frac{\sin bv}{b},
\frac{b^2(u^2+k^2-1)-1}{2b^2k}\),\; b,k\ne 1;$$

\item[{\rm (11)}] a surface defined by
$$\(\frac{-a^2(v^2+k^2+1)+1}{2a^2k}, \frac{\sinh au}{a}, \frac{\cosh au}{a},v,
\frac{a^2(k^2-v^2-1)-1}{2a^2k}\),\; a,k\ne 1;$$

\item[{\rm (12)}] a surface defined by
\begin{align*} & \Bigg(\frac{(u-v)^4}{24k}+\frac{u^2-v^2-k^2-1}{2k},\frac{1}{6}(u-v)^3+u, \frac{1}{2}(u-v)^2,
\\&\hskip.2in \frac{1}{6}(u-v)^3+v,\frac{(u-v)^4}{24k}+\frac{u^2-v^2+k^2-1}{2k}\Bigg),\; k\ne 0. \end{align*}
\end{itemize}
\end{Theorem}

\subsection{Special case: parallel Lorentzian surfaces in $H^3_1$}\label{subsec16.3}

For parallel Lorentzian surfaces in $H^3_1(-1)$,  Theorem \ref{T41} implies the following.

\begin{Corollary} If $M$ is a  parallel Lorentzian surface in $H^3_1(-1)\subset \mathbb E^4_1$, then $M$ is congruent to one of the following eight types of surfaces:
\begin{itemize}

\item[{\rm (1)}] a de Sitter space $S^2_1$ defined by
$(a, b \sinh u, b \cosh u \sin v, b \cosh u \cos v)$ with $ a^2 -b^2 = 1;$

\item[{\rm (2)}]  the surface $ \(u^2 -v^2 -\frac{5}{4}, u, v, u^2 -v^2 -\frac{3}{4}\);$

\item[{\rm (3)}] an anti-de Sitter space $H^2_1$ defined by
$(a \sin u, a \cos u \cosh v, a \cos u \sinh v, b)$ with $ a^2 -b^2 = 1;$

\item[{\rm (4)}] a surface $S^1_1\times H^1$ defined by
$(a \sinh u, b \cosh v, a \cosh u, b \sinh v)$ with $ b^2 -a^2 = 1;$

\item[{\rm (5)}]  a surface $H^1_1\times S^1$ defined by
$(a \cos u, a \sin u, b \cos v, b \sin v)$ with $ a^2 -b^2 = 1;$

\item[{\rm (6)}] a surface defined by
\begin{align*}&\hskip-.5in \big(\cos u \cosh v - \tan k \sin u \sinh v, \sec k \sin u \cosh v,\\&
\cos u \sinh v- \tan k \sin u \cosh v, \sec k \sin u \sinh v\big),\;\; 
\cos k\ne 0;\end{align*}

\item[{\rm (7)}] the surface defined by
$$\(\cos v -\frac{u-v}{2} \sin v, \sin v + \frac{u-v}{2} \cos v,\frac{u-v}{2} \sin v, \frac{u-v}{2} \cos v\);$$

\item[{\rm (8)}] the surface defined by 
$$\(\cosh v - \frac{u+v}{2} \sinh v, \frac{u+v}{2} \cosh v,\sinh v - \frac{u+v}{2} \cosh v, \frac{u+v}{2} \sinh v\).$$
\end{itemize}
\end{Corollary}

\section{Parallel spatial surfaces in $S^4_2$}\label{sec17}

\subsection{Four-dimensional manifolds with neutral metrics}\label{subsec17.1}

 The metrics of neutral signature $(- - + +)$ appear in many geometric and physics problems in the last 25 years. It has been realized that the theory of integrable systems and the techniques from the Seiberg-Witten theory  can be successfully used to study Kaehler-Einstein and self-dual metrics as well as the self-dual Yang-Mills equations in neutral signature. 
Riemannian manifolds with neutral signature are of special interest since it retains many interesting parallels with Riemannian geometry. 
 Such parallels are particularly evident in four dimensions, where Hodge's star operator is involutory for both positive-definite and neutral signatures. Both signatures possess the decomposition of two-forms into self-dual and anti-self-dual parts without the need to complexify as in the Lorentzian case. 
 
 As an interplay between indefiniteness and parallels with Riemannian geometry for neutral signature, the curvature decomposition in four dimensions for the two signatures allows one to deduce a neutral analogue of the Thorpe-Hitchin inequality for compact Einstein 4-manifolds (cf. e.g., \cite{ML01}). 
 Also, the development of the geometry of neutral signature in the work of H. Ooguri and C. Vafa \cite{OV90} showed that neutral signature arises naturally in string theory as well. 

Para-Kaehler manifolds provide further interesting examples of metrics of neutral signature. Such manifolds play some important roles in super-symmetric field theories as well as in string theory (see, for instance, \cite{C06,CLS06,CMMS04,chen10}).

\subsection{Classification of parallel Lorentzian surfaces in $S^4_2$}\label{subsec17.2}

Complete classification of parallel Lorentzian surfaces in neutral pseudo-sphere $S^4_2(1)$ was obtained by Chen in \cite{Chen10} as follows.

\begin{Theorem} \label{T42} There exist 24 families of parallel Lorentzian surfaces in the neutral pseudo 4-sphere $S^4_2(1)\subset \mathbb E^5_2$:

\begin{itemize}
\item[{\rm (1)}]  a totally geodesic de Sitter spacetime $S^2_1(1)\subset S^4_2(1)\subset \mathbb E^5_2$;

\item[{\rm (2)}]  a flat  surface in a totally geodesic $S^3_1(1)\subset S^4_2(1)$ defined by
\begin{align}\notag & \(\sqrt{a^2+b^2-1}, a\sinh u ,a \cosh u,b \cos v ,b\sin v\),\; a, b>0, a^2+b^2\geq 1; \end{align}

\item[{\rm (3)}]    a flat   surface defined by 
\begin{equation}\begin{aligned}\notag  &\hskip.04in \Big( a \cos u  \sinh v+b\sin u \cosh v,
 \sqrt{a^2+b^2} \sin u \sinh v,  \sqrt{a^2+b^2}  \sin u \cosh v,  \\& \hskip.3in
  a \cos u\cosh v+b\sin u \sinh v,\sqrt{1-a^2}\, \Big),\; a\in (0,1]; \end{aligned}\end{equation} 

\item[{\rm (4)}]   a flat  surface defined by 
$\(a \cos u, a\sin u , 
 b \cos v, b \sin v ,\sqrt{1+a^2-b^2}  \),\; a,b>0,\; b^2\leq 1+a^2;$

\item[{\rm (5)}]    a flat  surface defined by 
\begin{equation}\begin{aligned}\notag & \(ku, pu^2+\text{$  \frac{(1-b^2)\v}{k^2}-\frac{k^2}{4\v}$}, b\sin v, b\cos v,pu^2+\text{$  \frac{(1-b^2)\v}{k^2}+\frac{k^2}{4\v}$}\),\;b,k,p,\v\ne 0;\end{aligned} \end{equation}

\item[{\rm (6)}]    a flat   surface defined by 
$ \(\sqrt{b^2-a^2-1}, a\cosh u, a\sinh u, b\cos v, b\sin v\),\; a,b>0,\; b^2\geq 1+a^2 ;$

\item[{\rm (7)}]   a flat  surface defined by 
\begin{equation}\begin{aligned}\notag & \(pu^2+\text{$  \frac{(b^2-1)\v}{k^2}+\frac{k^2}{4\v}$}, b\sinh v, b\cosh v,ku,pu^2+\text{$  \frac{(b^2-1)\v}{k^2}-\frac{k^2}{4\v}$}\),\;b,k,p,\v\ne 0;\end{aligned} \end{equation}

\item[{\rm (8)}]    a flat  surface given by 
$\big(a \cosh u  ,b\sinh v ,a\sinh u, b \cosh v,\sqrt{1+a^2-b^2}\big),\, a,b>0,\, b^2\leq 1+a^2;$

\item[{\rm (9)}]  a marginally trapped surface of constant curvature one defined by 
$$\( \text{$\frac{xy}{x+y},\frac{2}{x+y},\frac{x-y}{x+y},\frac{2+xy}{x+y}$},0\),\;\; x+y\ne 0;$$

\item[{\rm (10)}]   a flat surface  defined by  $\big(x+xy,y-xy, x-y+xy,1+xy,0\big);$

\item[{\rm (11)}]  a surface of positive curvature $c^2$  defined by 
$$\(\text{$\frac{xy-c^2}{c^2(x+y)}, \frac{2\sqrt{1-c^2}\,y}{c^2(x+y)},\frac{xy+c^2}{c^2(x+y)},\frac{c^2(x+y)-2y}{c^2(x+y)}$},0\),\; c\in (0,1),\; x+y\ne 0;$$

\item[{\rm (12)}]   a surface of positive curvature $c^2$  defined by 
$$\(0,\text{$\frac{xy-c^2}{c^2(x+y)},\frac{xy+c^2}{c^2(x+y)},\frac{c^2(x+y)-2y}{c^2(x+y)}, \frac{2\sqrt{c^2-1}\,y}{c^2(x+y)}$}\),\; c>1,\; x+y\ne 0;$$

\item[{\rm (13)}]   a surface of negative curvature $-c^2$  defined by 
\begin{equation}\begin{aligned}\notag & \text{$\frac{1}{c}$}\Big( \cosh u-\sinh u \tanh  v , \sinh u\tanh v,
  \sinh u -\cosh u \tanh v ,{\sqrt{1+c^2}} ,0 \Big),\, c>0; \end{aligned} \end{equation}

\item[{\rm (14)}]   a flat surface  defined by 
\begin{equation}\begin{aligned}\notag &\Bigg(\frac{1+8c^2+2v}{4c}\cos u+ \frac{1+v}{2c}\sin u , \frac{4c^2-1}{4c}\cos u+
\(c+ \frac{v}{2c}\)\sin u, \\&\hskip.1in  \(\frac{1}{4c}+2c+\frac{v}{2c}\) \cos u+ \frac{v\sin u}{2c},  \frac{4c^2+1}{4c}\cos u+\frac{1+2c^2+v }{2c}\sin u  ,0\Bigg),\; c>0; \end{aligned} \end{equation}

\item[{\rm (15)}]    a flat surface  defined by 
\begin{equation}\begin{aligned}\notag &\Big(e^{u}- \frac{(2c-v)e^{-u}}{8c}, \frac{v e^{u}}{4}- \frac{e^{-u}}{2c},e^{u}+ \frac{(2c-v)e^{-u}}{8c}, \frac{v e^{u}}{4}+\frac{e^{-u}}{2c},0\Big),\; c>0;
\end{aligned} \end{equation}

\item[{\rm (16)}]    a flat surface  defined by 
\begin{equation}\begin{aligned}\notag &\Big( x+\frac{y}{2}+\frac{2c^2y^3}{3},xy+\frac{c^2 y^4}{6},x-\frac{y}{2}+\frac{2c^2 y^3}{3},c y^2,1+xy+\frac{c^2 y^4}{6}\Big),\; c> 0; \end{aligned} \end{equation}

\item[{\rm (17)}]   a flat surface  defined by 
\begin{equation}\begin{aligned}\notag \(a v \sinh u+b\cosh u, a v \cosh u,a v \cosh u+b \sinh u, a v \sinh u, \sqrt{1+b^2}\),\, a, b\ne 0;\end{aligned} \end{equation} 

\item[{\rm (18)}]    a  flat surface  defined by 
$(a\sin u-b v \cos u, a\cos u+b v \cos u, b v \cos u, b v \sin u, \sqrt{1+a^2}),\; a, b\ne 0; $

\item[{\rm (19)}]   a flat surface  defined by 
\begin{equation}\begin{aligned}\notag & \Big( v\cos u+ \frac{\sin u}{c}, v\sin u-  \frac{\cos u}{c},v\cos u- \frac{\sin u}{c},  v\sin u+ \frac{\cos u}{c},1\Big),\, c>0;\end{aligned} \end{equation}

\item[{\rm (20)}]   a flat surface  defined by 
\begin{equation}\begin{aligned}\notag & \( \cos u\cos v-\frac{\sin u\sin v}{c},\cos u\sin v+\frac{\sin u\cos v}{c},
\right. \cos u\cos v+ \frac{\sin u\sin v}{c},\\& \left.\hskip.8in \cos u\sin v-\frac{\sin u\cos v}{c},1\hskip-.01in \),\;  c>0;\end{aligned} \end{equation}

\item[{\rm (21)}]   a  flat surface  defined by 
\begin{equation}\begin{aligned}\notag & \( e^{v}\cos  u+\frac{e^{-v}\sin u}{c},e^{-v}\cos u -\frac{e^{v}\sin u}{c}, e^{v}\cos u-\frac{e^{-v}\sin u}{c}, e^{-v}\cos u+\frac{e^{v}\sin u}{c},
1 \),\; c>0;\end{aligned} \end{equation}

\item[{\rm (22)}]    a flat surface  defined by 
 $ \( e^u+a e^{-u}v, e^u v-a e^{-u}, e^u-a e^{-u}v, e^u v+a e^{-u},1\),\; a\ne 0;$
\item[{\rm (23)}]   a flat surface  defined by 
$ \( e^{u}-a e^{-u}, e^v+a e^{-v}, e^u+a e^{-u}, e^v-a e^{-v},1\), \; a\ne 0; $
\item[{\rm (24)}]   a flat surface  defined by 
$ (a \cosh u \cos v,a\cosh u\sin v,a\sinh u,\cos v,a\sinh u\sin v,\sqrt{1+a^2}), \, a>0.$
\end{itemize}

Conversely, every parallel  immersion  $L:M\to S^4_2(1) \subset \mathbb E^5_2$ of a Lorentzian surface $M$ into the pseudo 4-sphere $S^4_2(1)$ is congruent to an open portion a surface obtained from  one of 24 families of surfaces described above.
\end{Theorem}

\subsection{Classification of parallel Lorentzian surfaces in $H^4_2$}\label{subsec17.3}

 Complete classification of parallel Lorentzian surfaces in neutral pseudo hyperbolic 4-space $H^4_2(-1)\subset \mathbb E^5_3$ was obtained by Chen in \cite{c10.4}, in which he proved that there exist 53 families of parallel Lorentzian surfaces in neutral pseudo hyperbolic 4-space $H^4_2(-1)$.  

Among the 53 families we have:
 one family of totally geodesic anti-de Sitter space-time; one family of marginally trapped surfaces of curvature one;
one family of untrapped flat surfaces;
one family of untrapped  surfaces of positive curvature; 
one family of untrapped  surfaces of negative curvature;
two family of trapped surfaces of negative curvature; 
two families of flat minimal surfaces;
7 families of untrapped flat surfaces;
8 families of marginally trapped flat  surfaces;
9 families of flat surface which can be either trapped or untrapped; 
and 20 families of trapped flat surfaces.

Conversely, every parallel Lorentzian surface in  $H^4_2(-1)$ is congruent to an open portion of a surface obtained from one of the 53 families.

\section{Parallel spatial surfaces in $S^4_3$ and in $H^4_3$}\label{sec18}

Parallel Lorentzian surfaces in $S^4_3(1)$ and in $H^4_3(-1)$ were completely classified by Chen in \cite{c10.5}.

\subsection{Classification of parallel spatial surfaces in $S^4_3$}\label{subsec18.1}

Chen proved in \cite{c10.5} that there are 21 families of parallel Lorentzian surfaces in 
 $S^4_3(1)\subset \mathbb E^5_3$. 
 Among the 21 families, we have: 
 the totally geodesic  de Sitter space-time $S^2_1(1)\subset S^4_3(1)$; 
one family of minimal flat surfaces in $S^4_3(1)$; 
a totally umbilical flat surfaces lying in a totally geodesic $S^3_2(1)\subset S^4_2(1)$;
one family of totally umbilical de Sitter space $S^2_1(c^2)$ in a totally geodesic $S^3_2(1)\subset S^4_2(1)$;
one family of totally umbilical anti-de Sitter space $H^2_1(-c^2)$ lying in a totally geodesic $S^3_2(1)\subset S^4_2(1)$;
four families  of CMC flat surfaces lying in a totally geodesic $S^3_2(1)\subset S^4_2(1)$; 
and 12 families of flat minimal surfaces.

 Conversely, every parallel Lorentzian surface in  $S^4_3(1)\subset \mathbb E^5_3$ is congruent to an open portion of a surface obtained from one of the 21 families.

\subsection{Classification of parallel spatial surfaces  in $H^4_3$}\label{subsec18.2}

For parallel Lorentzian surfaces in  $H^4_3(-1)\subset \mathbb E^5_4$, Chen proved  in \cite{c10.5} the following classification theorem.

\begin{Theorem}  \label{T43} There are six families of parallel Lorentzian surfaces in  $H^4_3(-1)\subset \mathbb E^5_4$:
\begin{itemize}
\item[{\rm (1)}]  A totally geodesic  anti-de Sitter space $H^2_1(-1)\subset H^4_3(-1)$;

\item[{\rm (2)}]  A flat minimal surface in a totally geodesic $H^3_2(-1)\subset H^4_3(-1)$ defined by
\begin{align}\notag &\frac{1}{\sqrt{2}} \Bigg(\! \sin\( ax+\frac{y}{a} \)\hskip-.02in , \cos \( ax+\frac{y}{a} \)\hskip-.02in , \cosh \( ax -\frac{y}{a} \)\hskip-.02in ,\sinh \( ax- \frac{y}{a} \)\hskip-.02in ,0 \Bigg),\;\; a>0;\end{align}

\item[{\rm (3)}]  A totally umbilical anti-de Sitter space $H^2_1(-c^2)$ in a totally geodesic $H^3_2(-1)\subset H^4_3(-1)$ given by 
\begin{equation}\begin{aligned}\notag &\hskip .0in \frac{1}{c}   \Bigg( 0,\sqrt{c^2-1}, \tanh\Big( \frac{cx+cy}{\sqrt{2}} \Big) ,\sinh (\sqrt{2}cy)\tanh \Big( \frac{cx+cy}{\sqrt{2}} \Big)-\cosh(\sqrt{2}cy),\\& \hskip.3in  \sinh(\sqrt{2}cy)-\cosh(\sqrt{2}cy)\tanh \Big(\frac{cx+cy}{\sqrt{2}} \Big) \Bigg),\;\;  c>1; \end{aligned} \end{equation}

\item[{\rm (4)}] A CMC flat surface in a totally geodesic $H^3_2(-1)$ given by 
\begin{equation}\begin{aligned}\notag &\hskip-.6in \Bigg( \frac{\sqrt{ \sqrt{1+b^2}-b}}{\sqrt{2}\sqrt[\uproot{2}4]{1+b^2}}   \cos \!\Big(\frac{\!\sqrt{\sqrt{1+b^2}+b}(a^2 x+\sqrt{1+ b^2}y)}{a}\Big) ,\\&\hskip-.4in   \frac{\sqrt{\sqrt{1+b^2}-b}}{\sqrt{2}\sqrt[\uproot{2}4]{1\!+\!b^2}}   \sin \Big(\frac{\sqrt{\! \sqrt{1+b^2}+b}(a^2 x+\sqrt{1+b^2}y)}{a}\Big) ,\\& \hskip-.2in \frac{\sqrt{\sqrt{1+b^2}+b}}{\sqrt{2}\sqrt[\uproot{2}4]{1+b^2}}  \cosh \!\Big(\frac{\sqrt{ \sqrt{1+b^2}-b}(a^2 x-\sqrt{1+b^2}y)}{a}\Big) , \\&  \frac{\sqrt{\sqrt{1+b^2}+b}}{\sqrt{2}\sqrt[\uproot{2}4]{1+b^2}}   \sin \Big(\frac{\sqrt{\sqrt{1+b^2}-b}(a^2 x-\sqrt{1+b^2}y)}{a}\Big)\! \Bigg),\;\; a,b,c>0; \end{aligned} \end{equation}

\item[{\rm (5)}]  A non-minimal flat surface given by 
\begin{equation}\begin{aligned}\notag &\hskip -.3in  \frac{1}{\sqrt{2(1+b^2)}} \Bigg(\sqrt{2}b ,  \cos \Big( kx+ \frac{k^3}{\g^2}y\Big) , \sin \Big( kx+ \frac{k^3}{\g^2}y\Big), \cosh\Big( kx-\frac{k^3}{\g^2}y\Big),\sinh \Big( kx- \frac{k^3}{\g^2}y\Big)  \Bigg) \end{aligned} \end{equation}
with $k=\sqrt[\uproot{2} 4]{(1+b^2)\g^2},\;  b,\g>0$;

\item[{\rm (6)}]  A non-minimal flat surface given by 
\begin{equation}\begin{aligned}\notag &\hskip -.5in \Bigg(\frac{b\v }{\sqrt{\d^2+(1+ b^2)\v^2}}, 
 \frac{\sqrt{ \sqrt{1+b^2}(\d^2+ \v^2)-b\d \sqrt{\d^2+\v^2}} }{\sqrt{2}\sqrt[\uproot{2}4]{1+b^2}\sqrt{\d^2+(1+ b^2)\v^2}}\cos\! \big(\lambda (\sqrt{1+b^2}x+ \sqrt{\d^2+\v^2}y\big) ,\\& \hskip.1in   \frac{\sqrt{\! \sqrt{1+b^2}(\d^2+\v^2)- b\d \sqrt{\d^2+\v^2}}}{\sqrt{2}\sqrt[\uproot{2}4]{1\!+\!b^2}\sqrt{\d^2+(1+ b^2)\v^2}}  \sin \!\big(\lambda (\sqrt{1+b^2}x+\sqrt{\d^2+\v^2}y\big) ,\\& \hskip.3in  \frac{\sqrt{\! \sqrt{1+b^2}(\d^2+\v^2)+b\d \sqrt{\d^2+\v^2}} }{\sqrt{2}\sqrt[\uproot{2}4]{1+ b^2}\sqrt{\d^2+(1+b^2)\v^2}} \cosh\! \big(\mu (\sqrt{1+b^2}x-\sqrt{\d^2+\v^2}y\big),\\&\hskip.4in   \frac{\sqrt{\! \sqrt{1+b^2}(\d^2+\v^2)+b\d \sqrt{\d^2+\v^2}} }{\sqrt{2}\sqrt[\uproot{2}4]{1\!+\!b^2}\sqrt{\d^2\!+\!(1\!+\! b^2)\v^2}}  \sinh\! \big(\mu (\sqrt{1+b^2}x-\sqrt{\d^2+\v^2}y\big)\! \Bigg) \end{aligned} \end{equation}
with $\d,\v\ne 0, b>0$ and
\begin{align}\notag &\lambda =\frac{\sqrt{\sqrt{1+b^2}\sqrt{\d^2+\v^2}+b \d}}{\sqrt{\d^2+\v^2}},\, \;\; \mu= \frac{\sqrt{\sqrt{1+b^2}\sqrt{\d^2+\v^2}-b \d}}{\sqrt{\d^2+\v^2}}. \end{align}

\end{itemize}
Conversely, every parallel Lorentzian surface in $H^4_3(-1)$ is congruent to an open portion of one of the  six families of surfaces described above.
\end{Theorem}

\section{Parallel Lorentz surfaces in ${\mathbb C}^2_1$, $CP^2_1$ and $CH^2_1$}\label{sec19}

\subsection{Hopf fibrations}

Let $\mathbb C^n=\{(z_1,\ldots,z_n): z_1,\ldots, z_n\in {\bf C}\}$ be the complex $n$-space. If $\mathbb C^n$ endows with the metric given by the real part of the Hermitian form
\begin{equation} 
b_{j,n}((z_1,\ldots ,z_n),(w_1,\ldots ,w_n)) = -\sum_{k=1}^j \bar
z_k w_k + \sum_{k=j+1}^n \bar z_k w_k,
\end{equation}
then we obtain a flat indefinite Kaehler manifold of complex index $j$, denoted by $\mathbb C^n_j$. In particular, $\mathbb C^n_1$ is a flat Lorentzian Kaehler manifold. 

For any real number $c>0$, the differentiable manifold
\begin{equation} \label{1.4} S^{2n+1}_2(c) = \{ z\in \mathbb C^{n+1}_1  : b_{1,n+1}(z,z)=1/c \},\end{equation}
with the induced metric, is an indefinite real space form of constant sectional curvature $c>0$. 
The {\it Hopf-fibration}: $\pi : S^{2n+1}_2(c) \to CP^n_1(4c) : z \mapsto z\cdot \mathbb C^{\mathbf{\ast}}$ with $ \mathbb C^{\mathbf{\ast}}={\mathbb C}\setminus \{0\}$ is a submersion and there is a unique Lorentzian Kaehler metric
on $CP^n_1(4c)$ such that $\pi$ is a Riemannian submersion. The space $CP^n_1(4c)$ equipped with this metric is a Lorentzian Kaehler manifold of positive holomorphic sectional curvature $4c$. 

Similarly,  for any real number $c<0$, the differentiable manifold
\begin{equation} \label{1.6}
H^{2n+1}_3(c) = \{ z\in \mathbb C^{n+1}_2 : b_{2,n+1}(z,z)=1/c \},
\end{equation}
with the induced metric, is an indefinite real space form of constant sectional curvature $c<0$. The Hopf-fibration: $\pi : H^{2n+1}_3(c) \to \mathbb C H^n_1(4c) : z \mapsto
z\cdot \mathbb C^{\mathbf{\ast}}$
is a submersion and there is a unique Lorentzian Kaehler metric
on $CH^n_1(4c)$ such that $\pi$ is a Riemannian submersion. The
space $CH^n_1(4c)$ equipped with this metric is a Lorentzian
Kaehler manifold of negative holomorphic sectional curvature $4c$.

The manifolds $\mathbb C^n_1$, $CP^n_1(4c)$ and $CH^n_1(4c)$ are called
\emph{complex Lorentzian space forms}. The Riemann curvature tensor of a complex Lorentzian space form of constant
holomorphic sectional curvature $4c$ takes the form
\begin{equation} \label{1.8}
\tilde R(X,Y) = c(X \wedge Y + JX \wedge JY - 2\<JX,Y\>J),
\end{equation}
where $X$ and $Y$ are arbitrary tangent vectors at an arbitrary
point and $\wedge$ is defined by $$(X \wedge Y)Z=\<Y,Z\>X-\<X,Z\>Y.$$

\begin{Remark}\label{rem9}
The mapping $$\psi: \mathbb C^3_1 \to \mathbb C^3_2 : (z_1,z_2,z_3) \mapsto
(z_3,z_2,z_1)$$ maps $S^5_2(c)$ to $H^5_3(-c)$ and, via the
Hopf-fibrations, it induces a conformal mapping with factor $-1$
between $CP^2_1(4c)$ and $CH^2_1(-4c)$.
\end{Remark}

\subsection{Classification of parallel Lorentzian surface in ${\mathbb C}^2_1$}\label{subsec19.1}

For parallel Lorentzian surface in ${\mathbb C}^2_1$, we have the following result from 
\cite{CDV10} by Chen, Dillen and Van der Veken.

\begin{Theorem} \label{T44}
A parallel Lorentzian surface $M$ 
in ${\mathbb C}^2_1$ is isometric to an open part of one of the following
nine types of surfaces:
\begin{itemize}
\item[{\rm (1)}] a Lorentzian totally geodesic surface;
\item[{\rm (2)}] a Lorentzian product of parallel curves;
\item[{\rm (3)}] a complex circle, given by
$(a+ib)\big(\cos(x+iy),\sin(x+iy)\big)$ with $a,b\in {\bf R},\; (a,b)\neq(0,0);$
\item[{\rm (4)}] a $B$-scroll over the null cubic in $\mathbb E^3_1\subseteq \mathbb C^2_1$;
\item[{\rm (5)}] a $B$-scroll over the null cubic in $\mathbb E^3_2\subseteq \mathbb C^2_1$;
\item[{\rm (6)}] a surface given by
$$\frac{e^{-iy}}{\sqrt{2}}\big(i(1+a)-x-ay, i(1-a)+x+ay\big),\ \mbox{with }a\in {\bf R};$$
\item[{\rm (7)}] a surface with light-like mean curvature vector given by
$(q(x,y),x,y,q(x,y))$
with $q(x,y)=ax^2+bxy+cy^2+dx+ey+f$ and $a,b,c,d,e,f\in {\bf R}$;
\item[{\rm (8)}] a totally umbilical de Sitter space $S^2_1$ in $\mathbb E^3_1\subseteq \mathbb C^2_1$, given by
$a(0,\sinh x, \cosh x \cos y, \cosh x \sin y)$ with $a\in {\bf R}\setminus\{0\};$
\item[{\rm (9)}] a totally umbilical anti-de Sitter space $H^2_1$ in
$\mathbb E^3_2\subseteq {\mathbb C}^2_1$ given by
$a(\sin x, \cos x \cosh y, \cos x \sinh y, 0)$ with $a\in {\bf R}\setminus\{0\}.$
\end{itemize}
Conversely, each of the surfaces listed above is a Lorentzian
surface with parallel second fundamental form in $\mathbb C^2_1$.
\end{Theorem}

\subsection{Classification of parallel Lorentzian surface in $CP^2_1$}\label{subsec19.3}

First we mention the following result from \cite{CDV10}.

\begin{Lemma}\label{L5} Every parallel  Lorentzian surface in $CP^2_1(4)$ and in $CH^2_1(-4)$ is Lagrangian.
\end{Lemma}

The next classification of parallel Lorentzian surface in $CP^2_1$ was obtained by Chen, Dillen and Van der Veken in \cite{CDV10}.

\begin{Theorem} \label{T45}
Let $M$ be a Lorentzian surface in $CP^2_1(4)$ with parallel
second fundamental form. Then there are two possibilities:

\noindent\emph{(I)} $M$ is an open part of the totally geodesic,
Lagrangian surface $RP^2_1(1)\subseteq CP^2_1(4)$.

\noindent\emph{(II)} $M$ is flat, and the immersion is congruent
to $\pi\circ L$, where $\pi:S^5_2(1)\to CP^2_1(4)$ is the
Hopf-fibration and $L:M^2_1\to S^5_2(1)\subseteq\mathbb C^3_1$ is locally
one of the following twelve maps:

\begin{itemize}
\item[{\rm (1)}]  $L= \displaystyle{ \frac{1}{\sqrt 3} \( \sqrt 2 e^{\frac i2 x}\sinh\bigg(\frac{\sqrt 3}{2}y\Bigg),\ \sqrt 2 e^{\frac i2 x}\cosh\Bigg(\frac{\sqrt 3}{2}y\Bigg),\ e^{-ix} \) };$

\item[{\rm (2)}]  $L=  \displaystyle{ \( \frac{e^{\frac i2 (2x+y+\sqrt{1+4a}y)}}{(1+4a)^{1/4}},\ \frac{e^{\frac i2 (2x+y-\sqrt{1+4a}y)}}{(1+4a)^{1/4}},\ e^{iy} \) },\;\; a>-\frac 14;$

\item[{\rm (3)}]  $L=  \displaystyle{ \Bigg( \frac{(2-ie^{-\sqrt{4a-1}y})e^{ix+\frac 12(i+\sqrt{4a-1})y}}{2\sqrt[4]{4a-1}},  }
 \displaystyle{ \frac{(2+ie^{-\sqrt{4a-1}y})e^{ix+\frac 12(i+\sqrt{4a-1})y}}{2\sqrt[4]{4a-1}},\ e^{iy} \Bigg) },\;\; a>\frac 14;$

\item[{\rm (4)}]  $ L=  \displaystyle{ \frac{1}{\sqrt 2} \( e^{i(x+\frac y2)}(1+iy),\ e^{i(x+\frac y2)}(1-iy),\ \sqrt 2 e^{iy} \)};$

\item[{\rm (5)}]  $ L= \displaystyle{ \Bigg( \frac{\sqrt{a(2-a-b)}\,e^{i(bx+\frac{(1-b)y}{a(2-a-b)})}}{\sqrt{(a-b)(a+2b-2)}},\ \frac{\sqrt{b(2-a-b)}e^{i(ax+\frac{(1-a)y}{b(2-a-b)})}}{\sqrt{(a-b)(2a+b-2)}}, } \displaystyle{\frac{\sqrt{ab}\,e^{i((2-a-b)x+\frac{a+b-1}{ab}y)}}{\sqrt{(2a+b-2)(a+2b-2)}} \Bigg) }$
 with $a>b>2-a-b>0$ or  $\;0>a>b>2-a-b;$
 
 \item[{\rm (6)}] $ L=  \displaystyle{ \Bigg( \frac{\sqrt{b(a+b-2)}e^{i(ax+\frac{(1-a)y}{b(2-a-b)})}}{\sqrt{(a-b)(2a+b-2)}}, \frac{\sqrt{a(a+b-2)}e^{i(bx+\frac{(1-b)y}{a(2-a-b)})}}{\sqrt{(a-b)(a+2b-2)}},}  \displaystyle{\frac{\sqrt{ab}\,e^{i((2-a-b)x+\frac{a+b-1}{ab}y)}}{\sqrt{(2a+b-2)(a+2b-2)}} \Bigg) }$ with $a>b>0$ and $a+b>2;$

 \item[{\rm (7)}] $  L= \displaystyle{ \Bigg( \frac{\sqrt{-ab}\,e^{i((2-a-b)x+\frac{a+b-1}{ab}y)}}{\sqrt{(2a+b-2)(a+2b-2)}},\ \frac{\sqrt{b(2-a-b)}\,e^{i(ax+\frac{(1-a)y}{b(2-a-b)})}}{\sqrt{(a-b)(2a+b-2)}},} \displaystyle{ \frac{\sqrt{a(a+b-2)}e^{i(bx+\frac{(1-b)y}{a(2-a-b)})}}{\sqrt{(a-b)(a+2b-2)}} \Bigg) },$ with $ a>0>b>2-a-b;$

 \item[{\rm (8)}] $  L= \displaystyle{ \(
               \(\frac{2i\sqrt{(2a-1)(1-a)}}{2-3a}+\frac{2a^2(a-1)x+(2a-1)y}{2a\sqrt{(2a-1)(1-a)}}\)e^{i(ax+\frac{y} {2a})}, \right. }$
            $$ \displaystyle{ \left.\frac{(2a^2(a-1)x+(2a-1)y)e^{i(ax+\frac{y}{2a})}}{2a\sqrt{(2a-1)(1-a)}},\ \frac{a\,e^{i(2(1-
               a)x+\frac{2a-1}{a^2}y)}}{3a-2} \) },
              \;a\in(\tfrac 12 ,1)\setminus\{\tfrac 23\};$$

\item[{\rm (9)}] $ L =  \displaystyle{ \( \frac{(2a^2(a-1)x+(2a-1)y)e^{i(ax+\frac{y}{2a})}}{2a\sqrt{(2a-1)(a-1)}} , \right. }\displaystyle{ \( \frac{2a^2(a-1)x+(2a-1)y}{2a\sqrt{(2a-1)(a-1)}}+\frac{2i\sqrt{(2a-1)(a-1)}}{3a-2} \) \times
                } $
             $$\times e^{i(ax+\frac{y}{2a})}, \displaystyle{ \left. \frac{a\,e^{i(2(1-a)x+\frac{2a-1}{a^2}y)}}{3a-2} \) },\; a\in{\bf R}\setminus([\tfrac 12,1]\cup\{0\}); $$
               
\item[{\rm (10)}] $ L = \displaystyle{ \frac{e^{\frac{i}{12}(8x+9y)}}{24} \( 1+(8x-9y)^2+432iy,\ 2(8x-9y+12i),  1-(8x-9y)^2-432iy \) };$

\item[{\rm (11)}] $ L= \displaystyle{ \(
\frac{\sqrt{1-a}\,e^{i(ax+\frac{(a^2-b^2-a)y}{2(a-1)(a^2+b^2)})}}
               {b\sqrt{2a-1}\sqrt{(3a-2)^2+b^2}}\Big(2b(1-2a)\cosh\big(bx+\frac{b(2a-1)y}{2(a-1)(a^2+b^2)}\big) \right. }$
            $$ \displaystyle{ +i(3a^2-b^2-2a)\sinh\big(bx+\frac{b(2a-1)y}{2(a-1)(a^2+b^2)}\big)\Big), }$$
             $$ \displaystyle{ \frac{\sqrt{1-a}\sqrt{a^2+b^2}\,e^{i(ax+\frac{(a^2-b^2-a)y}{2(a-1)(a^2+b^2)})}}{b\sqrt{2a-1}}\sinh\big(bx+\frac{b(2a-1)y}{2(a-1)(a^2+b^2)}\big), }$$
           $$ \displaystyle{ \left. \frac{\sqrt{a^2+b^2}\,e^{i(2(1-a)x+\frac{2a-1}{a^2+b^2}y)}}{\sqrt{(3a-2)^2+b^2}} \) }, \mbox{ with } a\in(\tfrac 12, 1) \mbox{ and } b\in{\bf R}\setminus\{0\} ;$$
           
\item[{\rm (12)}] $ L= \displaystyle{ \( \sqrt{\frac{a-1}{2a-1}}\frac{e^{i(ax+\frac{(a^2-b^2-a)y}{2(a-1)(a^2+b^2)})}}{b\sqrt{(3a-2)^2+b^2}}\Big(2b(1-2a)\sinh\big(bx+\frac{b(2a-1)y}{2(a-1)(a^2+b^2)}\big) \right. }$
            $$ \displaystyle{ +i(3a^2-b^2-2a)\cosh\big(bx+\frac{b(2a-1)y}{2(a-1)(a^2+b^2)}\big)\Big), }$$
                         $$ \displaystyle{ \sqrt{\frac{a-1}{2a-1}}\frac{\sqrt{a^2+b^2}\,e^{i(ax+\frac{(a^2-b^2-a)y}{2(a-1)(a^2+b^2)})}}{b}\cosh\big(bx+\frac{b(2a-1)y}{2(a-1)(a^2+b^2)}\big),}$$
        $$ \displaystyle{ \left. \frac{\sqrt{a^2+b^2}\,e^{i(2(1-a)x+\frac{2a-1}{a^2+b^2}y)}}{\sqrt{(3a-2)^2+b^2}} \) }, \mbox{ with } a\in {\bf R}\setminus[\tfrac 12, 1] \mbox{ and } b\in {\bf R}\setminus\{0\}.$$
\end{itemize}\end{Theorem}

\subsection{Classification of parallel Lorentzian surface in $CH^2_1$}\label{subsec19.4}

It follows from Remark \ref{rem9} that one obtains immediately the classification of parallel Lorentzian surfaces in $CH^2_1(-4)$ from Theorem \ref{T45} via the mapping:
$$\psi: \mathbb C^3_1 \to \mathbb C^3_2 : (z_1,z_2,z_3) \mapsto
(z_3,z_2,z_1)$$
since $\psi$ gives rise to a conformal mapping with factor $-1$ between $CP^2_1(4)$ and $CH^2_1(-4)$. Hence, besides totally geodesic Lagrangian surface $RH^2_1(-1)\subset CH^2_1(-4)$, there are twelve families of flat parallel Lorentzian surfaces in $CH^2_1(-4)$.

\section{Parallel surfaces in warped product $I\times_f R^n(c)$}\label{sec20}

\subsection{Basics on Robertson-Walker space-times}\label{subsec20.1}

In the theory of general relativity,  a {\it Robertson-Walker space-time} is a warped product
 \begin{align}\label{4.1} L^4_1(c,f)=(I\times R^3(c),g),\;\;   g=-dt^2 +f^2(t)g_c,\end{align}
  of an open interval  $I$ and  a  Riemannian 3-manifold $(R^3(c),g_c)$ of constant  curvature $c$, while the warping function $f$ describes the expanding or contracting of our Universe (cf. \cite{O,book17}).

A Robertson-Walker space-time possesses two relevant geometrical features. On one hand, its fibers have  constant curvature. Hence the space-time is {\it spatially homogeneous}. 
On the other hand, it has a time-like vector field $K=f(t)\partial_t$ which satisfies $\nabla_XK=f'(t)X$ for any $X$. In particular, we have  ${\mathcal L}_K  g=2 f'  g$,  where ${\mathcal L}_K$ is the Lie derivative along $K$. Hence the canonical time-like vector field $K$  is a conformal vector field. These properties of $K$ show a certain symmetry on $L^n_1(c,f)$. 

One may also consider a higher dimensional Robertson-Walker space-time as 
\begin{align}\label{4.2} L^n_1(c,f):=(I\times R^{n-1}(c),g),\;\;   g=-dt^2 +f^2(t)g_c,\end{align}
where $R^{n-1}(c)$ is a Riemannian $(n-1)$-manifold of constant  curvature $c$ for $n> 5$. 

A {\it rest space} or a {\it space-like slice} in $L^n_1(c,f)$ is a space-like hypersurface given by $t$ constant. Thus a rest space in $L^{n}_1(c,f)$ is a fiber
$$S(t_0)=\{t_0\}\times_{f(t_0)} R^{n}(c),\;\; t_0\in I.$$ Hence a rest space
$S(t_0)$ in $L^{n}_1(c,f)$ is an $(n-1)$-manifold of constant curvature whose metric tensor is  given by $f^2(t_0)g_k$.

 A  pseudo-Riemannian submanifold $N$ of a Robertson-Walker space-time $L^n_1(c,f)$ is called  {\it transverse} if it is contained in a {\it rest space} $S(t_0)$ for some $t_0\in I$.  
A pseudo-Riemannian  submanifold $N$ of $L^n_1(c,f)$ is called a $\mathcal H$-{\it submanifold} if the tangent field $\frac{\partial}{\partial t}$, known as the {\it comoving observer field},  is tangent to $N$ at each point on $N$.

\subsection{Parallel submanifolds of Robertson-Walker space-times}\label{subsec20.2}

For parallel submanifolds of $L^n_1(c,f)$, we have the next classification result from \cite{ChenV13.1,book17}.
 
 \begin{Theorem}\label{T:46} If a Robertson-Walker space-time $L_1^n(c,f)$ does not contain any open subsets of constant curvature, then  a  $k$-dimensional pseudo-Riemannian submanifold of $L_1^n(c,f)$ is a parallel submanifold if and only if it is one of the following:

\begin{itemize}
  \item[{\rm (a)}]  
  A transverse submanifold lying in a rest space $S(t_0)$ of $L_1^n(c,f)$ as a parallel submanifold.

\item[{\rm (b)}] An $\mathcal H$-submanifold which is locally a warped product $ I\times_{ f} P^{k-1}$, where  $I$ is an open interval and $P^{k-1}$  is a submanifold of  $R^{n-1}(c)$. Further, 
\vskip.03in

 \hskip-.2in {\rm (b.1)} if $f'\ne 0$ on $I$, then $I \times_f P^{k-1}$ is totally geodesic in $L_1^m(k,f)$;
\vskip.03in

 \hskip-.2in {\rm (b.2)} if $f'= 0$ on $I$, then $P^{k-1}$ is a parallel submanifold of  $R^{n-1}(c)$.
\end{itemize}
\end{Theorem}
 
Similar result holds for submanifolds in a  warped product $I\times_f R^{n-1}(c)$ with the Riemannian warped product metric $g=dt^2+f^2(t)g_c$ {\rm (cf. \cite{CW08,VV08})}.

\section{Thurston's eight three-dimensional model geometries}\label{sec21}

The uniformization theorem for 2-dimensional surfaces says that every simply-connected Riemann surface is conformally equivalent to one of the three Riemann surfaces: the open unit disk, the complex plane, or the Riemann sphere. This result implies that every Riemann surface admits a Riemannian metric of constant curvature. 

 Roughly speaking,  for closed 3-manifolds W. Thurston's Geometrization Conjecture states that every closed 3-manifold can be decomposed in a canonical way into pieces that each have one of eight types of geometric structure  locally (see, \cite{Thu97}). In 2005, G. Perelman \cite{P03} provided a proof of Thurston's  geometrization conjecture via Ricci flow with surgery.

\vskip.05in
The  eight Thurston's 3-dimensional model geometries are the following.

\vskip.05in
(1) {\it Euclidean geometry $\mathbb E^3$}. 
\vskip.05in

(2) {\it Spherical geometry $S^3$}. 

\vskip.05in
(3) {\it Hyperbolic geometry $H^3$}. 

\vskip.05in
(4) {\it The geometry of $S^2\times {\mathbb R}$}. 

\vskip.05in
(5) {\it The geometry of $H^2\times {\mathbb R}$}. 

\vskip.05in
(6) {\it The geometry  $\widetilde{SL_2}({\mathbb R})$}. The 3-dimensional Lie group of all $2\times 2$ real matrices with determinant one is denoted by $SL_2({\mathbb R})$; and $\widetilde{SL_2}({\mathbb R})$ denotes its universal covering. $\widetilde{SL_2}({\mathbb R})$ is a unimodular Lie group with a special left invariant metric.
Examples of these manifolds in this geometry include the manifold of unit vectors of the tangent bundle of a hyperbolic surface and, more generally, the Brieskorn homology spheres.

\vskip.05in
(7) {\it Nil geometry $Nil_3$}. The group $Nil_3$ is a 3-dimensional unimodular Lie group with a special left invariant metric consisting of real matrices of the form
\begin{align}\notag \begin{pmatrix} 1\; & x\; & y\\ 0&1&z\\0&0&1\end{pmatrix}\end{align}
under multiplication. This group, also known as the {\it Heisenberg group}, is nilpotent. 

\vskip.05in

(8) {\it Sol geometry $Sol_3$}. This group $Sol_3$ has the least symmetry of all the eight
geometries as the identity component of the stabilizer of a point is trivial.

\vskip.1in

We mentioned earlier in \S1 that the complete classification of parallel surfaces in $\mathbb E^3$ was obtained by V. F. Kagan;  the complete classifications of parallel surfaces in 
 $S^3$ and in $H^3$ were given in  \S5.4 and \S5.5, respectively; the classifications of parallel surfaces in $S^2\times \mathbb R$ and in $H^2\times \mathbb R^3$ were given in \S20.
 
 In this section, we will deal the classification of parallel surfaces in $Sol_3, \widetilde{SL_2}({\mathbb R})$ and $Nil_3$ in \S22.2, \S22.4 and \S22.5, respectively.

\section{Parallel surfaces in three-dimensional Lie groups}\label{sec22}

\subsection{Milnor's classification of 3-dimensional unimodular Lie groups}\label{subsec22.1}

A Lie group $G$ is called {\it unimodular} if its left-invariant Haar measure is also right-invariant. In \cite{Milnor76}, J. Milnor provides an infinitesimal reformulation of unimodularity for 3-dimensional Lie groups. We recall it briefly as follows:

Let ${\mathfrak g}$ be a 3-dimensional oriented Lie algebra equipped with an inner product $\<\;,\>$. Define the vector product operation $\times :{\mathfrak g}\times {\mathfrak g}\to {\mathfrak g}$ as the skew-symmetric bilinear map which is
uniquely determined by the following three conditions:

\begin{itemize}
\item[{\rm (a)}] $\<X,X\times Y\>=\<Y,X\times Y\>=0$,

\item[{\rm (b)}] $|X\times Y|^2=\<X,X\>\<Y,Y\>-\<X,Y\>^2$,

\item[{\rm (c)}]  if $X$ and $Y$ are linearly independent, then $\det (X,Y,X\times Y)>0,$
\end{itemize}
for all $X, Y\in {\mathfrak g}$. The Lie-bracket $[\;\cdot\;,\;\cdot\;]$ on $ {\mathfrak g}$ is a skew-symmetric bilinear map. By comparing these two operations, one obtains a linear endomorphism $L_{\mathfrak g}$ which is uniquely determined by the formula
$$[X, Y] = L_{\mathfrak g}(X \times Y ),\;\; X, Y\in {\mathfrak g}.$$

If $G$ is an oriented 3-dimensional Lie group equipped with a left-invariant Riemannian metric, then the metric induces an inner product on the Lie algebra ${\mathfrak g}$. With
respect to the orientation on ${\mathfrak g}$ induced from $G$, the endomorphism field $L_{\mathfrak g}$  is uniquely
determined.  

J. Milnor proved in \cite{Milnor76} that the unimodularity of $G$ is characterized as follows.

\begin{Theorem}\label{T47} Let $G$ be an oriented 3-dimensional Lie group with a left-invariant
Riemannian metric. Then $G$ is unimodular if and only if the endomorphism $L_{\mathfrak g}$ is self-adjoint with respect to the metric.
\end{Theorem}

If $G$ is a 3-dimensional unimodular Lie group with a left-invariant metric, then there exists an orthonormal basis $\{e_1, e_2, e_3\}$ of the Lie algebra ${\mathfrak g}$ such that
$$[e_1, e_2] = c_3e_3,\;\; [e_2, e_3] = c_1e_1,\;\; [e_3, e_1] = c_2e_2,\;\;  c_i\in {\mathbb R}.$$

Milnor obtained the following classification of 3-dimensional unimodular Lie groups.

\vskip.05in
\begin{center}
\begin{tabular}{|c|c|c|}
\hline
$(c_1,c_2,c_3)$ & Simply-connected Lie group & Property \\ \hline
$(+,+,+)$ & $SU(2)$ & Compact and simple \\ \hline
$(+,+,-)$ & $\widetilde{SL}(2,\mathbb R)$ & Non-compact and simple \\ \hline
$(+,+,0)$ & $\tilde{E}(2)$ & Solvable \\ \hline
$(+,-,0)$ & $E(1,1)$ & Solvable \\ \hline
$(+,0,0)$ & Heisenberg group & Nilpotent \\ \hline
$(0,0,0)$ & $(\mathbb E^3,+)$ & Abelian \\ \hline
\end{tabular}
\end{center}
$$\hbox{Three-dimensional unimodular Lie groups classified by J. Milnor}$$

\vskip.05in
\noindent Here $E(1, 1)$ denotes the the group of orientation-preserving rigid motions of Minkowski plane, $E(2)$ denotes the group of orientation-preserving rigid motions of Euclidean plane and $\tilde{E}(2)$ is the universal covering of $E(2)$.

\subsection{Parallel  surfaces in  the motion group $E(1,1)$}\label{subsec22.2}

Let $E(1,1)$ be the motion group of the Minkowski plane:
$$ E(1,1)=\left\{ \left(\begin{array}{ccc} e^{z} & 0 & x\\ 0 & e^{-z} & y\\ 0 & 0 & 1
\end{array}\right)   : x,y,z \in \mathbb{R} \right\}.
$$
The Lie algebra $\mathfrak{e}(1,1)$ is given by
$\mathfrak{e}(1,1)= \left\{ \left(\begin{array}{ccc} w & 0 & u\\ 0 & -w & v\\ 0 & 0 & 0
\end{array}\right)  : u,v,w \in \mathbb{R} \right\}.$
Consider the basis
$$F_1=\left(\begin{array}{ccc} 0 & 0 & 1\\ 0 & 0 & 0\\ 0 & 0 &0\end{array}
\right), \quad 
F_2= \left(\begin{array}{ccc} 0 & 0 & 0\\ 0 & 0 & 1\\ 0 & 0 &0\end{array}\right), \quad 
F_3= \left(\begin{array}{ccc} 1 & 0 & 0\\ 0 & -1 & 0\\ 0 & 0 &0 \end{array} \right) $$
of $\mathfrak{e}(1,1)$. Then the left-translated vector fields of
$\{F_1,F_2,F_3\}$ are given by
$$ f_1=e^{z}\frac{\partial}{\partial x},\quad f_2=e^{-z}\frac{\partial}{\partial y},\quad
f_3=\frac{\partial}{\partial z}. $$
The dual coframe field is
$\omega^1=e^{-z}dx,\, \omega^2=e^{z}dy,\, \omega^3=dz.$
Now we take the following left-invariant vector fields $u_1,u_2,u_3$:
$$ u_1=\frac{1}{\sqrt{2}}(-f_1+f_2),\quad u_2=\frac{1}{\sqrt{2}}(f_1+f_2),\quad u_3=f_3. $$
This left-invariant frame field satisfies the  relations
$ [u_1,u_2]=0,\, [u_2,u_3]=u_1,\, [u_3,u_1]=-u_2. $
We equip $E(1,1)$ with a left-invariant Riemannian metric such
that $\{e_1,e_2,e_3\}$, with $e_i={u_i}/{\lambda_i}$, is
orthonormal, where $\lambda_1,\lambda_2,\lambda_3$ are positive
constants. The resulting Riemannian metric is
$$ g_{(\lambda_1,\lambda_2,\lambda_3)} =\frac{\lambda_1^2}{2}(-\omega^{1}+\omega^{2})^{2}
+\frac{\lambda_2^2}{2}(\omega^{1}+\omega^{2})^{2} +\lambda_{3}^{2}(\omega^{3})^2. $$

V. Patrangenaru proved the following result in \cite{Pa}.

\begin{Theorem}\label{T48} A left-invariant metric on $E(1,1)$ is isometric to one of the
metrics $g_{(\lambda_1,\lambda_2,\lambda_3)}$ with
$\lambda_{1}\geq \lambda_{2}>0$ and $\lambda_{3}=\frac{1}{\lambda_{1}\lambda_{2}}$.
\end{Theorem}

\subsection{Parallel surfaces in the Heisenberg group $Sol_3$}\label{subsec22.3}

If we put $g(\lambda_1,\lambda_2)=g_{(\lambda_1,\lambda_2,\frac{1}{\lambda_1\lambda_2})}$, then the homogeneous 3-manifold $Sol_3=(E(1,1),g_{(1,1)})$ is one of Thurston's eight model spaces. Hence, $Sol_3$ has a natural 2-parametric
deformation family
$\{(E(1,1),g(\lambda_1,\lambda_2)) \ \vert \ \lambda_{1}\geq \lambda_{2}>0 \}.$

In \cite{IV07}, J. Inoguchi and J. Van der Veken classified parallel surfaces in  $Sol_3=(E(1,1),g(\lambda_1,\lambda_2))$ as follows. 

\begin{Theorem}\label{T49} Let $M$ be a parallel surface in $Sol_3=(E(1,1),g(\lambda_1,\lambda_2))$. Then $M$ is one of the following:
\begin{itemize}
\item[{\rm (a)}] an integral surface of the distribution spanned by $\left\{{\partial}/{\partial x},{\partial}/{\partial y}\right\}$,
\item[{\rm (b)}] an integral surface of the distribution spanned by $\left\{{\partial}/{\partial x},{\partial}/{\partial z}\right\}$ or $\left\{{\partial}/{\partial y},{\partial}/{\partial z}\right\}$, \end{itemize}
the latter case only occurring if $\lambda_1=\lambda_2$. 
Moreover, the surfaces described in $\mathrm{(a)}$ are flat and minimal, but not
totally geodesic and the surfaces in $\mathrm{(b)}$ are totally
geodesic and have constant Gaussian curvature $-\lambda_1^4$.
\end{Theorem}

\subsection{Parallel  surfaces in  the motion group $E(2)$}\label{subsec22.4}

The Euclidean motion group $E(2)$ is given  by the following matrix group:
$$ E(2)= \left\{ \left(\begin{array}{ccc}\cos \theta & -\sin \theta & x
\\ \sin \theta & \cos \theta & y\\ 0 & 0 & 1\end{array}
\right) : x,y \in \mathbb{R},\ \theta\in S^1\right\}.$$
The universal covering group of $E(2)$ is $\mathbb{R}^3$ with multiplication
$$(x,y,z)\cdot (x^\prime,y^\prime,z^\prime)= (x+x^\prime\cos z-y^\prime\sin z\ , y+x^\prime\sin z +y^\prime\cos z\ , z+z^\prime).$$
Take positive constants $\lambda_1,\lambda_2$ and $\lambda_3$ and
a left-invariant frame
$$e_1= \frac{1}{\lambda_2} \left(\! -\sin z\frac{\partial}{\partial x}+ \cos z \frac{\partial}{\partial y} \right), \quad 
e_2=\frac{1}{\lambda_3} \frac{\partial}{\partial z}, \quad 
e_3=\frac{1}{\lambda_1} \left (\cos z\frac{\partial}{\partial x}+ \sin z \frac{\partial}{\partial y} \right).
$$
Then this frame satisfies the commutation relations:
$[e_1,e_2]=c_{1}e_3,\, [e_2,e_3]=c_{2}e_1,\, [e_3,e_1]=0, $
with $c_1=\frac{\lambda_1}{\lambda_2\lambda_3}$ and
$c_2=\frac{\lambda_2}{\lambda_1\lambda_3}$. 
The left-invariant
Riemannian metric determined by the condition that $\{e_1,e_2,e_3\}$ is orthonormal is given by
$$g_{(\lambda_1,\lambda_2,\lambda_3)}=\lambda_1^2(\cos z\,dx+\sin z\,dy)^2+\lambda_2^2(-\sin z\,dx+\cos z\,dy)^2+\lambda_3^2\,dz^2.$$

We have the following result on $\widetilde{E(2)}$ from \cite{Pa}.

\begin{Proposition}\label{P1} A left-invariant metric on $\widetilde{E(2)}$ is isometric to
one of the metrics $g_{(\lambda_1,\lambda_2,\lambda_3)}$ with $\lambda_{1} > \lambda_{2}>0$ and $\lambda_{3}=\frac{1}{\lambda_{1}\lambda_{2}}$, or $\lambda_1=\lambda_2=\lambda_3=1$. 
In particular, $\widetilde{E(2)}$ with metric $g_{(1,1,1)}$ is isometric to Euclidean $3$-space $\mathbb{E}^{3}$.
\end{Proposition}

J. Inoguchi and J. Van der Veken classified parallel surfaces in $\widetilde{E(2)}$  in \cite{IV07} as follows. 

\begin{Theorem}\label{T50} The only parallel surfaces in $\widetilde{E(2)}$ are integral surfaces of the distribution spanned by
$\left\{{\partial}/{\partial x},{\partial}/{\partial y}\right\}$. 
These surfaces are flat and minimal, but not totally geodesic.
\end{Theorem}

\subsection{Parallel  surfaces in  $SU(2)$}\label{subsec22.5}

The group ${SU}(2)$ is diffeomorphic to $S^3$, since
$${SU}(2)= \left \{ \ \left (\begin{array}{cc} x_0+\sqrt{-1}x_3 & -x_2+\sqrt{-1}x_1
\\x_2+\sqrt{-1}x_1 & x_0-\sqrt{-1}x_3\end{array}\right ) : x_0^2+x_1^2 +x_2^2+x_3^2=1 \right\}.$$

From \cite{Pa},  we have the following proposition which describes all possible left-invariant metrics on ${SU}(2)$.

\begin{Proposition} \label{prop3} Any left-invariant metric on ${SU}(2)$ is isometric to one of
the following metrics $g_{(\lambda_1,\lambda_2,\lambda_3)}$ with
$\lambda_i\in\mathbb R$ and $\lambda_{1}\geq \lambda_{2}\geq
\lambda_{3}>0$:$$
g_{(\lambda_1,\lambda_2,\lambda_3)}=
\frac{4}{\lambda_{2}\lambda_{3}}\sigma_{1}^{2}+
\frac{4}{\lambda_{3}\lambda_{1}}\sigma_{2}^{2}+
\frac{4}{\lambda_{1}\lambda_{2}}\sigma_{3}^{2},
$$
on the unit three-sphere $S^{3}(1)=\{ (x_0,x_1,x_2,x_3)\in \mathbb{E}^{4} :
x_0^{2}+x_1^{2}+x_2^{2}+x_3^{2}=1\},$
where
\begin{align*}&\sigma_{1}= -x_1dx_0+x_0dx_1-x_3dx_2+x_2dx_3,\\
&\sigma_{2}= -x_2dx_0+x_3dx_1-x_0dx_2+x_1dx_3,\\
&\sigma_{3}= -x_3dx_0+x_2dx_1-x_1dx_2+x_0dx_3.
\end{align*}
The dimension $d(\lambda_1,\lambda_2,\lambda_3)$ of the isometry
group of $({SU}(2), g_{(\lambda_1,\lambda_2,\lambda_3)})$
is
$$d(\lambda_1,\lambda_2,\lambda_3)= \left\{\begin{array}{ccc}
3 & \textrm{if} & \lambda_{1}>\lambda_{2}>\lambda_{3},\\
4 & \textrm{if} & \lambda_{1}=\lambda_{2}>\lambda_{3}\ \textrm{or}
\ \lambda_{1}>\lambda_{2}=\lambda_{3},\\ 6 & \textrm{if} & \lambda_{1}=\lambda_{2}=\lambda_{3}.\end{array}\right.$$
\end{Proposition}

Let $\mathfrak{su}(2)$ denote the Lie algebra of ${SU}(2)$. Take the following quaternionic basis $\{i,j,k\}$ of $\mathfrak{su}(2)$:
$$i= \left (\begin{array}{cc} 0 & \sqrt{-1} \\\sqrt{-1} & 0 \end{array}\right ), \ \
 j= \left (\begin{array}{cc} 0 & -1 \\ 1 & 0 \end{array} \right ), \ \ 
 k= \left (\begin{array}{cc}\sqrt{-1} & 0 \\ 0 & -\sqrt{-1}\end{array}\right ).$$
We denote the left-translated vector fields of $i,j,k$ by $E_1,E_2,E_3$. Then the commutation relations of$\{E_1,E_2,E_3\}$ are given by
$[E_1,E_2]=2E_3,\, [E_2,E_3]=2E_1,\, [E_3,E_1]=2E_2. $
Choose strictly positive real constants $\lambda_1$, $\lambda_2$,
$\lambda_3$ and define
$$e_{1}=\frac{1}{\lambda_{2}\lambda_{3}} i,\ e_{2}=\frac{1}{\lambda_{3}\lambda_{1}} j,\ 
e_{3}=\frac{1}{\lambda_{1}\lambda_{2}} k. $$
Then $ [e_1,e_2]=c_{3}e_{3},\, [e_2,e_3]=c_{1}e_{1},\, [e_3,e_1]=c_{2}e_{2}, $
with
$c_{1}={2}/{\lambda_{1}^{2}},c_{2}={2}/{\lambda_{2}^{2}},c_{3}={2}/{\lambda_{3}^{2}}.$
The left-invariant metric $\overline{g}_{(c_1,c_2,c_3)}$, defined
by the condition that $\{e_1,e_2,e_3\}$ is an orthonormal basis,
is
$$
\overline{g}_{(c_1,c_2,c_3)}=4\left(\frac{1}{c_2c_3}\omega_{1}^{2}+
\frac{1}{c_1c_3}\omega_{2}^{2}+
\frac{1}{c_1c_2}\omega_{3}^{2}\right),
$$
where $\{\omega_1,\omega_2,\omega_3\}$ is the dual coframe field
of
$\{E_1,E_2,E_3\}$.

The following result from \cite{Pa} describes all left-invariant metrics on $SU(2)$

\begin{Proposition}\label{P3} A left-invariant metric on ${SU}(2)$ is isometric to one of the metrics $\overline{g}_{(c_1,c_2,c_3)}$, with $c_1,c_2,c_3\geq
0$. Moreover, the dimension of the isometry group is $\geq 4$ if and
only if at least two of the parameters $c_i$ coincide.
\end{Proposition}

The next non-existence result was proved  by J. Inoguchi and J. Van der Veken in \cite{IVan08}.

\begin{Theorem}\label{T51}  There are no parallel surfaces in $SU(2)$ equipped with a left-invariant metric
with 3-dimensional isometry group.
\end{Theorem}

\subsection{Parallel surfaces in the real special linear group ${SL}(2,\mathbb R)$}\label{subsec22.6}

The group $SL(2,\mathbb R)$ is defined as the following subgroup of $GL(2,\mathbb R)$:
$$SL(2,\mathbb  R) =\left\{\begin{pmatrix} a & b\\ c& d\end{pmatrix}: ad-bc=1\right\}.$$
This group is isomorphic to the following subgroup of $GL(2,\mathbb C)$:
$$SU(1, 1) =\left\{\begin{pmatrix} \alpha & \beta\\ \bar\beta & \bar\alpha\end{pmatrix}: |\alpha|-|\beta|^2=1\right\}$$
via the isomorphism
$SL(2,\mathbb R)\to SU(1, 1):\begin{pmatrix} a & b\\ c& d\end{pmatrix}\mapsto
\begin{pmatrix} i & 1\\ 1 &i\end{pmatrix}\begin{pmatrix} a & b\\ c& d\end{pmatrix}
\begin{pmatrix}- i & 1\\ 1 &-i\end{pmatrix}.$
The Lie algebra of $SU(1, 1)$ is explicitly given by
$$\mathfrak{s}\mathfrak{u}(1, 1) =\left\{\begin{pmatrix} iu & v-iw\\ v+iw & -iu\end{pmatrix}: u,v,w\in \mathbb R\right\}.$$
We take the following split-quaternionic basis of the Lie algebra su(1, 1):
$${\bf i}=\begin{pmatrix} i & 0\\ 0 &-i\end{pmatrix},\;\; 
{\bf j}'=\begin{pmatrix} 0 & -i\\ i &0\end{pmatrix}, \;\; {\bf k}'=\begin{pmatrix} 0 & 1\\ 1 &0\end{pmatrix}.$$
Denote the left-translated vector fields of $\{{\bf j}',{\bf k}',{\bf i}\}$ by $\{E_1, E_2, E_3\}$ and choose strictly positive real constants $\lambda_1,\lambda_2,\lambda_3$ and define
$$e_1=\frac{1}{\lambda_2\lambda_3}E_1,\; \; e_2=\frac{1}{\lambda_1\lambda_3}E_2,\;\; e_3=\frac{1}{\lambda_1\lambda_2}E_3.$$
Then we have
$[e_1,e_2]=c_3e_3,\, [e_2,e_2]=c_1e_1,\, [e_3,e_1]=c_23e_2$ with $c_1=2/\lambda_1^2, c_2=2/\lambda_2^2$ and $c_2=-2/\lambda_3^2$.
The left-invariant Riemannian metric $g{(c_1,c_2,c_3)}$  by the
condition that $\{e_1,e_2,e_3\}$ is an orthonormal basis is
\begin{align*}
g{(c_1,c_2,c_3)} =4 \(-\frac{1}{c_2c_3}\omega_1^2-\frac{1}{c_1c_3}\omega_2^2
+\frac{1}{c_1c_2}\omega_3^2\),\end{align*}
where $\{\omega_1,\omega_2,\omega_3\}$ is the dual coframe  field of $\{E_1.E_2,E_3\}$.
This three-parameter family of Riemannian metrics exhausts all
left-invariant metrics on ${SL}(2,\mathbb R)$ as shown in the next proposition from \cite{Pa}.

\begin{Proposition} \label{prop6} Any left-invariant metric on ${SU}(1,1)$ is isometric to one of the metrics ${g}{(c_1,c_2,c_3)}$ with $c_{3}<0<c_{2}\leq c_{1}$. Moreover, this metric gives rise to an isometry group of dimension $4$ if and only if $c_1=c_2$.
\end{Proposition}

We consider ${SL}(2,\mathbb R)$ equipped with a left-invariant metric such that the dimension of the isometry group is only $3$. With the notations given above, we have that $c_1>c_2>0>c_3$.

The following classification theorem for parallel surfaces in ${SL}(2,\mathbb R)$ was proved by J. Inoguchi and J. Van der Veken in \cite{IVan08}.

\begin{Theorem}\label{TY52} Consider ${SL}(2,\mathbb R)$, equipped with a left-invariant metric with $c_1>c_2>0>c_3$. Parallel surfaces only occur if $c_2=c_1+c_3$.
Moreover, they are integral surfaces of the distribution spanned by
$\{\cos\theta\, e_1+\sin\theta\,e_3,e_2\}$, where $\theta$ is a
constant, satisfying $\tan^2\theta=-{c_3}/{c_1}$. These surfaces
are totally geodesic and of constant Gaussian curvature given by $c_1c_3<0$.
\end{Theorem}

\subsection{Parallel surfaces in  non-unimodular three-dimensional Lie groups}\label{subsec22.7}

Let $G$ be a non-unimodular 3-dimensional Lie group with a left-invariant metric. Then the \textit{unimodular kernel} $\mathfrak{u}$ of the Lie algebra $\mathfrak{g}$ of $G$ is defined
by $\mathfrak{u}=\{X \in \mathfrak{g} : \mathrm{Tr}\; \mathrm{ad}(X)=0\},$
where $\mathrm{ad}:\mathfrak{g}\to \mathrm{End}(\mathfrak{g})$ is a homomorphism defined by
$\mathrm{ad}(X)Y=[X,Y].$ Then $\mathfrak{u}$ is an ideal of $\mathfrak{g}$ containing the ideal $[\mathfrak{g},\mathfrak{g}]$.

On $\mathfrak{g}$, we can take an orthonormal basis $\{e_1,e_2,e_3\}$ such that
(a) $\langle e_{1}, X\rangle=0,\, X \in \mathfrak{u}$ and
(b) $\langle [e_1,e_2], [e_1,e_3]\rangle=0$.
The commutation relations of this basis are given by
$$[e_1,e_2]=a e_{2}+b e_{3},\ [e_2,e_3]=0,\ [e_1,e_3]=c e_2+d e_3, $$
with $a+d\not=0$ and $ac+bd=0$. Under a suitable homothetic change
of the metric, we may assume that $a+d=2$. Then the constants $a$,
$b$, $c$ and $d$ are represented as
$$a=1+\xi,\ b=(1+\xi)\eta,\ c=-(1-\xi)\eta,\ d=1-\xi,$$
where $(\xi,\eta)$ satisfies the condition $\xi,\eta\geq 0$ and $\xi^2+\eta^2\not=0$. 

The next was also proved by J. Inoguchi and J. Van der Veken in \cite{IVan08}.

\begin{Proposition} \label{P5} The non-unimodular Lie group $G$ is locally symmetric if and only if $\xi=0$ or $(\xi,\eta)=(1,0)$.
\end{Proposition}

Since the parallel surfaces in ${H}^3$, ${H}^{2}(-4)\times \mathbb{R}$ and $\widetilde{{SL}}(2,\mathbb R)$ (with four-dimensional isometry group) are already classified, we shall restrict our attention to such surfaces in the non-unimodular Lie groups, satisfying $\xi\notin\{0,1\}.$

 The following theorem of J. Inoguchi and J. Van der Veken from  \cite{IVan08} provides the classification of parallel surfaces in the corresponding Lie groups.

\begin{Theorem}\label{T53}  Let $G$ be a non-unimodular Lie group with structure constants $(\xi,\eta)$. Assume that $\xi\notin\{0,1\}$. Then the only parallel surfaces in $G$ are:
\begin{itemize}
\item[{\rm (1)}] Integral surfaces of the distributions spanned by $\{e_1,e_2\}$, respectively
$\{e_1,e_3\}$. These surfaces are totally geodesic and of constant negative curvature $-(1+\xi)^2$ respectively $-(1-\xi)^2$.
\item[{\rm (2)}] Integral surfaces of the distribution spanned by $\{e_2,e_3\}$. These surfaces are flat and of constant mean curvature 1.
\end{itemize}
The former case only occurs when $\eta=0$.
\end{Theorem}

  \subsection{Parallel surfaces in the Heisenberg group $Nil_3$}\label{subsec22.8}

The following classification  theorem of parallel surfaces in the Heisenberg group $Nil_3$ was proved  by M. Belkhelfa, F. Dillen and J.  Inoguchi in \cite{BDI02}.

\begin{Theorem}\label{T54} The only parallel surfaces in the Heisenberg group $Nil_3$ are open parts of vertical planes and vertical round cylinders.
\end{Theorem}

\begin{Remark}  The oscillator group was introduced and first studied by R. F. Streater in \cite{St67} and owes its name to the fact that its Lie algebra coincides with the one generated by the differential operators associated to the harmonic oscillator problem. Generalizing this construction, oscillator groups
have been defined in any even dimension greater or equal to four. Since their
introduction, the oscillator groups have been intensively studied from several
different points of view, both in differential geometry and in mathematical
physics. Beside direct extensions with Euclidean
groups, the oscillator groups are the only simply connected non-Abelian solvable
Lie groups admitting a bi-invariant Lorentzian metric. 

In \cite{CV13.1} G. Calvaruso and J. Van der Veken obtained the complete classification and explicitly description of totally geodesic and parallel hypersurfaces of four-dimensional oscillator groups, equipped with a one-parameter family of left-invariant Lorentzian metrics.
\end{Remark}

\section{Parallel surfaces in three-dimensional Lorentzian Lie groups}\label{sec23}

Homogeneous Lorentzian 3-spaces $(N,g)$ where classified by G. Calvaruso  in \cite{C1}. Unless they are symmetric, they are Lie groups equipped with left-invariant Lorentzian metrics.

\subsection{Three-dimensional Lorentzian Lie groups}\label{subsec23.1}

G. Calvaruso  in \cite{C1} classified 3-dimensional simply connected, complete homogeneous Lorentzian manifold as the following theorem.

\begin{Theorem}\label{classif} Let $(N,g)$ be a 3-dimensional connected, simply connected,
complete homogeneous Lorentzian manifold. If $(N,g)$ is not
symmetric, then $N=G$ is a 3-dimensional Lie group and $g$ is
left-invariant. Moreover, there exists a pseudo-orthonormal frame
field $\{ e_1,e_2,e_3\}$, with $e_3$ time-like, such that the Lie
algebra of $G$ is one of the following seven types.
\begin{itemize}

\item[\rm {(1)}] {Type $\mathfrak{g} _1$}:
\begin{equation}\label{g1}
\left[e_1,e_2 \right]=\alpha e_1-\beta e_3, \;\;
\left[ e_1,e_3\right]=-\alpha e_1-\beta e_2, \;\; \left[e_2,e_3\right]=\beta e_1 +\alpha e_2 +\alpha  e_3, \;\; \alpha \neq 0. \end{equation}

In this case, $G= O(1,2)$ or $G=SL(2,\mathbb R)$ if $\beta \neq 0$, while
$G=E(1,1)$ if $\beta=0$.

\item[\rm {(2)}]  {Type $\mathfrak{g} _2$}:
\begin{equation}\label{g2}
\left[e_1,e_2 \right]=\gamma e_2-\beta e_3, \;\;\left[ e_1,e_3\right]=-\beta e_2+\gamma e_3, \;\;  \left[e_2,e_3\right]=\alpha e_1, \;\; \gamma\neq 0 .
\end{equation}
In this case, $G= O(1,2)$ or $G=SL(2,\mathbb R)$ if $\alpha \neq 0$, while
$G=E(1,1)$ if $\alpha=0$.

\item[\rm {(3)}]  {Type $\mathfrak{g} _3$}:
\begin{equation}\label{g3}
\left[e_1,e_2 \right]=-\gamma e_3,\;\; \left[ e_1,e_3\right]=-\beta e_2,\;\; \left[e_2,e_3\right]=\alpha e_1 . 
\end{equation}%
The following Table I lists all the Lie groups $G$ which admit a Lie
algebra $\mathfrak{g}_3$, taking into account the different
possibilities for $\alpha$, $\beta$ and $\gamma$:
\begin{center}
\begin{tabular}{|c|c|c|c|}
\hline $G$ &  $\alpha$ & $\beta$ & $\gamma$ \\
\hline  $O(1,2)$ or $SL(2,\mathbb R)$ & $+$ & $+$ & $+$\\
\hline  $O(1,2)$ or $SL(2,\mathbb R)$ & $+$ & $-$ & $-$\\
\hline  $SO(3)$ or $SU(2)$ & $+$ & $+$ & $-$\\
\hline  $E(2)$ & $+$ & $+$ & $0$\\
\hline  $E(2)$ & $+$ & $0$ & $-$\\
\hline  $E(1,1)$ & $+$ & $-$ & $0$\\
\hline  $E(1,1)$ & $+$ & $0$ & $+$\\
\hline  $H_3$ & $+$ & $0$ & $0$\\
\hline  $H_3$ & $0$ & $0$ & $-$\\
\hline  $\mathbb R \oplus \mathbb R \oplus \mathbb R$ & $0$ & $0$ & $0$\\
\hline
\end{tabular} \nopagebreak \\ \nopagebreak Table I $\vphantom{\displaystyle\frac{a}{2}}$
\end{center}

\item[\rm {(4)}]  Type $\mathfrak g_4$:
\begin{equation}\label{g4}
\left[e_1,e_2 \right]=- e_2 + (2 \eta - \beta) e_3, \;\; \left[ e_1,e_3\right]=-\beta e_2 + e_3,  \;\; \left[e_2,e_3\right]=\alpha e_1, \;\; \eta = \pm 1.  
\end{equation}
The following Table II describes all Lie groups $G$ admitting a Lie
algebra $\mathfrak g_4$:
\begin{gather*}
\begin{array}{c}
\begin{tabular}{|c|c|c|}
\hline $G$ &  $\alpha$ & $\beta$  \\
\hline  $O(1,2)$ or $SL(2,\mathbb R)$ & $\neq 0$ & $\neq \eta$ \\
\hline  $E(1,1)$ & $0$ & $\neq \eta$ \\
\hline  $E(1,1)$ & $<0$ & $\eta$ \\
\hline  $E(2)$ & $>0$ & $\eta$ \\
\hline  $H_3$ & $0$ & $\eta$ \\
\hline
\end{tabular}
\end{array}
 \\
 \text{Table II} \vphantom{\displaystyle\frac{a}{2}}
\end{gather*}
\item[\rm {(5)}]  Type $\mathfrak g_5$:
\begin{equation}\label{g5} \left[e_1,e_2 \right]=0, \;\; \left[ e_1,e_3\right]=\alpha e_1+\beta e_2, \;\; \left[e_2,e_3\right]=\gamma e_1 +\delta e_2, \;\; \alpha +\delta \neq 0, \, \alpha \gamma +\beta \delta =0. 
\end{equation}

\item[\rm {(6)}]  Type $\mathfrak g_6$:
\begin{equation}\label{g6}
\left[e_1,e_2 \right]=\alpha e_2 +\beta e_3, \;\;
\left[ e_1,e_3\right]=\gamma e_2+\delta e_3,  \ \left[e_2,e_3\right]= 0, \;\;  \alpha +\delta
\neq 0, \, \alpha \gamma -\beta \delta =0. \end{equation}

\item[\rm {(7)}]  Type $\mathfrak g_7$:
\begin{equation}\label{g7}
\left[e_1,e_2 \right]=-\alpha e_1-\beta e_2 -\beta e_3, \;\;
\left[ e_1,e_3\right]=\alpha e_1+\beta e_2 +\beta e_3, \;\; \left[e_2,e_3\right]=\gamma e_1 +\delta e_2 +\delta  e_3 ,\end{equation}
with $ \alpha +\delta \neq 0, \, \alpha \gamma =0.$ 
\end{itemize}

Lie algebras of types $\mathfrak g_1$, $\mathfrak g_2$, $\mathfrak
g_3$ and $\mathfrak g_4$ correspond to unimodular groups, whereas
Lie algebras of types $\mathfrak g_5$, $\mathfrak g_6$ and
$\mathfrak g_7$ correspond to non-unimodular groups.
\end{Theorem}

 G. Calvaruso  determined in \cite{C2} those 3-dimensional
Lorentzian Lie groups $(G,g)$ which have constant sectional curvature and which are symmetric.

By a {\it 3-dimensional Lorentzian Lie group $G_i$} we mean a connected, simply connected 3-dimensional Lie group $G$, equipped with a left-invariant Lorentzian metric $g$ and having Lie algebra $\mathfrak g_i$.

\subsection{Classification of parallel surfaces in three-dimensional Lorentzian Lie groups}\label{subsec23.2}

Let $(N,g)$ be a 3-dimensional homogeneous Lorentzian manifold and $M$ a surface in $N$.  We denote by $\xi$ a fixed normal vector field on the surface, with $\<\xi,\xi\>=\varepsilon$. Here, either $\varepsilon=-1$ or $\varepsilon=1$, according to the surface being
either Riemannian or Lorentzian, respectively. We  call $\xi$ an {\it $\varepsilon$-unit normal vector field}.

Parallel surfaces in 3-dimensional Lorentzian Lie groups were classified by G. Calvaruso and J. Van der Veken in \cite{CV10.2}. More precisely, under the notations of Theorem \ref{classif}, they proved the following.

\begin{Theorem}\label{T56} Let $M$ be a parallel surface in a 3-dimensional Lorentzian Lie group $G_1$. Then, $\beta=0$, $\xi=e_1+be_2+be_3$ and the vector fields $E_1=(be_1-e_2)/\sqrt{1+b^2}$ and $E_2=(be_1+b^2e_2+(1+b^2)e_3)/\sqrt{1+b^2}$ form a pseudo-orthonormal basis for the tangent plane at every point. Moreover, the function $b$ satisfies $E_1(b)=E_2(b)$
and
$$E_1\left(\frac{E_1b}{\sqrt{1+b^2}}-\frac{2b}{1+b^2}\alpha\right)
+2\left(\frac{E_1b}{\sqrt{1+b^2}}-\frac{2b}{1+b^2}\alpha\right)\left(\frac{b}{\sqrt{1+b^2}}E_1b-\frac{\alpha}{\sqrt{1+b^2}}\right)=0.$$
The surface is flat and parallel. Moreover, it is totally geodesic
in the case that
$E_1b=E_2b={2 b\alpha}/{\sqrt{1+b^2}}.$
\end{Theorem}

\begin{Theorem}\label{T57} Let $M$ be a parallel surface in a three-dimensional Lorentzian Lie group $G_2$.  Then, one of the following statements holds.
\begin{itemize}
\item[{\rm (a)}] $M$ is an integral surface of the distribution spanned by $\{e_2,e_3\}$. This case only occurs if $\alpha=0$ and $M$ is parallel, flat and minimal, but not totally geodesic. 

\item [{\rm (b)}] $M$ is an integral surface of the distribution spanned by $\{e_1, ce_2+be_3\}$, where $b$ and $c$ are real constants satisfying $b^2-c^2=\varepsilon=\pm 1, \, bc=-{\varepsilon\beta}/{(2\gamma)}.$ This case only occurs if $\alpha=2\beta$ and $M$ is totally geodesic.
\end{itemize}
\end{Theorem}

\begin{Theorem}\label{T58}
Let $M$ be a parallel surface in a non-symmetric three-dimensional Lorentzian Lie group $G_3$. Then, one of the following statements holds. 
\begin{itemize}
\item[{\rm (a)}] $M$ is an integral surface of the distribution
spanned by $\{e_2, e_3\}$. This case only occurs if $\gamma=0$ and
$M$ is flat and minimal, but not totally geodesic. 

\item[{\rm (b)}] $M$ is an integral surface of the distribution spanned
by $\{e_2,e_3\}$. This case only occurs if $\alpha=0$ and  $M$ is
flat and minimal, but not totally geodesic. 

\item [{\rm (c)}] $M$ is an integral surface of the distribution
spanned by $\{e_1, e_3\}$. This case only occurs if $\beta=0$ and
$M$ is flat and minimal, but not totally geodesic. 

\item [{\rm (d)}] $M$ is an integral surface of the distribution
spanned by $\{E_1=e_1,E_2=ce_2+be_3\}$, where $b$ and $c$ are
functions on  $M$ satisfying $b^2-c^2=\varepsilon$ and
$E_1b=\beta c, \, E_1c=\beta b, \, E_2b=k_1\varepsilon c, \,E_2c=k_1\varepsilon b,$
for some real constant $k_1$. This case only occurs if
$\beta=\gamma$ and $M$ is flat. 

\item[{\rm (e)}] $M$ is an integral surface of the distribution spanned by $\{ce_2+be_3,e_1\}$. Here, $b$ and $c$ are real constants satisfying
$b^2={\gamma\varepsilon}/{(\gamma-\beta)}, \, c^2={\beta\varepsilon}/{(\gamma-\beta)}.$
This case only occurs if $\alpha=\beta+\gamma$ and $\beta\neq\gamma$
and $M$ is totally geodesic. 

\item [{\rm (f)}] $M$ is an integral surface of the distribution
spanned by $\{E_1=ce_1+ae_3,E_2=e_2\}$, where $a$ and $c$ are
functions on the surface satisfying $a^2-c^2=\varepsilon$ and
$E_1a=k_2\varepsilon c, \, E_1c=k_2\varepsilon a, \, E_2a=-\alpha c, \, E_2c=-\alpha a,$
for some real constant $k_2$. This case only occurs if
$\alpha=\gamma$ and $M$ is flat. 

\item[{\rm (g)}] $M$ is an integral surface of the distribution spanned by $\{ce_1+ae_3,e_2\}$. Here, $a$ and $c$ are real constants satisfying
$a^2=-{\gamma\varepsilon}/{(\alpha-\gamma)}, \, c^2=-{\alpha\varepsilon}/{(\alpha-\gamma)}.$
This case only occurs if $\beta=\alpha+\gamma$ and
$\alpha\neq\gamma$ and $M$ is totally geodesic. \vskip.05in
\item[{\rm (h)}] $M$ is an integral surface of the distribution spanned
by $\{E_1=be_1-ae_2,E_2=e_3\}$, where $a$ and $b$ are functions
satisfying $a^2+b^2=1$ and
$$E_1a=\frac{k_3b}{b^2-a^2}, \;\; E_1b=-\frac{k_3a}{b^2-a^2}, \;\; E_2a=\frac{b\alpha}{b^2-a^2}, \;\;E_2b=-\frac{a\alpha}{b^2-a^2},$$
for some real constant $k_3$. This case only occurs if
$\alpha=\beta$ and $M$ is flat. 

\item[{\rm (i)}] $M$ is an integral surface of the distribution spanned by $\{be_1-ae_2,e_3\}$, where $a$ and $b$ are constants satisfying $a^2=-{\beta}/{(\alpha-\beta)}, \, b^2={\alpha}/{(\alpha-\beta)}.$ This case only occurs if $\gamma=\alpha+\beta$ and $\alpha\neq\beta$
and $M$ is totally geodesic.
\end{itemize}
\end{Theorem}

\begin{Theorem}\label{T59} Let $M$ be a parallel surface in a non-symmetric three-dimensional Lorentzian Lie group $G_4$. Then one of the following statements holds. 
\begin{itemize}
\item[(a)] $M$ is an integral surface of the distribution spanned by $\{e_2,e_3\}$. This case only occurs if $\alpha=0$. $M$ is parallel, flat and minimal, but not totally geodesic. \vskip.05in
\item[(b)] $M$ is an integral surface of the distribution spanned by $\{e_1,ce_2+be_3\}$, where $b$ and $c$ are constants satisfying $b^2-c^2=\varepsilon$ and $\beta b^2+2bc+(\beta-2\eta)c^2=0$. $M$ is totally geodesic and has constant Gaussian curvature
$G=-\varepsilon(\beta-\eta)$. 
\end{itemize}
\end{Theorem}

\begin{Theorem}\label{T60}
Let $M$ be a parallel surface in a non-symmetric three-dimensional
Lorentzian Lie group $G_5$. Then $M$ is one of the surfaces listed
below. 
\begin{itemize}
\item[{\rm (a)}] $M$ is an integral surface of the distribution spanned
by $e_1$ and $e_2$. $M$ is flat but not totally geodesic.

\item[{\rm (b)}] $M$ is an integral surface of the distribution spanned
by $e_2$ and $e_3$. This case only occurs if either $\beta=\gamma=0$
or $\gamma=\delta=0$. In the first case, $M$ is totally geodesic and
has constant Gaussian curvature $K=-\delta^2\leq 0$. In the second
case, $M$ is flat and minimal, but not necessarily totally geodesic.

\item[{\rm (c)}] $M$ is an integral surface of the distribution spanned by $e_1$ and $e_3$. This case only occurs if either $\alpha=\beta=0$ or $\beta=\gamma=0$. In the first case, $M$ is flat and minimal, but not necessarily totally geodesic. In the second case, $M$ is totally
geodesic and has constant Gaussian curvature $K=\alpha^2\geq 0$.

\item[{\rm (d)}] $M$ is an integral surface of the distribution spanned
by $\{E_1=e_1, E_2=ce_2+be_3\}$, where $b$ and $c$ are functions
satisfying $b^2-c^2=\varepsilon$ and
$E_1b=E_1c=0, \, E_2b=c(k_1-c\delta), \, E_2c=b(k_1-c\delta),$
for some real constant $k_1$. This case only occurs if
$\alpha=\beta=0$ and $M$ is flat. 

\item[{\rm (e)}] $M$ is an integral surface of the distribution spanned
by $\{E_1=ce_1+ae_3,E_2=e_2\}$, where $a$ and $c$ are functions satisfying $a^2-c^2=\varepsilon$ and
$E_1a=-\varepsilon c(a^2c\alpha-k_2), \, E_1c=-\varepsilon a(a^2c\alpha-k_2), \,E_2a=E_2c=0,$
for some real constant $k_2$. This case only occurs if $\gamma=\delta=0$ and $M$ is flat.
\end{itemize}
\end{Theorem}

\begin{Theorem}\label{61}
Let $M$ be a parallel surface in a three-dimensional Lorentzian Lie
group $G_6$. Then, one of the following statements holds.\begin{itemize}
\item[{\rm (a)}] $M$ is an integral surface of the distribution spanned
by $e_1$ and $e_2$. This case only occurs if either $\alpha=\beta=0$
or $\beta=\gamma=0$. In the first case, $M$ is parallel, flat and
minimal, but not necessarily totally geodesic. In the second case,
$M$ is totally geodesic. 
\item[{\rm (b)}] $M$ is an integral surface of the distribution spanned
by $e_2$ and $e_3$. $M$ is parallel and flat, but not necessarily
totally geodesic. 
\item[{\rm (c)}] $M$ is an integral surface of the distribution spanned
by $e_1$ and $e_3$. This case only occurs if either $\beta=\gamma=0$
or $\gamma=\delta=0$. In the first case, $M$ is totally geodesic. In
the second case, $M$ is parallel, flat and minimal, but not
necessarily totally geodesic. 
\item[{\rm (d)}] $M$ is an integral surface of the distribution spanned
by $\{E_1=ce_1+ae_3,E_2=e_2\}$, where $a$ and $c$ are functions
satisfying $a^2-c^2=\varepsilon$ and
$E_1a=c(k_1-\delta a),  \, E_1c=a(k_1-\delta a), \, E_2a=E_2c=0,$
for some real constant $k_1$. This case only occurs if
$\alpha=\beta=0$ and $M$ is parallel and flat. 

\item[{\rm (e)}] $M$ is an integral surface of the distribution spanned
by $\{E_1=be_1-ae_2,E_2=e_3\}$, where $a$ and $b$ are functions
satisfying $a^2+b^2=1$ and
$E_1a=b(k_2+\alpha a), \, E_1b=-a(k_2+\alpha b), \, E_2a=E_2c=0, $
for some real constant $k_2$. This case only occurs if
$\gamma=\delta=0$ and $M$ is parallel and flat.
\end{itemize}
\end{Theorem}

\begin{Theorem}\label{T62} Let $M$ be a parallel surface in a non-symmetric three-dimensional
Lorentzian Lie group $G_7$. Then, $M$ is one of surfaces listed below. 
\begin{itemize}
\item[{\rm (a)}] $M$ is an integral surface of the distribution spanned by $\{e_2,e_3\}$. This case only occurs if either $\beta=\gamma=0$ or $\gamma=\delta=0$. In the first case, $M$ is totally geodesic. In the second case, $M$ is parallel and flat, but not necessarily totally geodesic. 

\item[{\rm (b)}] $M$ is an integral surface of the distribution spanned by $\{E_1=e_1,E_2=ce_2+be_3\}$, where $b$ and $c$ are functions satisfying $b^2-c^2=\varepsilon$ and $E_1b=E_1c=0, \, E_2b=c((b-c)\delta-k_1), \, E_2c=b((b-c)\delta-k_1)$
for some real constant $k_1$. This case only occurs if $\alpha=\beta=0$. $M$ is flat, but not necessarily totally geodesic.

\item[{\rm (c)}] $M$ is an integral surface of the distribution spanned by $E_1=(be_1-e_2)/\sqrt{1+b^2}$ and $E_2=(be_1+b^2e_2+(1+b^2)e_3)/\sqrt{1+b^2}$, where $b$ is a
function satisfying $E_1(b)=E_2(b)$ and
\begin{align*}
E_1\left(\frac{E_1b}{\sqrt{1+b^2}}+\frac{b(\alpha-\delta)}{1+b^2}\right)
+2\left(\frac{E_1b}{\sqrt{1+b^2}}+\frac{b(\alpha-\delta)}{1+b^2}\right)
\left(\frac{bE_1b}{\sqrt{1+b^2}}-\frac{\delta}{\sqrt{1+b^2}}\right)=0.
\end{align*}
The surface is flat and parallel. Moreover, it is totally geodesic in the special case that
$E_1b=E_2b={b(\delta-\alpha)}/{\sqrt{1+b^2}}.$
\end{itemize}
\end{Theorem}

\section{Parallel  surfaces in reducible three-spaces}\label{sec24}

\subsection{Classification of parallel surfaces in reducible three-spaces}\label{subsec24.1}

Parallel submanifolds of the a Robertson-Walker space-time $I\times_f  R^n(c)$ have been treated in \S20. In \cite{CV13.2}, G. Calvaruso and J. Van der Veken studied  parallel surfaces in 3-dimensional reducible spaces ${\mathbb M}^2\times {\mathbb E^1}$. More precisely, they proved the following results.

\begin{Theorem}\label{T63} Let $M$ be a parallel surface in a reducible 3-dimensional Riemannian
manifold ${\mathbb M}^2\times {\mathbb E}^1$. Then one of the following three cases holds:

\begin{itemize}
\item[{\rm (1)}] $M$ is isometric to an open portion of a surface of type ${\mathbb M}^2\times \{t_0\}$ for some $t_0\in {\bf R}$;

\item[{\rm (2)}]  $M$ is isometric to an open portion of a surface of type $\gamma\times {\mathbb E}^1$, where $\gamma$ is a curve of constant geodesic curvature in $M$;

\item[{\rm (3)}]  ${\mathbb M}^2\times {\mathbb E}^1$ is flat, and $M$ is isometric to an open portion of a standard sphere $S^2\subset \mathbb E^3$.
\end{itemize}\end{Theorem}

The following is a consequence of Theorem \ref{T63}.

\begin{Corollary}\label{C8} The pair $(S^2,{\mathbb E}^3)$ is the only proper parallel surface in a reducible Riemannian 3-space.
\end{Corollary}

For parallel surfaces in a reducible 3-dimensional Lorentzian manifold, G. Calvaruso and J. Van der Veken obtained the following.

\begin{Theorem}\label{T64}  Let $M$ be a  parallel surface in a reducible 3-dimensional Lorentzian manifold ${\mathbb M}_1^2\times {\mathbb E^1}$ (respectively ${\mathbb M}^2\times {\mathbb E}_1^1$). Then one of the
following holds.
\begin{itemize}
\item[{\rm (1)}] $M$ is isometric to an open portion of a surface of type ${\mathbb M}_1^2\times \{t_0\}$ (respectively ${\mathbb M}^2\times \{t_0\}$) for some real number $t_0$.

\item[{\rm (2)}] $M$ is isometric to an open portion of a surface of type $\gamma\times {\mathbb E^1}$ (respectively $\gamma\times {\mathbb E}^1_1$) where $\gamma$ is a non-degenerate curve of constant geodesic curvature in ${\mathbb M}_1^2$ (respectively ${\mathbb M}^2$).
\item[{\rm (3)}]  The ambient space is  flat, and $M$ is isometric to an open portion of
one of the following surfaces:
{\rm (a)} a hyperbolic plane $H^2$;
{\rm (b)} an indefinite sphere $S^2_1$;
{\rm (c)} the null scroll $N^2_1$.
\end{itemize}\end{Theorem}

As a consequence of Theorem \ref{T64}, G. Calvaruso and J. Van der Veken obtained the following.

\begin{Corollary}\label{C9} The pairs $(H^2,{\mathbb E}^3_1 )$, $(S^2_1,{\mathbb E}^3_1)$ and $(N^2_1,{\mathbb E}^3_1)$ are the only proper parallel surfaces in a reducible Lorentzian 3-space.
\end{Corollary}

\subsection{Parallel surfaces in Walker three-manifolds}\label{subsec24.2}

A particularly interesting class of pseudo-Riemannian manifolds are ones which admit a
parallel null vector field. The study of such metrics in the 3-dimensional Lorentzian setting was initiated by M. Chaichi, E. Garc\'{i}a-R\'{i}o and M. E. V\'{a}zquez-Abal in \cite{ChGRVA}.    W. Batat and S. J. Hall named such manifolds as {\it Walker manifolds}  in \cite{BH17}.

Complete classification of  parallel surfaces of an arbitrary reducible 3-manifold, both in Riemannian and Lorentzian was obtained by G. Calvaruso and J. Van der Veken in \cite{CV13.2}.  It turns out that the Euclidean space $\mathbb E^3$ and the Minkowski space 
$\mathbb E^3_1$ are the only cases admitting parallel surfaces which
are non-trivial, in the sense that they do not reflect the reducibility of the space itself. 
Since the reducibility of a pseudo-Riemannian manifold corresponds to the existence of a parallel non-null vector field, it is natural to study parallel surfaces in a Lorentzian 3-manifold which admits a parallel null vector field, i.e., in a Walker 3-manifold.
 G. Calvaruso and J. Van der Veken provided in \cite{CV10.1} a complete classification of parallel surfaces in Walker 3-manifolds.

In \cite{BH17}, W. Batat and S. J. Hall proved that totally umbilical nondegenerate surfaces
in a Walker 3-manifold with metric $g= \epsilon dx^2+f(x,y)dy^2+2dtdy$ with $\epsilon=\pm 1$ and satisfying $f_{xx}\ne 0$ are either one of a totally geodesic family described by G. Calvaruso and J. Van der Veken in \cite{CV10.1} or the ambient manifold must be locally conformally flat (here the surface can also be totally geodesic).

\section{Bianchi-Cartan-Vranceanu spaces}\label{sec25}

\subsection{Basics on Bianchi-Cartan-Vranceanu spaces}\label{subsec25.1}

The simply-connected homogeneous 3-manifolds are classified according to the dimension of their isometry group which is equal to 3, 4 or 6. If it is  6, one obtains the real space forms. 
The Bianchi-Cartan-Vranceanu spaces are homogeneous Riemannian 3-manifolds with isometry group of dimension 4 or 6.
Such spaces, denoted by $\widetilde{\mathcal M}^3(\lambda,\mu)$, are given by a two-parameter family of Riemannian 3-manifolds $({\mathcal M},g_{\lambda,\mu})$ where the underlying 3-manifolds $\widetilde{\mathcal M}^3$ are ${\mathbb R}^3$ if $\mu\geq 0$; and
$$\widetilde{\mathcal M}^3=\left\{(x,y,z)\in {\mathbb R}^3: x^2+y^2<-\frac{1}{\mu}\right\}\;\;\; {\rm if}\;\; \mu<0.$$ 
The metrics $\tilde g_{\lambda,\mu}$ on $\widetilde{\mathcal M}^3$ are given by
\begin{align}&\label{32}  g_{\lambda,\mu}
=\frac{{\rm d}x^2+{\rm d}y^2}{\{1+\mu (x^2+y^2)\}^2} 
+\( {\rm d}z+\frac{\lambda (y{\rm d}x-x{\rm d}y)}{2\{ 1+\mu(x^2+y^2)\}}\)^2. \end{align}
The 2-parameter family $\tilde g_{\lambda,\mu}$ is called the {\it Bianchi-Cartan-Vranceanu metrics}. The metrics above are defined over the whole 3-space $\mathbb R^3$ for $\mu>0$ and over the region $x^2+y^2<-1/\mu$ for $\mu<0$. 

Consider the following Riemannian surface with constant Gaussian curvature $4\mu$:
$$\widetilde{\mathcal M}^2(\mu)=\left(\left\{(x,y)\in\mathbb R^2 : 1+\mu (x^2+y^2)>0 \right\}\ ,\ \frac{dx^2+dy^2}{(1+\mu (x^2+y^2))^2}\right).$$
Then the mapping
$$\pi:\widetilde{\mathcal M}^3({\lambda,\mu})\to \widetilde{\mathcal M}^2(\mu):(x,y,z)\mapsto (x,y)$$
is a Riemannian submersion, referred to as the \emph{Hopf-fibration}. For $\mu=4\lambda^2\neq 0$, this mapping coincides with the ``classical'' Hopf-fibration
$\pi:S^3\left(\mu\right)\to S^2(4\mu)$. 

In the following, by a \emph{Hopf-cylinder} we mean the inverse image of a curve in $\widetilde{\mathcal M}^2(\mu)$ under $\pi$. By a \emph{leaf} of the
Hopf-fibration, we mean a surface which is everywhere orthogonal to the fibres. 

The family of Bianchi-Cartan-Vranceanu spaces $\widetilde{\mathcal M}^3(\lambda,\mu)$ includes six of the eight Thurston's  3-dimensional geometries except $Sol_3$ and the hyperbolic space $H^3$.
The family of the Riemannian metrics given by \e{32} includes all 3-dimensional homogeneous metrics whose group of isometries has dimension 4 or 6, except for those with negative constant curvature.

For two given real numbers $\lambda,\mu$, the Bianchi-Cartan-Vranceanu space $\widetilde{\mathcal M}^3(\lambda,\mu)$ is the following 3-spaces (cf. e,g., \cite{Ve07,Ve08,OC19}).

\begin{itemize}

  \item[{\rm (1)}] If $\lambda=\mu=0$, it is the Euclidean 3-space.
  
  \item[{\rm (2)}] If $\lambda=0$, $\mu\ne 0$,  it is the product of real line and a surface of constant curvature $4\lambda$.

 \item[{\rm (3)}] If $\lambda\ne 0$, $\lambda^2=4\mu$,  it is a space of positive constant curvature.

 \item[{\rm (4)}] If $\lambda\ne 0$, $\mu>0$,  it is $SU(2)\setminus \{\infty\}$.

 \item[{\rm (5)}] If $\lambda\ne 0$, $\mu<0$,  it is $\widetilde{SL_2}({\mathbb R})$ with a left-invariant metric.

 \item[{\rm (6)}] If $\lambda\ne 0$, $\mu=0$, it is the Heisenberg group $Nil_3$ with a left-invariant metric.
\end{itemize}

\subsection{$B$-scrolls}\label{subsec25.2}

For every $\gamma$ in  the unit 3-sphere $S^3(1)$, one can define the Frenet frame $\{T,N, B\}$ provided the geodesic curvature $\kappa$ does not vanish.
The {\it $B$-scroll} of a curve $\gamma$ in  the unit 3-sphere $S^3(1)$ is the surface described by moving the geodesic through $\gamma(s)$ in the direction of spherical binormal $B(s)$ along $\gamma$. A curve in  $S^3(1)$ of constant geodesic curvature  and constant torsion $\pm 1$ is called a {\it twisted spherical spiral}. The $B$-scroll of a twisted spherical spiral has parallel second fundamental form (cf. \cite{Di92}), so it is a parallel surface in $S^3(1)$.

If $\gamma$ is a closed curve in $S^2(\frac{1}{2})$, then the Hopf cylinder $\pi^{-1}(\gamma)$ is called a {\it  Hopf torus}.
 A $B$-scroll of a twisted spherical spiral is a Hopf cylinder (torus) over a curve  with constant curvature in $S^2(\frac{1}{2})$ (cf. \cite{BDI02}).

\subsection{Parallel surfaces in Bianchi-Cartan-Vranceanu spaces}\label{subsec25.3}

If $4\mu=\lambda^2$, then $\widetilde{\mathcal M}^3(\lambda,\mu)$ is a real space form whose parallel surfaces are already known. In the next theorem,  M. Belkhelfa, F. Dillen and J.  Inoguchi \cite{BDI02} classified  parallel surfaces in Bianchi-Cartan-Vranceanu spaces $\widetilde{\mathcal M}^3(\lambda,\mu)$ with  $4\mu\ne \lambda^2$.

\begin{Theorem}\label{T65} Let $\widetilde{\mathcal M}^3(\lambda,\mu)$ be a Bianchi-Cartan-Vranceanu space with $4\mu\ne \lambda^2$. 
\begin{itemize}
 \item[{\rm (1)}] If $\lambda\ne 0$, then the only parallel surfaces in $\widetilde{\mathcal M}^3(\lambda,\mu)$ are Hopf cylinders over curves with constant curvature in $\widetilde{\mathcal M}^2(\mu)$.   
 \item[{\rm (2)}]  If $\lambda=0$, then the only parallel surfaces in $\widetilde{\mathcal M}^3(\lambda,\mu)$ with $\mu\ne 0$ are totally geodesic leaves and Hopf cylinders over circles with constant geodesic curvature in $\widetilde{\mathcal M}^2(\mu)$.
\end{itemize}
\end{Theorem}

\section{Parallel surfaces in homogeneous three-spaces}\label{sec26}

\subsection{Homogeneous three-spaces}\label{subsec2.1}

A Riemannian manifold $M$ is said to be {\it homogeneous} if for any two points $p$ and $q$ of $M$ there exists an isometry of $M$ which carries $p$ into $q$. It is clear that these spaces are a natural generalization of real space forms. A parallel submanifold is called {\it proper parallel} it is non-totally geodesic. In dimension 3, the classification of
these spaces is well known as follows.

\begin{Theorem} \label{T66} Let $M^3$ be a simply connected homogeneous Riemannian manifold with isometry group $I(M^3)$, i.e., $I(M^3)$ acts transitively on $M^3$.
Then $\dim I(M^3)\in\{3,4,6\}$ and moreover:
\begin{itemize}

\item[{\rm (i)}] if $\dim I(M^3)=6$, then $M^3$ is a real space form of constant sectional curvature $c$, i.e., Euclidean space $\mathbb E^3$, hyperbolic space $H^3(c)$ or a
three-sphere $S^3(c)$,

\item[{\rm (ii)}] if $\dim I(M^3)=4$, then $M^3$ is a Bianchi-Cartan-Vranceanu space (different from $\mathbb E^3$ and $S^3(c)$), i.e. a Riemannian product $H^2(c)\times\mathbb R$
or $S^2(c)\times\mathbb R$, or one of following Lie groups, equipped with a left-invariant metric yielding a four-dimensional isometry group: the special unitary group ${SU}(2)$, the universal
covering of the special linear group $\widetilde{{SL}}(2,\mathbb R)$ or the Heisenberg group
$\mathrm{Nil}_3$,

\item[{\rm (iii)}] if $\dim I(M^3)=3$, then $M^3$ is a general
three-dimensional Lie group with left-invariant metric.
\end{itemize}
\end{Theorem}

\subsection{Classification of parallel surfaces in homogeneous three-spaces}\label{subsec26.2}

In \cite{IVan08,V07}, J. Inoguchi and J. Van der Veken classified parallel surfaces in homogeneous 3-spaces in the next two theorems.

For totally geodesic surfaces in a  3-dimensional homogeneous Riemannian manifold, we have:

\begin{Theorem}\label{T67}  Let $(M^3, g)$ be a  3-dimensional homogeneous Riemannian manifold. Then
$M^3$ admits totally geodesic surfaces if and only if $M^3$ is locally isometric to one of the
following spaces:
\begin{itemize}
\item[{\rm (1)}]  a real space form $S^3, \mathbb E^3$ or $H^3$,
\item[{\rm (2)}]  a Riemannian product space $S^2\times \mathbb E^1$  or $H^2\times \mathbb E^1$,
\item[{\rm (3)}] $SL(2,{\mathbb R})$ with a left-invariant metric determined by the condition $c_2 = c_1 + c_3$, or equivalently
$\mu_2=0$,
\item[{\rm (4)}]  the Minkowski motion group $E(1, 1)$ with Riemannian 4-symmetric metric,
including the model space $Sol_3$,
\item[{\rm (5)}]  a non-unimodular Lie group with structure constants $(\xi,\eta)$ satisfying 
$\xi\notin \{0,1\}$ and $\eta=0$.
\end{itemize}\end{Theorem}

For proper parallel surfaces in a $3$-dimensional homogeneous Riemannian manifold, J. Inoguchi and J. Van der Veken \cite{IVan08} proved the following.

\begin{Theorem}\label{T68} Let $(M^3, g)$ be a $3$-dimensional homogeneous Riemannian manifold. Then $M^3$ admits proper parallel surfaces if and only if $M^3$ is locally isometric to one of the
following spaces:
\begin{itemize}
\item[{\rm (1)}]  a real space form $S^3, \mathbb E^3$ or $H^3$,
\item[{\rm (2)}]  a Bianchi-Cartan-Vranceanu space,
\item[{\rm (3)}]  the Minkowski motion group E(1, 1) with any left-invariant metric, including the model space $Sol_3$,
\item[{\rm (4)}]  the Euclidean motion group $E(2)$ with any left-invariant metric,
\item[{\rm (5)}]  a non-unimodular Lie group with structure constants $(\xi,\eta)$  satisfying $\xi\notin \{0,1\}$.
\end{itemize}\end{Theorem}

\section{Parallel surfaces in symmetric Lorentzian three-spaces}\label{sec27}

Symmetric spaces are one of the most important topics in Riemannian geometry. In the Lorentzian setting, their study goes back to the work of M. Cahen and N. Wallach \cite{CW70} in the 1970s. 

\subsection{Symmetric Lorentzian  three-spaces}\label{subsec27.1}

It is well known that the curvature of a 3-dimensional pseudo-Riemannian manifold $(N,g)$ is completely determined by the Ricci tensor, denoted by $Ric$, defined for any point $p \in N$ and any $X,Y \in T_p N$ by
 \begin{equation}\label{Ric} Ric(X,Y) _p = \sum _{i=1} ^3 \varepsilon _i g(R(X,e_i)Y,e_i),
\end{equation}
where $R$ is the Riemann curvature tensor, $\{e_1,e_2,e_3\}$ is a pseudo-orthonormal basis of $T_p N$ and $ \varepsilon _i = g_p (e_i,e_i)= \pm 1$ for all $i$. Throughout this section, if not stated otherwise, we shall assume that $e_3$ is {\em time-like}, i.e., $\varepsilon _1 = \varepsilon _2 = - \varepsilon _3 =1$.

Due to the symmetries of the curvature tensor, the Ricci tensor $Ric$ is symmetric \cite{O}. Thus, the {\em Ricci operator} $Q$, defined by $g(QX,Y)=Ric(X,Y)$, is self-adjoint. In the Riemannian case, there always exists an orthonormal basis diagonalizing $Q$, but in the Lorentzian case four different cases can occur \cite{O}, and  there exists a pseudo-orthonormal basis $\{e_1,e_2,e_3\}$, with $e_3$ time-like, such that $Q$ takes one of the following canonical forms, called {\em Segre types}:
 \begin{equation}\begin{aligned}\label{segre}&
 {\rm Segre \, type} \, \{11,1 \}: \left(
\begin{array}{ccc}
a & 0 & 0 \\ 0 & b & 0 \\ 0 & 0 & c
\end{array}   \right), \;\;
  {\rm Segre \, type} \, \{1 z \bar{z} \}: \left(
\begin{array}{ccc}
a & 0 & 0 \\ 0 & b & c \\ 0 & -c & b
\end{array}   \right), \; c \neq 0,
 \\&   {\rm Segre \, type} \, \{21\}: \left(
\begin{array}{ccc}
a & 0 & 0 \\ 0 & b & \eta \\ 0 & -\eta & b - 2\eta
\end{array}   \right), \, \eta = \pm 1, \;\;
{\rm Segre \, type} \, \{3 \}: \left(
\begin{array}{ccc}
b & a & -a \\ a & b & 0 \\ a & 0 & b
\end{array}   \right), \, a \neq 0.
\end{aligned}\end{equation}
When $(N,g)$ is homogeneous, the Ricci operator $Q$ has the same Segre type at any point $p \in N$ and has constant eigenvalues.

G. Calvaruso studied homogeneous Lorentzian 3-manifolds $(N^3,g)$  in \cite{C1,C2}. For symmetric ones, he proved that 3-dimensional symmetric spaces can only occur for some Segre types of the Ricci operator $Q$. More precisely, he proved the following:

\medskip\noindent
I) For Segre type $\{11,1\}$, $(N,g)$ is symmetric if and only if
\begin{itemize}

\item[{\rm (i)}] $a=b=c$. Then, $(N,g)$ is an Einstein manifold and hence it has constant sectional curvature. If $N$ is connected and simply connected, then
$(N,g)$ is isometric to one of the Lorentzian space forms: either
$\mathbb S _1 ^3$, $\mathbb R_1 ^3$ or $\mathbb H ^3_1$.

\item[{\rm (ii)}] $a = b \neq c$.
Then, $N$ is reducible as a direct product $M^2 \times \mathbb R _1$, where $M^2$ is a Riemannian surface of constant curvature. If $N$ is connected and simply connected, $(N,g)$ is then isometric to either $\mathbb S ^2 \times \mathbb R_1$ or  $\mathbb H ^2 \times \mathbb R_1$.

\item[{\rm (iii)}] $a \neq b = c$. Then, $N$ is reducible as a direct product $\mathbb R \times M^2 _1$, where $M^2 _1$ is a Lorentzian surface of constant sectional curvature. When $N$ is connected and simply connected, $(N,g)$ is isometric to either $\mathbb R \times \mathbb S ^2 _1$ or  $\mathbb R \times \mathbb H ^2 _1$.
\end{itemize}

\medskip\noindent
II) For Segre type $\{21 \}$, $(N,g)$ is symmetric if and only if
$a -b = \eta$ and, with respect to a suitable pseudo-orthonormal
frame field $\{e_1,e_2,e_3\}$, the Levi Civita connection of $(N,g)$
is completely described by
\begin{equation}\label{nabla21}
\begin{array}{lll}
\nabla _{e_1} e_1 = A e_2 - A e_3 , & \nabla _{e_2} e_1 =  B e_2 - B e_3, & \nabla _{e_3} e_1 = C e_2 -C e_3,  \\
\nabla _{e_1} e_2 = -A e_1, & \nabla _{e_2} e_2 =  -B e_1, & \nabla _{e_3} e_2 = -C e_1,  \\
\nabla _{e_1} e_3 = -A e_1, & \nabla _{e_2} e_3 =  -B e_1, & \nabla
_{e_3} e_3 = -C e_1,
\end{array}
\end{equation}
where $A,B,C$ are smooth functions. Put $u = e_2-e_3$. Then, $\nabla_{e_i} u =0$ for all $i$, that is, $u$ is a {\em parallel null vector field}. Three-dimensional symmetric spaces admitting a
parallel null vector field were described in \cite{ChGRVA} in terms of local coordinates. In fact, a three-dimensional locally symmetric Lorentzian manifold $(N,g)$, having a parallel null vector field, admits local coordinates $(t,x,y)$ such that, with respect to the local frame field $\{ \frac{\partial}{\partial t},\frac{\partial}{\partial x}, \frac{\partial}{\partial y} \}$, the
Lorentzian metric $g$ and the Ricci operator are respectively given by
\begin{equation}\label{Rio} g=\left(\begin{array}{ccc} 0 &  0 & 1   \\ 0 & \varepsilon & 0 \\
1 & 0 & f \end{array} \right),  \quad
 Q=\left(\begin{array}{ccc} 0 &  0 & -\varepsilon \alpha   \\ 0 & 0 & 0 \\ 0 & 0 & 0\end{array}
\right),\end{equation}
where $\varepsilon = \pm 1$, $u = \frac{\partial}{\partial t}$ and
\begin{equation}\label{Rio1} f(x,y)=x^2 \alpha +x \beta (y) + \xi (y),
\end{equation}
for a constant $\alpha \in \mathbb R$ and a functions $\beta,\xi$ (cf. \cite{ChGRVA}). It is easy to build a (local) pseudo-orthonormal frame field from $\{ \frac{\partial}{\partial t},
\frac{\partial}{\partial x}, \frac{\partial}{\partial y} \}$, and to check that, apart from the flat case $\alpha f \neq 0$, the Ricci operator $Q$ described by (\ref{Rio}) is of degenerate Segre type
$\{ 21 \}$, with $\lambda =0$ as the only Ricci eigenvalue, of multiplicity three,  associated to a 2-dimensional eigenspace.

\medskip\noindent
III) For Segre types either $\{1z \bar{z} \}$ or $\{3 \}$, $(N,g)$ is never symmetric.

\medskip\noindent
Therefore, we have the following classification result from \cite{C2}.

\begin{Theorem}\label{T69}
A connected, simply connected three-dimensional symmetric  Lorentzian space $(N,g)$ is either

\medskip
\emph{(i)} a Lorentzian space form $\mathbb S^3 _1$, $\mathbb R^3_1$ or $\mathbb H^3 _1$, or

\medskip
\emph{(ii)} a direct product $\mathbb R \times \mathbb S ^2 _1$, $\mathbb R \times \mathbb H ^2 _1$, $\mathbb S ^2 \times \mathbb R_1$ or $\mathbb H ^2 \times \mathbb R_1$, or

\medskip
\emph{(iii)} a space with a Lorentzian metric $g$ locally described by {\rm (\ref{Rio})-(\ref{Rio1})}.
\end{Theorem}

\subsection{Classification of parallel surfaces in symmetric Lorentzian three-spaces}
\label{subsec27.2}

Three-dimensional Lorentzian manifolds admitting a parallel null vector field were first studied in  \cite{ChGRVA}, in  which the attention was focused on local properties. G.  Calvaruso and J. Van der Veken described in \cite{CVan09} a global model carrying a metric  described by (\ref{Rio})-(\ref{Rio1}) as follows. 

First they showed that the curvature components with respect to the pseudo-orthonormal frame field $\{e_1,e_2,e_3\}$ for which (\ref{nabla21}) holds and then apply (\ref{Ric}) to obtain its Ricci components. Since the Ricci operator must be of degenerate Segre type $\{21\}$ (that is, with $a=b-\eta$), standard calculations lead to the following system of partial differential equations:
\begin{equation}\label{syst21}
\left\{
\begin{array}{l}
\vphantom{\displaystyle\frac{a}{a}}  -e_1 (B) + e_2 (A) + e_1 (C) -e_3 (A)- (B-C)^2 =b -\eta, \\
\vphantom{\displaystyle\frac{a}{a}}  -e_1 (B) + e_2 (A) -A^2-B^2+BC=b, \\
\vphantom{\displaystyle\frac{a}{a}} e_1 (C) - e_3 (A) +A^2 -C^2 +BC =b- 2 \eta,  \\
\vphantom{\displaystyle\frac{a}{a}} e_2 (C) - e_3 (B)+ A(B-C) =0,\\
- e_1 (B)+ e_2 (A) -A^2 -B^2 +BC=\eta .
\end{array}
\right.
\end{equation}
System (\ref{syst21}) implies $a=b-\eta=0$ (which also follows from (\ref{Rio})-(\ref{Rio1})), and the remaining equations reduce to
\begin{equation}\label{syst21bis}
\left\{
\begin{array}{l}
\vphantom{\displaystyle\frac{a}{a}}  e_1 (B) - e_2(A) = -A^2-B^2+BC -\eta, \\
\vphantom{\displaystyle\frac{a}{a}}  e_1 (C) - e_3 (A)= -A^2+C^2-BC -\eta, \\
\vphantom{\displaystyle\frac{a}{a}} e_2 (C) - e_3 (B)=A(C-B).
\end{array}
\right.\end{equation}
Then they proved that,  for any smooth function $\omega$, with respect to the following new pseudo- orthonormal frame field
\begin{equation}\label{newbasis}
e' _1 = e_1 + \omega e_2 -\omega e_3, \; \vphantom{\displaystyle\frac{A}{B}}  e' _2 = -\omega e_1 + (1-\frac{\omega^2}{2}) e_2+ \frac{\omega^2}{2} e_3, \; \vphantom{\displaystyle\frac{A}{B}}  e' _3 =  -\omega e_1 - \frac{\omega^2}{2} e_2+ (1+\frac{\omega^2}{2}) e_3,
\end{equation}
the Ricci operator still keeps the same components than with respect to $\{e_1,e_2,e_3\}$. 
It follows from (\ref{nabla21}) and (\ref{newbasis}) that, with respect to $\{e'_1,e'_2,e'_3\}$, the Levi Civita connection satisfies
\begin{equation}\label{newnabla21}\begin{array}{lll}
\nabla _{e'_1} e'_1 = A' e'_2 - A' e'_3 , & \nabla _{e'_2} e'_1 =  B' e'_2 - B' e'_3, & \nabla _{e'_3} e'_1 = C' e'_2 -C' e'_3,  \\
\nabla _{e'_1} e'_2 = -A' e'_1, & \nabla _{e'_2} e'_2 =  -B' e'_1, & \nabla _{e'_3} e'_2 = -C' e'_1,  
\\ \nabla _{e'_1} e'_3 = -A' e'_1, & \nabla _{e'_2} e'_3 =  -B' e'_1, & \nabla _{e'_3} e'_3 = -C' e'_1,
\end{array}\end{equation}
where
$$\begin{array}{l} A' = A +e_1 \omega, \\ \vphantom{\displaystyle\frac{a}{a}} B'= A
\omega +\omega e_1 \omega-(1-\frac{\omega ^2}{2})B-(1-\frac{\omega
^2}{2})e_2\omega-\frac{\omega ^2}{2}C -\frac{\omega ^2}{2}
e_3\omega, \\ \vphantom{\displaystyle\frac{a}{a}} C'=A \omega
+\omega e_1 \omega +\frac{\omega ^2}{2}B+\frac{\omega
^2}{2}e_2\omega-(1+\frac{\omega ^2}{2})C -(1+\frac{\omega ^2}{2})
e_3\omega.\end{array}$$
Thus, by choosing $\omega$ to be a solution of the system of differential equations
\begin{equation}\label{sigmaomega}
A+e_1\omega=k, \;\; d e_2\omega -e_3\omega=C-B,
\end{equation}
where $k$ is a real constant, we can always specify the pseudo-orthonormal frame field $\{e_1,e_2,e_3\}$ in such a way that  $A=k$ and $B=C$. In this case, system of equations (\ref{syst21bis}) reduces to
\begin{equation}\label{newsyst21} e_1 B = -k^2-\eta,\;\; e_2 B - e_3 B= 0,
\end{equation}
and the Lie brackets $[e_i,e_j]$ are easily determined as follows:
\begin{equation}\label{Lie21}
\begin{array}{l} [ e_1 , e_2 ] = [e_1,e_3 ] =-k e_1 -B(e_2-e_3), \;\; [ e_2 , e_3 ] = 0.\end{array}
\end{equation}

In \cite{CVan09}, G.  Calvaruso and J. Van der Veken proved the following.

\begin{Theorem}\label{T70}
Let $(N,g)$ be a connected, simply connected 3-dimensional Lorentzian manifold. Then the necessary and sufficient condition for $(N,g)$ to be symmetric and to have a Ricci operator of (degenerate) Segre type $\{21\}$, is the existence of a global pseudo-orthonormal frame field $\{e_1,e_2,e_3\}$, with $e_3$ timelike, a real constant $k$ and a smooth function $B$, satisfying {\em (\ref{newsyst21})-(\ref{Lie21})}.
\end{Theorem}

The following classification of parallel surfaces in a symmetric Lorentzian 3-space was also obtained by G. Calvaruso and J. Van der Veken in \cite{CVan09}.

\begin{Theorem}\label{T71}
Let $M$ be a parallel surface in a symmetric Lorentzian three-space
$(\tilde N,\tilde g)$ carrying a parallel null vector field,
described by {\em (\ref{newsyst21})-(\ref{Lie21})}. Then $M$ is an
integral surface of the distribution spanned by $\{e_2,e_3\}$, on
which $B$ is constant. Moreover, $M$ is always flat and $M$ is
totally geodesic if and only if $B=0$ on it. If $B$ is non-constant
on all integral surfaces of the distribution spanned by
$\{e_2,e_3\}$, then $(\tilde N, \tilde g)$ does not admit any
parallel surfaces.
\end{Theorem}

\section{Three natural extensions of parallel submanifolds}\label{sec28}

\subsection{Submanifolds with parallel mean curvature vector}\label{sec28.1}

One natural extension of the class of parallel submanifolds ($\bar\nabla h=0$) is the class of submanifolds with parallel mean curvature vector, i.e., $\bar\nabla ({\rm Tr}\, h)=0$ or equivalently $DH=0$. Trivially, both minimal submanifolds  and parallel submanifolds have parallel mean curvature vector automatically.
Further, a hypersurface of any Riemannian manifold has parallel mean curvature vector if and only if it has constant mean curvature.

Euclidean hypersurfaces with constant mean curvature are important since they are critical points of some natural functionals. In fact, a hypersurface of constant mean curvature in a Euclidean space is a solution to a variational problem. With respect to any volume-preserving variation of a domain $D$ in a Euclidean space the mean curvature of $M=\partial D$ is constant if and only if the volume of $M$ is critical, where $\partial D$ is the boundary of $D$. 

The condition of submanifolds to have parallel mean curvature vector in higher dimensional Euclidean spaces is very interesting as well since it is equivalent to a critical points of being variational problem; namely, their Gauss maps are harmonic maps (see \cite{RV}).

During the last 50 years, there are many research done on submanifolds with parallel mean curvature vector.  Among others, for submanifolds with parallel mean curvature vector in real space forms, see \cite{CL72,Ho73,chen73,Smyth, Re74,Yau74,CH1975,c1980}; for surfaces with parallel mean curvature vector in complex space forms, see \cite{Og,KZ,Hir06,Ken16,Hi,FP15,Fet12};  for surfaces with parallel mean curvature vector in indefinite space forms, see \cite{c6,cent,CG,Kyushu,cv3,FH10}; for surfaces  with parallel mean curvature vector in  homogeneous spaces or symmetric spaces, see \cite{MTV17,FT14}; for surfaces with parallel mean curvature vector in Sasakian space forms, see \cite{FR13}; and for surfaces with parallel mean curvature vector in reducible manifolds, see \cite{FR12,FR13.1,FR13.2,FR14}.
 For general references of submanifolds with parallel mean curvature vector, see \cite{Chen2010}.

\subsection{Higher order parallel submanifolds}\label{sec28.2}

Higher order parallel submanifolds, i.e., submanifolds that satisfy $\bar\nabla^k h=0$ for some positive integer $k$, were first studied by D. Del-Pezzo in \cite{Del1886} and then investigated by several authors after Del-Pezzo (see J. A. Schouten and D. J. Struik's 1938 book \cite{SS38} for details). This research topic was renewed in late 1980s by F. Dillen, V. Mirzoyan and \"U. Lumiste.  Since then, this interesting research topic has been studied by several differential geometers.

Among others, for  higher order parallel submanifolds  in real space forms, see  \cite{Di92,Lum87,Dillen90,Mirz78,Lum91,Mir98,KN12}; for higher order parallel surfaces in three-dimensional homogeneous spaces, see \cite{V07}; 
for higher order parallel surfaces in Bianchi-Cartan-Vranceanu spaces, see \cite{Ve08}; for higher order parallel surfaces in the Heisenberg group, see \cite{DV08}; and  for higher order parallel submanifolds of a complex space form, see \cite{DV90}. For some further results on higher order parallel submanifolds, see \"U Lumiste's 2000 survey article \cite{Lu00}.

\subsection{Semi-parallel submanifolds}\label{sec28.3}

The notion of semi-parallel submanifolds was introduced in 1985 by J. Deprez in \cite{De85}.
A submanifold $M$  of a Riemannian manifold $N$ is called {\it semi-parallel} if its second fundamental form $h$ satisfies 
$$ \tilde R(X,Y)\cdot h=(\bar \nabla_X \bar \nabla_Y-\bar \nabla_Y\bar \nabla_X-\bar \nabla_{[X,Y]})h=0,$$
where $\tilde R$ is the Riemann curvature tensor of $N$. Obviously, parallel submanifolds are semi-parallel. Hence, semi-parallel submanifolds are  natural extensions of parallel submanifolds as well.

In  \cite{De85}, J. Deprez applying the work of E. Backes \cite{Back83} on Euclidean Jordan triple systems to prove that totally geodesic surfaces are the only minimal semi-parallel surfaces in a Euclidean space.  Furthermore, he proved in \cite{De86} that every semi-parallel submanifolds of a Euclidean space is intrinsically a semi-symmetric Riemannian manifold. By a semi-symmetric Riemannian manifold $(M,g)$ we mean that the Riemann curvature curvature tensor of $(M,g)$  satisfies the condition $R\cdot R=0$, where the first tensor $R$ acts on the second one as a derivation.
In \cite{De85}, Deprez also classified semi-parallel surfaces in a Euclidean space.
Since then many articles were devoted to the study of semi-parallel submanifolds.

Among others, for semi-parallel submanifolds in real space forms, see \cite{De86,Dillen91,Mer91,Asp91,AM94,Lu96,Li01,Lum09,MMC12,BA15}; for semi-parallel submanifolds of  indefinite space forms, see \cite{Lumu97,Saf01}; for semi-parallel submanifolds in Kaehler manifolds, see \cite{NR98,DVV09,Kon11,OGM09}; for semi-parallel submanifolds in reducible spaces, see \cite{VV08}; for manifold with semi-parallel geodesic spheres or semi-parallel tubes, see \cite{BGV96}; for semi-parallel submanifolds in contact metric manifolds, see \cite{OGM14,AMU19}; and for semi-parallel submanifolds in other Riemannian manifolds, see \cite{Lu02,Lu03,Lu04}. For some further results on semi-parallel submanifolds, see \cite{Lu00}.

{\bf Acknowledgments}. The author thanks Joeri Van der Veken for his very careful reading of this article and for his many valuable suggestions to improve the presentation.

  \end{document}